\documentclass{article}
\usepackage{graphicx}
\usepackage{amssymb}
\usepackage[T1]{fontenc}
\usepackage{graphicx}
\usepackage{tikz-cd}
\usepackage{quiver}
\usepackage{mathtools}
\usepackage{stmaryrd}
\usepackage{mathpartir}
\usepackage{anyfontsize}
\usepackage{prftree}
\usepackage{color}
\usepackage{hyperref}
\usepackage[margin=1.4in]{geometry}
\usepackage[utf8]{inputenc}
\usepackage[all,2cell]{xy}
\usepackage{orcidlink}
\usepackage{graphicx}
\usepackage{amsmath, bm}
\usepackage{amsthm}
\usepackage{amsfonts}
\usepackage{amssymb}
\usepackage{bussproofs}
\usepackage{cleveref}
\usepackage{enumitem, ragged2e}
\definecolor{ceruleanblue}{rgb}{0.16, 0.32, 0.75}
\hypersetup{
    colorlinks=true,
    urlcolor=ceruleanblue,
    linkcolor=ceruleanblue,
    citecolor=ceruleanblue
}

\theoremstyle{plain}
\newtheorem{theorem}{Theorem}[section]
\newtheorem{cor}[theorem]{Corollary}

\newtheorem{proposition}[theorem]{Proposition}

\theoremstyle{definition}
\newtheorem{definition}[theorem]{Definition}
\newtheorem{remark}[theorem]{Remark}

\crefname{subsection}{Subsection}{Subsections}
\crefname{cor}{Corollary}{Corollaries}
\crefname{proposition}{Proposition}{Propositions}
\title{A biequivalence of path categories \\ and axiomatic Martin-L\"of type theories}
\author{Matteo Spadetto \orcidlink{0000-0002-6495-7405} \\
    \small University of Nantes, FR \\
    \small \href{mailto:matteo.spadetto.42@gmail.com}{matteo.spadetto.42@gmail.com}
}
\date{ }

\begin{document}

\newcommand{\judge}[2]{\boldsymbol{\lfloor}#1\boldsymbol{\rfloor}\;#2}
\newcommand{\judgectx}[1]{\boldsymbol{\lfloor}#1\boldsymbol{\rfloor}}
\newcommand{\type}{\textsc{Type}}
\newcommand{\ctx}{\textsc{Ctx}}
\newcommand{\liltriangle}{{\scalebox{0.5}{$\blacktriangledown$}}}
\newcommand{\lilbullet}{{\scalebox{0.6}{$\bullet$}}}
\newcommand{\id}{\mathsf{Id}}
\newcommand{\refl}{\mathsf{r}}
\newcommand{\reflex}{\mathsf{refl}}
\renewcommand{\j}{\mathsf{J}}
\newcommand{\h}{\mathsf{H}}
\newcommand{\sigmad}{\mathsf{\Sigma}}
\newcommand{\pair}{\mathsf{p}}
\newcommand{\pairing}{\mathsf{pair}}
\newcommand{\splitt}{\mathsf{split}}
\newcommand{\pid}{\mathsf{\Pi}}
\newcommand{\ev}{\mathsf{ev}}
\newcommand{\app}{\mathsf{app}}
\newcommand{\funext}{\mathsf{funext}}
\newcommand{\happly}{\mathsf{happly}}
\newcommand{\happ}{\mathsf{happ}}
\renewcommand{\t}{\mathsf{t}}
\newcommand{\grpd}{\mathbf{Grpd}}
\newcommand{\blank}{\raisebox{0.5ex}{\_}}
\newcommand{\abst}[1]{\mathsf{abst}_{#1}}
\newcommand{\s}{\mathsf{s}}
\newcommand{\successivo}{\mathsf{succ}}
\newcommand{\0}{\mathsf{0}}
\newcommand{\zero}{\mathsf{zero}}
\newcommand{\ind}{\mathsf{ind}}
\newcommand{\nat}{\mathbb{N}}

\title{A 2-categorical approach to the semantics\\of dependent type theory with computation axioms}

\maketitle

\begin{abstract}
Axiomatic type theory is a dependent type theory without computation rules. The term equality judgements that usually characterise these rules are replaced by computation \textit{axioms}, i.e., additional term judgements that are typed by identity types. This paper is devoted to providing an effective description of its semantics, from a higher categorical perspective: given the challenge of encoding intensional type formers into 1-dimensional categorical terms and properties, a challenge that persists even for axiomatic type formers, we adopt Richard Garner’s approach in the 2-dimensional study of dependent types. We prove that the type formers of axiomatic theories can be encoded into natural 2-dimensional category theoretic data, obtaining a presentation of the semantics of axiomatic type theory via 2-categorical models called \textit{display map 2-categories}. In the axiomatic case, the 2-categorical requirements identified by Garner for interpreting intensional type formers are relaxed. Therefore, we obtain a presentation of the semantics of the axiomatic theory that generalises Garner’s one for the intensional case.
Our main result states that the interpretation of axiomatic theories within display map 2-categories is well-defined and enjoys the soundness property. We use this fact to provide a semantic proof that the computation rule of intensional identity types is not admissible in axiomatic type theory. This is achieved via a revisitation of Hofmann and Streicher's \textit{groupoid model} that believes axiomatic identity types but does not believe intensional ones.
\end{abstract}

\text{ }

\medskip

\medskip

\setcounter{tocdepth}{1}
\tableofcontents

\section{Introduction and motivation}\label{section:Introduction and motivation} This paper presents an analysis of the semantics of dependent types from a categorical perspective, with the aim of applying it to the study of axiomatic type theories.

\subsection{The problem of the semantics of intensional type formers}\label{subsection:The problem of the semantics of intensional type formers} In formulating the semantics of theories of dependent types, we typically encounter two distinct, but related, approaches: \textit{syntactic} and \textit{categorical}. The syntactic approach directly mirrors the structure of the syntax of the theory, defining how judgements in the conclusions of inference rules are interpreted based on the interpretations assigned to the premises---in other words, we equip a model with a choice function for every inference rule: a function that chooses, for the interpretation of every instance of the premises, an interpretation of the corresponding instance of the conclusion. On the other hand, the categorical approach adds structure to the model, structure that allows to automatically recover choice functions analogous to those provided in the syntactic approach.

For example, given a display map category $(\mathbf{C},\mathcal{D})$---we will recall this notion later---consider the case of extensional $=$-types. As we mentioned, to model a theory featuring this type former we may just replicate its syntax within the “language” of $(\mathbf{C},\mathcal{D})$, according to which a \textit{display map} is the interpretation of a type judgement $A:\type$ and its \textit{sections} interpret the corresponding term judgements $t:A$. Hence, a display map interpreting a given type $A : \type$ is required to be endowed with the choice of another display map interpreting the identity type $\judge{x, x' : A}{x = x' : \type}$; in other words, the display map category $(\mathbf{C},\mathcal{D})$ is equipped with a choice function validating the \textit{formation rule} of extensional $=$-types. Analogously, additional choice functions will validate the other rules of this type former. But alternatively, as explained in \cite{MR1674451}, we can simply require that the diagonal arrow $\Gamma.A\to \Gamma.A.A[P_A]$ is isomorphic to some display map whenever $P_A$ is a display map $\Gamma.A\to \Gamma$ in $\mathcal{D}$, because this categorical requirement on $(\mathbf{C},\mathcal{D})$ is equivalent to equipping $(\mathbf{C},\mathcal{D})$ with the necessary choice functions to validate the extensional $=$-type former.
To provide the reader with another example, if we add $\Sigma$-types to the theory, then the categorical condition characterising its semantics is the following: it stipulates that display maps must be closed under composition up to isomorphism. This categorical characterisation of $\Sigma$-types is the basis of the analysis contained in \cite{MR0727078} and \cite{hofmannlccc}.

These are examples of how extensional type formers admit a concise categorical characterisation of their models: it consists of \textit{closure properties on the class of their display maps}. However, in case we drop extensional identity types and we ask ourselves how to handle $=$- and $\Sigma$-types, things become more complicated, and such a categorical characterisation of their semantics is harder to find.
To simplify the problem and recover a categorical characterisation of the semantics of intensional type formers, Garner proposes a 2-categorical formulation in \cite{MR2525957}. In this scenario, a synthetic and conceptually simple characterisation of identities, dependent sums, and other type constructors \textit{exists}. This stems from the wider range of potential closure properties---whether strong or weak---to require on the class $\mathcal{D}$, that we can access: properties that enable us to distinguish between varying strengths of type formers. Notably, if dimension 2 were omitted, these properties would collapse into the ones we have just discussed, recovering the extensional version of our type constructors.

In this paper, we adopt Garner's perspective to study the semantics of \textit{axiomatic type theory} from a categorical point of view, specifically through a 2-categorical one.

\subsection{Axiomatic dependent type theory}\label{subsection:Axiomatic dependent type theory}

In recent years, there has been a growing interest in various weakenings for theories of dependent
types, particularly those weakenings with respect to the strength of the computation rules of the type formers. A type former is meant to encode a piece of logic---e.g. the rules of existential quantification for a notion of $\Sigma$-type former---into the \textit{propositions-as-types} formalism of dependent type theories. It is the ensemble of a \textit{formation} rule, an \textit{introduction} rule, and an \textit{elimination} rule, together with one or more \textit{computation} rules that explain the interaction between terms obtained by introducing and eliminating, and that allow one to simplify or re-phrase proofs. In both extensional and intensional type theory (ETT and ITT respectively), computation rules consist of term equality judgements, of the form $t\equiv t'$, where $t:A$ and $t':A$ for some type judgement $A:\type$. However, term equality judgements are not the only available way to state that “two terms of a given type are equal”. If a theory has some flavour of the \textit{identity type} former, one may judge a specific term judgement $p:t=t'$, typed by the \textit{identity type} $t = t' : \type$, to state that “there is a proof that the two terms $t:A$ and $t':A$ are equal”. This leads us to a \textit{weaker} equality judgement, expressed as $p : t = t'$. We say that $t\equiv t'$ is a \textit{judgemental} equality, because it consists of a term equality judgement, and that $p:t=t'$ is a \textit{propositional} equality, because the judgement $t=t':\type$---unlike the judgement $t\equiv t'$---is a type, hence an actual statement or proposition, and $p$ is a proof, or a witness, that “the statement $t=t'$ holds”.

Judgemental equalities are stronger than propositional equalities because of $=$-introduction rule. Propositional equalities are as strong as judgemental equalities if $=$-types are extensional. On the other hand, propositional equalities are actually weaker than judgemental equalities if $=$-types are intensional: the groupoid model \cite{HofmannStreicher-Thegroupoidmodelref,MR1686862,MR4442637} is a model of intensional identity types that does not believe in the uniqueness of identity proof, whereas every extensional model will do so. Hence the groupoid model does not have extensional identity types and therefore if will not necessarily believe that $t\equiv t'$ whenever the type $t=t'$ is inhabited. As we are going to see, propositional equalities become even weaker if $=$-types are just \textit{axiomatic}.

Axiomatic $=$-types are “very intensional” $=$-types. They consist of the usual rules of intensional $=$-types, except for the judgemental equality of its computation rule, \textit{which is replaced by a propositional equality}: whenever we are given judgements: \begin{center} $\judge{x,x':A;\;p:x=x'}{C(x,x',p):\type}$ \;\;\;\; and \;\;\;\; $\judge{x:A}{c(x):C(x,x,r(x))}$ \end{center} and hence, by elimination, a term judgement $\judge{x,x':A;\;p:x=x'}{\j(c,x,x',p):C(x,x',p)}$, then, in place of asking that the term equality judgement $\judge{x:A}{\j(c,x,x,r(x))\equiv c(x)}$ holds, we only ask that an additional term judgement of the form: $$\judge{x:A}{\h(c,x):\j(c,x,x,r(x))=c(x)}$$ holds. This term judgement replacing the computation rule is called \textit{computation axiom}: it requires that the statement $\j(c,x,x,r(x))=c(x)$ holds. This is depicted in detail in Fig. \ref{propositional id}.

\begin{figure*}[t]
\begin{center}
\fbox{$\begin{alignedat}{2}
\text{Form Rule\,}&\inferrule{\;\;\quad\quad\quad{A : \type}\quad\quad\quad\;\;}{\judge{x,x' : A}{x=x' : \type}}&\text{Elim Rule\,}&\inferrule{A : \type\\\\\judge{x, x';\; p : x=x'}{C(x,x',p) : \type}\\\\\judge{x}{c(x) : C(x,x,\refl(x))}}{\judge{x, x';\; p}{\j(c,x,x',p) : C(x,x',p)}}\\
\\
\text{Intro Rule\,}&\inferrule{\;\;\quad\quad\quad{A : \type}\quad\quad\quad\;\;}{\judge{x : A}{\refl(x) : x=x}}&\;\;\;\text{Comp Axiom\,}&\inferrule{A : \type\\\\\judge{x, x';\; p : x=x'}{C(x,x',p) : \type}\\\\ \judge{x}{c(x) : C(x,x,\refl(x))}}{\judge{x}{\h(c,x):\j(c,x,x,\refl(x))= c(x)}}\\
\end{alignedat}$}
\caption{Axiomatic $=$-types}\label{propositional id}
\end{center}
\end{figure*}

In general, when a dependent type theory has a type former whose computation rules are replaced by computation axioms, one says that the type former is in \textit{axiomatic} form. This paper is about semantics of \textit{axiomatic type theory} (ATT), a weakening of ITT with computation axioms. We will be focusing on the following type formers: $=$-types, $\Sigma$-types, $\Pi$-types, function extensionality, $0$-types, $1$-types, $2$-types, and $\nat$-types.

There are several motivations and advantages to work on ATT:

\noindent$\succ$ Axiomatic type theory is \textit{objective} \cite{MR4352353,2021arXiv210200905V} in the way operations, such as the sum of natural numbers or the composition of identity proofs, can be defined. In ITT, depending on how the elimination rule of these type formers is applied, will satisfy either $n+0 \equiv n$ and $0 + n = n$ or alternatively $n+0 = n$ and $0 + n \equiv n$, if $n$ is a natural number---and similarly for composition of identity proofs. In ATT, such an ambiguity is removed from the outset: every notion of natural sum and every notion of identity proof composition will satisfy all these equalities just in propositional form.

At the same time, while it only has one notion of equality like ETT, axiomatic type theory does not deduce the uniqueness of identity proofs, and it is therefore compatible with homotopy type theory.

\noindent$\succ$ Axiomatic type formers are \textit{homotopy invariant}, as explained in \cite{MR3050430,MR3614859}. Whenever a type is homotopy equivalent to one derived from the formation rule of a given type former, it adheres to all the rules associated with that type former.

\noindent$\succ$ As shown in \cite{2021arXiv210200905V}, by replacing computation rules with computation axioms, the \textit{type checking}, i.e. the decidability of the derivability of a term judgement $t:A$, holds and can be done in quadratic time. This improves the non-elementary one of ordinary ITT.

\noindent$\succ$ Compared to ITT, axiomatic type theory enjoy a \textit{broader concept of semantics}, in the sense that it has more models---every model of ITT is also a model of ATT but, as we show in the last section, there are models of the latter that do not validate computation rules---making it simpler---in practice there are less constraints regarding judgemental equalities to be satisfied and checked---to construct concrete models, for example, to obtain independence results.

Nevertheless, ATT is as expressive as ITT and encodes the whole intuitionistic predicate logic.

\noindent$\succ$ Accordingly, several \textit{conservativity results} of intensional and extensional type theories over the axiomatic one hold: for example, ETT extended with function extensionality and uniqueness of identity proof rules is conservative over ATT \cite{Winterhalter2020,conservativity}. Building on Hofmann's work \cite{hofmannconservativity}, this topic has been studied extensively by Bocquet \cite{bocquet}, by Boulier and Winterhalter \cite{boulierwinterhalter,Winterhalter2020}, and by the author \cite{conservativity}, showing how the axiomatic theory does not lose much deductive strength with respect to the intensional and extensional ones. Therefore, one may to some extent interchangeably use axiomatic type formers in order to study intensional and extensional ones, and vice versa.

\subsection{Related work}\label{subsection:Related work} Hints of the emergence of an interest in propositional computation rules, or computation axioms, are scattered everywhere in the homotopy type theoretic literature, starting from the work of Awodey, Gambino, Sojakova \cite{MR3050430,MR3614859}, that presented an initial study of an axiomatic type former, namely axiomatic well-founded tree types, or axiomatic $W$-types. Coquand and others \cite{MR3124820, MR3781068} conducted initial analyses related to axiomatic $=$-types. Another recent paper introduced a univalent model of ITT where the computation axiom for identity types is validated, although its judgemental version is not \cite{MR3281415}. Following this, the type constructor has been thoroughly examined by van den Berg and Moerdijk \cite{MR3828037, MR3795638}, who introduced and explored a semantic concept for dependent type theories with propositional identity types, using the notion of a \textit{path category} and showing its link with axiomatic $=$-types. Otten and the author \cite{daniel} explain how this notion of semantics is actually an instance of the usual one via display map categories (formulated as full comprehension categories). Among other works on axiomatic type formers, we note the previously mentioned studies on conservativity over such type theories by Bocquet \cite{bocquet}, by Boulier and Winterhalter \cite{boulierwinterhalter}, and by the author \cite{conservativity}, Bocquet's work on the coherence property of the class of models of axiomatic type formers \cite{MR4481908}, and Vidmar's work on extending natural models \cite{MR3742564} to type theories with propositional expansion rules, or expansion axioms \cite{vidmar2018}.

\subsection{Contributions}\label{subsection:Contributions}

In this paper, we adopt Garner's perspective \cite{MR2525957} to study the semantics of ATT from a categorical point of view, specifically from a 2-categorical one. In detail, we prove that a \textit{display map 2-category} with some categorical structure---obtained reformulating the syntax of the axiomatic type formers in 2-categorical terms---is sufficient to recover the semantic counterparts of these type constructors as in the \textit{syntactic approach}, particularly inducing an ordinary (split) display map category---i.e. an actual model of the structural rules of dependent type theory---that models the theory of dependent types in axiomatic form---see Section \ref{section:Categorical formulation of the semantics of axiomatic dependent type theory} and Theorem \ref{from a 2-dimensional model to a 1-dimensional one}. This provides an answer to the problem of formulating models of ATT in categorical terms, i.e. within the \textit{categorical approach}, in a way that includes and generalises the class of intensional models identified by Garner: display maps are not necessarily normal isofibrations, as in the intensional case, but merely cloven isofibrations and the axiomatic form of $\Sigma$-types does not require that display maps be closed under composition up to injective equivalence, as for intensional models, but merely up to homotopy equivalence.

In broad terms, this work falls within the field of categorical semantics of dependent type theories, and aims to contribute to the central objective of developing a general and feasible notion of semantics for variants and generalisations of dependent type theory.

\subsection{Outline}\label{subsection:Outline} In Section \ref{section:Recap on display map categories} we recall the notion of a (split) display map category and briefly explain how such a structure provides a model of the structural rules of dependent type theory. In Section \ref{section:Syntactic approach to the semantics of axiomatic dependent type theory} we specialise this notion to model ATT, using the syntactic approach. In Section \ref{section:Categorical formulation of the semantics of axiomatic dependent type theory}, we formulate the notion of a display map 2-category as a structure that encodes the structural and logical aspects of ATT as 2-categorical data on display maps. Accordingly, we show that fulfilling such data allows a display map 2-category to induce a model of ATT as described in Section \ref{section:Syntactic approach to the semantics of axiomatic dependent type theory}. In Section \ref{sec:Revisiting the groupoid model}, we use this result to identify a model of ATT---based on the groupoid model---that is not a model of ITT. In Section \ref{sec:conclusion and future work}, we discuss the position of this work within contemporary research on the categorical semantics of generalised dependent type theory.

\section{Basics of display map categories}\label{section:Recap on display map categories}

In this section, we recall the notion of display map categories \cite{taylor1999practical, MR1674451,mossvonglehn2018} one of several equivalent categorical structures that provide a sound and complete semantics for dependent type theories. We specialise it to provide such a notion of semantics for ATT with $=$-types and $\Sigma$-types. We leave the semantics of the other axiomatic type formers for an extended version.

\begin{definition}\label{display map category} A \textit{(split) display map category} $(\mathbf{C},\mathcal{D})$ is a category $\mathbf{C}$ with a chosen terminal object $\epsilon$, together with a class $\mathcal{D}$ of arrows of $\mathbf{C}$, called \textit{display maps} and denoted as $\Gamma.A\to \Gamma$ and labeled as $P_A$, such that:\begin{itemize}
    \item[$\succ$] for every display map $\Gamma.A\to \Gamma$ and every arrow $f: \Delta \to \Gamma$, there is a choice of a display map $\Delta.A[f]\to\Delta$ such that the square: 
    \[\begin{tikzcd}[column sep=tiny,row sep=small]
	{\Delta.A[f]} & {\Gamma.A} \\
	\Delta & \Gamma
	\arrow[from=1-1, to=1-2]
	\arrow[from=1-1, to=2-1]
	\arrow["\lrcorner"{anchor=center, pos=0.125}, draw=none, from=1-1, to=2-2]
	\arrow[from=1-2, to=2-2]
	\arrow["f"{description}, from=2-1, to=2-2]
    \end{tikzcd}\] is a pullback; as the former display map is labeled as $P_A$, the latter will be labeled as $P_{A[f]}$; the arrow $\Delta.A[f]\to \Gamma.A$ will be denoted as $f.A$, or as $f^\lilbullet$ in absence of ambiguity; we will adopt the notation $A^\liltriangle$ in place of $A[P_A]$;
    \item[$\succ$] the equalities: $$\begin{alignedat}{2}
        P_{A[1_A]}&=P_A& (1_\Gamma).A&=1_{\Gamma.A}\\
        P_{A[fg]}&=P_{A[f][g]}&\quad\quad (f.A)(g.A[f])&=(fg).A
    \end{alignedat}$$ hold for every choice of composable arrows $f$ and $g$ and every display map $\Gamma.A\to \Gamma$, where $\Gamma$ is the target of $f$.
\end{itemize}
\end{definition} In this paper, the term “display map category” specifically refers to a “\textit{split} display map category”. That is, we always assume that the display map categories we discuss satisfy the second of the two conditions in Definition \ref{display map category}. Under this assumption, the notion of a display map category is entirely equivalent to other notions commonly used as models of dependent type theories, such as categories with attributes \cite{phdcartmell,moggi_1991, KAPULKIN2021106563}, categories with families \cite{clairambault_dybjer_2014,dybjer1996internal,phdhofmann,Hofmann1997}, and full split comprehension categories \cite{MR1201808,MR3372323}.

Such a display map category $(\mathbf{C},\mathcal{D})$ constitutes a sound model of the structural part of dependent type theory: as we mentioned in Section \ref{section:Introduction and motivation}, a context is encoded as an object $\Gamma$ of $\mathbf{C}$, a type in that context is interpreted as a display map $P_A$ of codomain $\Gamma$, and a term of that type as a section of $P_A$. Then we can observe e.g. that the fundamental structural rules: $$\inferrule{\text{ }}{\judgectx{\blank}}\;\;\;\inferrule{\judge{\gamma:\Gamma}{A(\gamma):\type}}{\judgectx{\gamma:\Gamma,x:A(\gamma)}}\;\;\;\inferrule{\judge{\gamma:\Gamma}{A(\gamma):\type}}{\judge{\gamma:\Gamma,x:A(\gamma)}{x:A(\gamma)}}$$ are validated: the terminal object $\epsilon$ interprets the \textit{empty context} $\blank$ and, if the object $\Gamma$ is the interpretation of the context $\gamma$ and the display map $P_A: \Gamma.A\to\Gamma$ is the interpretation of the type-in-context $A$, then the object $\Gamma.A$ interprets the \textit{extended context} $\gamma,x$ and the unique section $\delta_A:\Gamma.A\to \Gamma.A.A^\liltriangle$ of $P_{A^\liltriangle}$ such that the diagram: 
\[\begin{tikzcd}[sep=small]
	{\Gamma.A} \\
	{\Gamma.A.A^\liltriangle} & {\Gamma.A} \\
	{\Gamma.A} & \Gamma
	\arrow[from=1-1, to=2-1]
	\arrow[Rightarrow, no head, from=1-1, to=2-2]
	\arrow[from=2-1, to=2-2]
	\arrow[from=2-1, to=3-1]
	\arrow["\lrcorner"{anchor=center, pos=0.125}, draw=none, from=2-1, to=3-2]
	\arrow[from=2-2, to=3-2]
	\arrow["{P_A}"{description}, from=3-1, to=3-2]
\end{tikzcd}\] commutes interprets the \textit{variable term-in-context} $x$. Additionally, substitution---and hence weakening---is interpreted by using the cleavage of $(\mathbf{C},\mathcal{D})$: if the section $a:\Gamma\to\Gamma.A$ of $P_A$ is the interpretation of a term judgement $\judge{\gamma:\Gamma}{a(\gamma):A(\gamma)}$ and the arrow $f:\Delta\to\Gamma$ is the interpretation of a substitution $\judge{\delta:\Delta}{f(\delta):\Gamma}$, then the judgements: $$\judge{\delta:\Delta}{A(f(\delta)):\type}\text{\quad and \quad}\judge{\delta:\Delta}{a(f(\delta)):A(f(\delta))}$$ will be interpreted by the display map $P_{A[f]}$ and by its unique section $a[f]:\Delta\to\Delta.A[f]$ such that the square: 
\[\begin{tikzcd}[sep=scriptsize]
	\Delta & \Gamma \\
	{\Delta.A[f]} & {\Gamma.A}
	\arrow["f"{description}, from=1-1, to=1-2]
	\arrow["{a[f]}"', from=1-1, to=2-1]
	\arrow["a", from=1-2, to=2-2]
	\arrow["{f^\lilbullet}", from=2-1, to=2-2]
\end{tikzcd}\] commutes, respectively.

\subsection{Encoding parallel terms into a substitution}\label{subsection:Encoding parallel terms into a substitution} If $a$ and $b$ are sections of a display map $P_A:\Gamma.A\to\Gamma$, we denote as $a;b$ the unique arrow $\Gamma\to\Gamma.A.A^\liltriangle$ such that the diagram: 
\[\begin{tikzcd}[sep=scriptsize]
	\Gamma \\
	& {\Gamma.A.A^\liltriangle} & {\Gamma.A} \\
	& {\Gamma.A} & \Gamma
	\arrow[from=1-1, to=2-2]
	\arrow["b"{description}, curve={height=-18pt}, from=1-1, to=2-3]
	\arrow["a"{description}, curve={height=18pt}, from=1-1, to=3-2]
	\arrow["{P_A^\lilbullet}", from=2-2, to=2-3]
	\arrow["{P_{A^\liltriangle}}"', from=2-2, to=3-2]
	\arrow["\lrcorner"{anchor=center, pos=0.125}, draw=none, from=2-2, to=3-3]
	\arrow[from=2-3, to=3-3]
	\arrow[from=3-2, to=3-3]
\end{tikzcd}\] commutes. If $a$ and $b$ interpret given terms $a(\gamma):A(\gamma)$ and $b(\gamma):A(\gamma)$ in context $\gamma:\Gamma$, then the substitution $\judge{\gamma:\Gamma}{\gamma,a(\gamma),b(\gamma):\Gamma.A.A^\liltriangle}$ is interpreted by the arrow $a;b$.

We recall that the square: \[\begin{tikzcd}[column sep=small,row sep=scriptsize]
	\Delta & \Gamma \\
	{\Delta.A[f].A[f]^{\liltriangle}} & {\Gamma.A.A^\liltriangle}
	\arrow["f", from=1-1, to=1-2]
	\arrow["{a[f];b[f]}"', from=1-1, to=2-1]
	\arrow["{a;b}", from=1-2, to=2-2]
	\arrow["{f^{\lilbullet\lilbullet}}", from=2-1, to=2-2]
\end{tikzcd}\] commutes for every arrow $f:\Delta\to\Gamma$. This fact is crucial to prove that the stability conditions relative to $\h_c$ and $\sigma_c$ type-check, in Definitions \ref{def:Semantics of axiomatic =-types} and \ref{def:Semantics of axiomatic Sigma-types} respectively. Hence, we recall a proof of this fact in Appendix \ref{appendixI}.

In the next section we show how to use display map categories to model axiomatic type theory via the syntactic approach.

\section{Syntactic approach to the semantics of axiomatic type theory}\label{section:Syntactic approach to the semantics of axiomatic dependent type theory}

In this section we specialise the notion of display map categories to make them into models of axiomatic type theory. This consists in endowing a display map category with additional structure constituting the semantic counterpart of the type former. As we mentioned in \Cref{section:Introduction and motivation}, this additional structure may be formulated in alignment with the syntax, by means of choice functions, each associated to a given rule of the theory, assigning to the interpretation of the premises of that rule an interpretation of its consequence. These choice functions constitute an encoding of type formers into the given display map categories. If a display map category $(\mathbf{C},\mathcal{D})$ is endowed with such a function for every logical rule in the theory, then the interpretation in $(\mathbf{C},\mathcal{D})$ of every type judgement and every term judgement of the theory is defined and respects all the type equality judgements and all the term equality judgements of the theory. In other words, the interpretation of the theory in $(\mathbf{C},\mathcal{D})$ is defined and is \textit{sound}.

Referring to the rules of Fig. \ref{propositional id} for axiomatic $=$-types, we give the following:

\begin{definition}[Semantics of axiomatic $=$-types---syntactic formulation]\label{def:Semantics of axiomatic =-types}
    Let $(\mathbf{C},\mathcal{D})$ be a display map category.
    
    Let us assume that, for every object $\Gamma$ and every display map $P_A$ of codomain $\Gamma$, there is a choice of:
    \begin{itemize}
        \item[$\succ$] (\textit{Form Rule}) a display map $\Gamma.A.A^\liltriangle.\id_A\to\Gamma.A.A^\liltriangle$;
        \item[$\succ$] (\textit{Intro Rule}) a section: $$\reflex_A:\Gamma.A\to \Gamma.A.\id_A[\delta_A]$$ of the display map $\Gamma.A.\id_A[\delta_A]\to\Gamma.A$;
    \end{itemize} and that, for every object $\Gamma$, every display map $P_A$ of codomain $\Gamma$, every display map $P_C$ of codomain $\Gamma.A.A^\liltriangle.\id_A$, and every section $c:\Gamma.A\to\Gamma.A.C[\refl_A]$ of $P_{C[\refl_A]}$---where $\refl_A$ is the composition: 
\[\begin{tikzcd}[column sep=2.25em,row sep=tiny]
	{\Gamma.A} && {\Gamma.A.\id_A[\delta_A]} && {\Gamma.A.A^\liltriangle.\id_A}
	\arrow["{\reflex_A}"{description}, from=1-1, to=1-3]
	\arrow["{\delta_A^\lilbullet}"{description}, from=1-3, to=1-5]
\end{tikzcd}\]---there is a choice of:
    \begin{itemize}
        \item[$\succ$] (\textit{Elim Rule}) a section: $$\j_c : \Gamma.A.A^\liltriangle.\id_A \to\Gamma.A.A^\liltriangle.\id_A.C$$ of the display map $\Gamma.A.A^\liltriangle.\id_A.C\to\Gamma.A.A^\liltriangle.\id_A$;
        \item[$\succ$] (\textit{Comp Axiom}) a section: $$\h_c:\Gamma.A\to\Gamma.A.\id_{C[\refl_A]}[\j_c[\refl_A];c]$$ of the display map $\Gamma.A.\id_{C[\refl_A]}[\j_c[\refl_A];c]\to\Gamma.A$, where the arrow: $$\j_c[\refl_A];c: \Gamma.A\to\Gamma.A.C[\refl_A].C[\refl_A]^\liltriangle$$ is built as explained in Subsection \ref{subsection:Encoding parallel terms into a substitution}.
    \end{itemize} Moreover, let us assume\footnotetext[1]{We refer the reader to Appendix \ref{appendixI} for additional details on the type-checking of the stability conditions.} that the following \textit{stability conditions}: $$\begin{alignedat}{2}
        \id_A[f^{\lilbullet\lilbullet}]&=\id_{A[f]}&\quad\quad\j_c[f^{\lilbullet\lilbullet\lilbullet}]&=\j_{c[f^\lilbullet]}\\
        \reflex_A[f^\lilbullet]&=\reflex_{A[f]}&\h_c[f^\lilbullet]&=\h_{c[f^\lilbullet]}\\
    \end{alignedat}$$ hold\footnotemark[1] for every arrow $f:\Delta\to\Gamma$.

    \noindent Then we say that $(\mathbf{C},\mathcal{D})$ is \textit{endowed with axiomatic $=$-types}.
\end{definition}

\begin{figure*}[t]
\begin{center}
\fbox{$\begin{alignedat}{2}
\text{Form Rule\,}&\inferrule{{A : \type}\\\\\judge{x:A}{B(x):\type}}{{\Sigma_{x:A}B(x) : \type}}&\quad\text{Elim Rule\,}&\inferrule{\quad{A : \type}\quad\quad\judge{x:A}{B(x):\type}\quad\\\\ \judge{u:\Sigma_{x:A}B(x)}{C(u):\type} \\\\ \judge{x:A;\;y:B(x)}{c(x,y):C(\langle x,y\rangle)}}{\judge{u:\Sigma_{x:A}B(x)}{\splitt(c,u):C(u)}}\\
\\
\text{Intro Rule\,}&\inferrule{{A : \type}\\\\\judge{x:A}{B(x):\type}}{\judgectx{x:A;\;y:B(x)}\\\\{\langle x,y\rangle:\Sigma_{x:A}B(x)}}&\quad\text{Comp Axiom\,}&\inferrule{\quad{A : \type}\quad\quad\judge{x:A}{B(x):\type}\quad\\\\ \judge{u:\Sigma_{x:A}B(x)}{C(u):\type} \\\\ \judge{x:A;\;y:B(x)}{c(x,y):C(\langle x,y\rangle)}}{\judgectx{x:A;\;y:B(x)}\\\\{\sigma(c, x,y):\splitt(c,\langle x,y\rangle)= c(x,y)}}\\
\end{alignedat}$}
\caption{Axiomatic $\Sigma$-types}\label{propositional sigma}
\end{center}
\end{figure*}

Analogously, referring to the rules of Fig. \ref{propositional sigma} for axiomatic $\Sigma$-types, we give the following:

\begin{definition}[Semantics of axiomatic $\Sigma$-types---syntactic formulation]\label{def:Semantics of axiomatic Sigma-types}
    Let $(\mathbf{C},\mathcal{D})$ be a display map category endowed with axiomatic $=$-types.
    
    Let us assume that, for every object $\Gamma$, every display map $P_A$ of codomain $\Gamma$, and every display map $P_B$ of codomain $\Gamma.A$, there is a choice of:
    \begin{itemize}
        \item[$\succ$] (\textit{Form Rule}) a display map $\Gamma.\sigmad_A^B\to\Gamma$;
        \item[$\succ$] (\textit{Intro Rule}) a section: $$\pairing_A^B:\Gamma.A.B\to \Gamma.A.B.\sigmad_A^B[P_AP_B]$$ of the display map $\Gamma.A.B.\sigmad_A^B[P_AP_B]\to\Gamma.A.B$;
    \end{itemize} and that, for every object $\Gamma$, every display map $P_A$ of codomain $\Gamma$, every display map $P_B$ of codomain $\Gamma.A$, every display map $P_C$ of codomain $\Gamma.\sigmad_A$, and every section $c:\Gamma.A.B\to\Gamma.A.B.C[\pair_A^B]$ of $P_{C[\pair_A^B]}$---where $\pair_A^B$ is the composition: 
\[\begin{tikzcd}[column sep=small,row sep=tiny]
	{\Gamma.A.B} &&&& {\Gamma.A.B.\sigmad_A^B[P_AP_B]} &&&& {\Gamma.\sigmad_A^B}
	\arrow["{\pairing_A^B}"{description}, from=1-1, to=1-5]
	\arrow["{(P_AP_B)^\lilbullet}"{description}, from=1-5, to=1-9]
\end{tikzcd}\]---there is a choice of:
    \begin{itemize}
        \item[$\succ$] (\textit{Elim Rule}) a section: $$\splitt_c : \Gamma.\sigmad_A^B\to\Gamma.\sigmad_A^B.C$$ of the display map $\Gamma.\sigmad_A^B.C\to\Gamma.\sigmad_A^B$;
        \item[$\succ$] (\textit{Comp Axiom}) a section: $$\sigma_c:\Gamma.A.B\to\Gamma.A.B.\id_{C[\pair_A^B]}[\splitt_c[\pair_A^B];c]$$ of the display map $\Gamma.A.B.\id_{C[\pair_A^B]}[\splitt_c[\pair_A^B];c]\to\Gamma.A.B$, where the arrow: $$\splitt_c[\pair_A^B];c: \Gamma.A.B\to\Gamma.A.B.C[\pair_A^B].C[\pair_A^B]^\liltriangle$$ is built as explained in Subsection \ref{subsection:Encoding parallel terms into a substitution}.
    \end{itemize} Moreover, let us assume that the following \textit{stability conditions}: $$\begin{alignedat}{2}
        \sigmad_A^B[f]&=\sigmad_{A[f]}^{B[f^\lilbullet]}&\quad\quad\splitt_c[f^\lilbullet]&=\splitt_{c[f^{\lilbullet\lilbullet}]}\\
        \pairing_A^B[f^{\lilbullet\lilbullet}]&=\pairing_{A[f]}^{B[f^\lilbullet]}&\sigma_c[f^{\lilbullet\lilbullet}]&=\sigma_{c[f^{\lilbullet\lilbullet}]}\\
    \end{alignedat}$$ hold\footnotemark[1] for every arrow $f:\Delta\to\Gamma$.

    \noindent Then we say that $(\mathbf{C},\mathcal{D})$ is \textit{endowed with axiomatic $\Sigma$-types}.
\end{definition}

\begin{figure*}[t]
\begin{center}
\fbox{$\begin{alignedat}{2}
\text{Form Rule\,}&\inferrule{{A : \type}\\\\\quad\judge{x:A}{B(x):\type}\quad}{{\Pi_{x:A}B(x) : \type}}&\quad\text{Intro Rule\,}&\inferrule{{A : \type}\quad\quad\judge{x:A}{B(x):\type}\\\\\judge{x:A}{\upsilon(x):B(x)}}{{\lambda x.\upsilon(x):\Pi_{x:A}B(x)}}\\
\\
\text{Elim Rule\,}&\inferrule{{A : \type}\\\\\quad\judge{x:A}{B(x):\type}\quad}{\judgectx{z:\Pi_{x:A}B(x);\;x:A}\\\\{\ev(z,x):B(x)}}&\quad\text{Comp Axiom\,}&\inferrule{{A : \type}\quad\quad\judge{x:A}{B(x):\type}\\\\\judge{x:A}{\upsilon(x):B(x)}}{\judgectx{x:A}\\\\{\beta(\upsilon,x):\ev(\lambda x.\upsilon(x),x)=\upsilon(x)}}\\
\\
\\
\\
\\
&&\quad\text{Intro Rule\,}&\inferrule{{A : \type}\quad\quad\judge{x:A}{B(x):\type}}{\judgectx{z,z';\;q:\Pi_{x:A}\ev(z,x)=\ev(z',x)}\\\\{\funext(z,z',q):z=z'}}\\
\\
\text{Exp Axiom\,}&\inferrule{{A : \type}\\\\\quad\judge{x:A}{B(x):\type}\quad}{\boldsymbol{\lfloor}z,z';\;p:z=z'\boldsymbol{\rfloor}\\\\\eta^\Pi(z,z',p):p=\\\\\funext(z,z',\lambda x.\ev(p,x))}&\quad\text{Comp Axiom\,}&\inferrule{{A : \type}\quad\quad\judge{x:A}{B(x):\type}\\\\\text{ }}{\judgectx{z,z';\;q:\Pi_{x:A}\ev(z,x)=\ev(z',x)}\\\\\beta^\Pi(z,z',q):\\\\\lambda x.\ev(\funext(z,z',q),x)=q}
\end{alignedat}$}
\caption{Axiomatic $\Pi$-types \& axiomatic function extensionality}\label{propositional pi}
\end{center}
\end{figure*}

Now, referring to the rules of Fig. \ref{propositional pi} for axiomatic $\Pi$-types \& axiomatic function extensionality, we give the following:

\begin{definition}[Semantics of axiomatic $\Pi$-types \& axiomatic function extensionality---syntactic formulation]\label{def:Semantics of axiomatic Pi-types}
    Let $(\mathbf{C},\mathcal{D})$ be a display map category endowed with axiomatic $=$-types.

\medskip
    
    Let us assume that, for every object $\Gamma$, every display map $P_A$ of codomain $\Gamma$, and every display map $P_B$ of codomain $\Gamma.A$, there is a choice of:
    \begin{itemize}
        \item[$\succ$] (\textit{Form Rule}) a display map $\Gamma.\pid_A^B\to\Gamma$;
        \item[$\succ$] (\textit{Elim Rule}) a section: $$\ev_A^B:\Gamma.\pid_A^B.A[P_{\pid_A^B}]\to \Gamma.\pid_A^B.A[P_{\pid_A^B}].B[P_{\pid_A^B}^\lilbullet]$$ of the display map $\Gamma.\pid_A^B.A[P_{\pid_A^B}].B[P_{\pid_A^B}^\lilbullet]\to\Gamma.\pid_A^B.A[P_{\pid_A^B}]$;
    \end{itemize} and that, for every object $\Gamma$, every display map $P_A$ of codomain $\Gamma$, every display map $P_B$ of codomain $\Gamma.A$, and every section $\upsilon:\Gamma.A\to\Gamma.A.B$ of $P_B$ there is a choice of:
    \begin{itemize}
        \item[$\succ$] (\textit{Intro Rule}) a section: $$\abst{\upsilon} : \Gamma\to\Gamma.\pid_A^B$$ of the display map $\Gamma.\pid_A^B\to\Gamma$;
        \item[$\succ$] (\textit{Comp Axiom}) a section: $$\beta_\upsilon:\Gamma.A\to\Gamma.A.\id_{B}[\ev_A^B[\abst{\upsilon}^\lilbullet];\upsilon]$$ of the display map $\Gamma.A.\id_{B}[\ev_A^B[\abst{\upsilon}^\lilbullet];\upsilon]\to\Gamma.A$, where the arrow: $$\ev_A^B[\abst{\upsilon}^\lilbullet];\upsilon: \Gamma.A\to\Gamma.A.B.B^\liltriangle$$ is built as explained in Subsection \ref{subsection:Encoding parallel terms into a substitution}.
    \end{itemize} Moreover, let us assume that the following \textit{stability conditions}: $$\begin{alignedat}{2}
        \pid_A^B[f]&=\pid_{A[f]}^{B[f^\lilbullet]}&\quad\quad\abst{\upsilon}[f]&=\abst{\upsilon[f^{\lilbullet}]}\\
        \ev_A^B[f^{\lilbullet\lilbullet}]&=\ev_{A[f]}^{B[f^\lilbullet]}&\beta_\upsilon[f^{\lilbullet}]&=\beta_{\upsilon[f^{\lilbullet}]}\\
    \end{alignedat}$$ hold\footnotemark[1] for every arrow $f:\Delta\to\Gamma$.

    \noindent Then we say that $(\mathbf{C},\mathcal{D})$ is \textit{endowed with axiomatic $\Pi$-types}.

\medskip
    
    Additionally, let us assume that, for every object $\Gamma$, every display map $P_A$ of codomain $\Gamma$, every display map $P_B$ of codomain $\Gamma.A$, and every pair of sections $z,z'$ of $P_{\pid_A^B}$, there is a choice of:\begin{itemize}
        \item[$\succ$] (\textit{Intro Rule}) a section: $$\funext_q:\Gamma\to\Gamma.\id_{\pid_A^B}[z;z']$$ of the display map $\Gamma.A.\id_{\pid_A^B}[z;z']\to\Gamma.A$ for every section: $$q:\Gamma\to \Gamma.\pid_A^{\id_B[\ev_A^B[z^\lilbullet];\ev_A^B[z'^\lilbullet]]}$$ of the display map $\Gamma.\pid_A^{\id_B[\ev_A^B[z^\lilbullet];\ev_A^B[z'^\lilbullet]]}\to\Gamma$;
        \item[$\succ$] (\textit{Comp Axiom}) a section: $$\beta^\Pi_q:\Gamma\to\Gamma.\id_{\pid_A^{\id_B[\ev_A^B[z^\lilbullet];\ev_A^B[z'^\lilbullet]]}}[\;\happly_{z;z'}[\funext_q]\;;\;q\;]$$ of the display map $\Gamma.\id_{\pid_A^{\id_B[\ev_A^B[z^\lilbullet];\ev_A^B[z'^\lilbullet]]}}[\;\happly_{z;z'}[\funext_q]\;;\;q\;]\to\Gamma$ for every section: $$q:\Gamma\to \Gamma.\pid_A^{\id_B[\ev_A^B[z^\lilbullet];\ev_A^B[z'^\lilbullet]]}$$ of the display map $\Gamma.\pid_A^{\id_B[\ev_A^B[z^\lilbullet];\ev_A^B[z'^\lilbullet]]}\to\Gamma$;
        \item[$\succ$] (\textit{Exp Axiom}) a section: $$\eta^\Pi_p:\Gamma\to\Gamma.\id_{\id_{\pid_A^B}[z;z']}[\;p\;;\;\funext_{\happly_{z;z'}[p]}\;]$$ of the display map $\Gamma\to\Gamma.\id_{\id_{\pid_A^B}[z;z']}[\;p\;;\;\funext_{\happly_{z;z'}[p]}\;]$ for every section: $$p:\Gamma\to\Gamma.\id_{\pid_A^B}[z;z']$$ of the display map $\Gamma.\id_{\pid_A^B}[z;z']\to\Gamma$;
    \end{itemize} where:\begin{itemize}
        \item[$\succ$] $\mathsf{s}$ and $\mathsf{t}$ are $P_{(\pid_A^B)^\liltriangle}P_{\id_{\pid_A^B}}$ and $P_{\pid_A^B}^\lilbullet P_{\id_{\pid_A^B}}$ respectively;
        \item[$\succ$] the section: $$\happly:\Gamma.\pid_A^B.(\pid_A^B)^\liltriangle.\id_{\pid_A^B}\to\Gamma.\pid_A^B.(\pid_A^B)^\liltriangle.\id_{\pid_A^B}.\pid_{A[P_{\pid_A^B}\mathsf{s}]}^{\id_{B[(P_{\pid_A^B}\mathsf{s})^\lilbullet]}[\ev_A^B[\mathsf{s}^\lilbullet],\ev_A^B[\mathsf{t}^\lilbullet]]}$$ of the display map $\Gamma.\pid_A^B.(\pid_A^B)^\liltriangle.\id_{\pid_A^B}.\pid_{A[P_{\pid_A^B}\mathsf{s}]}^{\id_{B[(P_{\pid_A^B}\mathsf{s})^\lilbullet]}[\ev_A^B[\mathsf{s}^\lilbullet],\ev_A^B[\mathsf{t}^\lilbullet]]}\to\Gamma.\pid_A^B.(\pid_A^B)^\liltriangle.\id_{\pid_A^B}$ is defined in \Cref{howtohapply};
        \item[$\succ$] and $\happly_{z;z'}$ is $\happly[(z;z')^\lilbullet]$.
    \end{itemize} Moreover, let us assume that the following \textit{stability conditions}: $$\begin{aligned}
    \funext_q[f]&=\funext_{q[f]}\\
    \beta^\Pi_q[f]&=\beta^\Pi_{q[f]}\\
    \eta^\Pi_q[f]&=\eta^\Pi_{q[f]}
    \end{aligned}$$ hold\footnotemark[1] for every arrow $f:\Delta\to\Gamma$.

    \noindent Then we say that $(\mathbf{C},\mathcal{D})$ is \textit{endowed with axiomatic $\Pi$-types \& axiomatic function extensionality}.
\end{definition}

\begin{remark}\label{howtohapply} Let $\mathsf{s}$ and $\mathsf{t}$ be the arrows: $$P_{(\pid_A^B)^\liltriangle}P_{\id_{\pid_A^B}} \quad\quad\text{and}\quad\quad P_{\pid_A^B}^\lilbullet P_{\id_{\pid_A^B}}$$ respectively. By using the \textit{stability conditions} and the commutativity of Diagram (\ref{reindexingprimareindexingdopo}), the re-indexing of the display map: $$\Gamma.\pid_A^B.(\pid_A^B)^\liltriangle.\id_{\pid_A^B}.\pid_{A[P_{\pid_A^B}\mathsf{s}]}^{\id_{B[(P_{\pid_A^B}\mathsf{s})^\lilbullet]}[\ev_A^B[\mathsf{s}^\lilbullet],\ev_A^B[\mathsf{t}^\lilbullet]]}\to\Gamma.\pid_A^B.(\pid_A^B)^\liltriangle.\id_{\pid_A^B}$$ via the arrow $\refl_{\pid_A^B}:\Gamma.\pid_A^B\to\Gamma.\pid_A^B.(\pid_A^B)^\liltriangle.\id_{\pid_A^B}$ is proven to be the display map: $$\Gamma.\pid_A^B.\pid_{A[P_{\pid_A^B}]}^{\id_{B[P_{\pid_A^B}^\lilbullet]}[\ev_A^B;\ev_A^B]}\to\Gamma.\pid_A^B\;\;.$$ Now, if we build a section $c$ of the latter display map, then a section $\happly$ of the former display map can be defined as $\j_c$. Additionally, if we build a section $\upsilon$ of the display map: $$\Gamma.\pid_A^B.{A[P_{\pid_A^B}]}.{\id_{B[P_{\pid_A^B}^\lilbullet]}[\ev_A^B;\ev_A^B]}\to\Gamma.\pid_A^B.{A[P_{\pid_A^B}]}$$ then $c$ can be defined as $\abst{\upsilon}$. Finally, we define $\upsilon$ as the unique section $\refl_{B[P_{\pid_A^B}^\lilbullet]}[\ev_A^B]$ of this display map that makes the diagram: 
\[\begin{tikzcd}
	{\Gamma.\pid_A^B.{A[P_{\pid_A^B}]}} && {\Gamma.\pid_A^B.{A[P_{\pid_A^B}]}.{\id_{B[P_{\pid_A^B}^\lilbullet]}[\ev_A^B;\ev_A^B]}} \\
	\\
	{\Gamma.\pid_A^B.{A[P_{\pid_A^B}]}.B[P_{\pid_A^B}^\lilbullet]} && {\Gamma.\pid_A^B.{A[P_{\pid_A^B}]}.B[P_{\pid_A^B}^\lilbullet].B[P_{\pid_A^B}^\lilbullet]^\liltriangle.\id_{B[P_{\pid_A^B}^\lilbullet]}}
	\arrow[from=1-1, to=1-3]
	\arrow["{\ev_A^B}"{description}, from=1-1, to=3-1]
	\arrow["{(\ev_A^B;\ev_A^B)^\liltriangle}"{description}, from=1-3, to=3-3]
	\arrow["{\refl_{B[P_{\pid_A^B}^\lilbullet]}}", from=3-1, to=3-3]
\end{tikzcd}\] commute. Summarising, the section $\happly$ is defined as: $$\j_{\abst{\refl_{B[P_{\pid_A^B}^\lilbullet]}[\ev_A^B]}}$$ and it is the interpretation of the term judgement of type $\Pi_{x:A}^{\ev(z,x)=\ev(z',x)}$ in context: $$\judgectx{z,z':\Pi_{x:A}B(x);\;p:z=z'}$$ obtained by elimination via the term judgement $\judge{z:\Pi_{x:A}B(x)}{\lambda x.\refl(\ev(z,x))}$.

Now, by using again the \textit{stability conditions} and the commutativity of Diagram (\ref{reindexingprimareindexingdopo}), the re-indexing of the display map: $$\Gamma.\pid_A^B.(\pid_A^B)^\liltriangle.\id_{\pid_A^B}.\pid_{A[P_{\pid_A^B}\mathsf{s}]}^{\id_{B[(P_{\pid_A^B}\mathsf{s})^\lilbullet]}[\ev_A^B[\mathsf{s}^\lilbullet],\ev_A^B[\mathsf{t}^\lilbullet]]}\to\Gamma.\pid_A^B.(\pid_A^B)^\liltriangle.\id_{\pid_A^B}$$ via the arrow $(z;z')^\lilbullet p:\Gamma.\pid_A^B\to\Gamma.\pid_A^B.(\pid_A^B)^\liltriangle.\id_{\pid_A^B}$---for some section: $p:\Gamma\to\Gamma.\id_{\pid_A^B}[z;z']$ of the display map $\Gamma.\id_{\pid_A^B}[z;z']\to\Gamma$---is proven to be the display map: $$\Gamma.\pid_A^{\id_B[\ev_A^B[z^\lilbullet];\ev_A^B[z'^\lilbullet]]}\to\Gamma\;\;$$ hence $\happly_{z;z'}[p]$ is a section of this display map.

\end{remark}

\begin{figure*}[p]
\begin{center}
\fbox{$\begin{alignedat}{2}
\text{Form Rule\,}&\inferrule{\;\;\quad\quad\quad\text{ }\quad\quad\quad\;\;}{0:\type}&\quad\text{Elim Rule\,}&\inferrule{\judge{x:0}{C(x):\type}}{\judge{x:0}{\ind^0(x):C(x)}}\\
\\
\\
\\
\\
\text{Form Rule\,}&\inferrule{\;\;\quad\quad\quad\text{ }\quad\quad\quad\;\;}{1:\type}&\quad\text{Elim Rule\,}&\inferrule{\judge{x:1}{C(x):\type}\\c:C(\star)}{\judge{x:1}{\ind^1(c,x):C(x)}}\\
\\
\text{Intro Rule\,}&\inferrule{\;\;\quad\quad\quad\text{ }\quad\quad\quad\;\;}{\star:1}&\quad\text{Comp Axiom\,}&\inferrule{\judge{x:1}{C(x):\type}\\c:C(\star)}{{\beta^1(c):\ind^1(c,\star)=c}}\\
\\
\\
\\
\\
\text{Form Rule\,}&\inferrule{\;\;\quad\quad\quad\text{ }\quad\quad\quad\;\;}{2:\type}&\quad\text{Elim Rule\,}&\inferrule{\judge{x:2}{C(x):\type}\\\\\quad\; c:C(\bot)\\d:C(\top) \quad\;}{\judge{x:2}{\ind^2(c,d,x):C(x)}}\\
\\
\text{Intro Rule\,}&\inferrule{\;\;\quad\quad\quad\text{ }\quad\quad\quad\;\;}{\bot:2\\\\\top:2}&\quad\text{Comp Axiom\,}&\inferrule{\judge{x:2}{C(x):\type}\\\\\quad\; c:C(\bot)\\d:C(\top) \quad\;}{\beta^2_\bot(c,d):\ind^2(c,d,\bot)=c\\\\\beta^2_\top(c,d):\ind^2(c,d,\top)=d}\\
\\
\\
\\
\\
\text{Form Rule\,}&\inferrule{\;\;\quad\quad\quad\text{ }\quad\quad\quad\;\;}{\nat:\type}&\quad\text{Elim Rule\,}&\inferrule{\judge{n:\nat}{C(n):\type}\\c:C(\0)\\\\\;\;\judge{n:\nat,y:C(n)}{d(n,y):C(\s(n))}\;\;}{\judge{n:\nat}{\ind^\nat(c,d,n):C(n)}}\\
\\
\text{Intro Rule\,}&\inferrule{\;\;\quad\quad\quad\text{ }\quad\quad\quad\;\;}{\0:\nat\\\\\judge{n:\nat}{\s(n):\nat}}&\quad\text{Comp Axiom\,}&\inferrule{\judge{n:\nat}{C(n):\type}\\c:C(\0)\\\\\;\;\judge{n:\nat,y:C(n)}{d(n,y):C(\s(n))}\;\;}{\beta^\nat_\0(c,d):\ind^{\nat}(c,d,\0)=c\\\\\judge{n:\nat}{\beta^\nat_\s(c,d,n):}\\\\{\ind_{\nat}(c,d,\s(n))=d(n,\ind^{\nat}(c,d,n))}}
\end{alignedat}$}
\caption{Axiomatic $0$-, $1$-, $2$-, and $\nat$-types}\label{propositional nat}
\end{center}
\end{figure*}

Finally, referring to the rules of Fig. \ref{propositional nat} for axiomatic axiomatic $0$-, $1$-, $2$-, $\nat$-types, we give the following:

\begin{definition}[Semantics of axiomatic $0$-, $1$-, $2$-, $\nat$-types---syntactic formulation]\label{def:Semantics of axiomatic Nat-types}
    Let $(\mathbf{C},\mathcal{D})$ be a display map category endowed with axiomatic $=$-types.

\begin{itemize}
    
    \item Let us assume that, for every object $\Gamma$, there is a choice of:
    \begin{itemize}
        \item[$\succ$] (\textit{Form Rule}) a display map $\Gamma.0\to\Gamma$;
    \end{itemize} and that, for every display map $P_C$ of codomain $\Gamma.0$, there is a choice of:
    \begin{itemize}
        \item[$\succ$] (\textit{Elim Rule}) a section: $$\ind^0 : \Gamma.0\to\Gamma.0.C$$ of the display map $\Gamma.0.C\to\Gamma.0$.
    \end{itemize} Moreover, let us assume that the following \textit{stability conditions}: $$\begin{alignedat}{2}
        0[f]&=0&\quad\quad\ind^0[f^\lilbullet]&=\ind^0
    \end{alignedat}$$ hold\footnotemark[1] for every arrow $f:\Delta\to\Gamma$.
    
    \noindent Then we say that $(\mathbf{C},\mathcal{D})$ is \textit{endowed with axiomatic $0$-types}.

    \item Let us assume that, for every object $\Gamma$, there is a choice of:
    \begin{itemize}
        \item[$\succ$] (\textit{Form Rule}) a display map $\Gamma.1\to\Gamma$;
        \item[$\succ$] (\textit{Intro Rule}) a section: $$\star:\Gamma\to \Gamma.1$$ of the display map $\Gamma.1\to\Gamma$;
    \end{itemize} and that, for every display map $P_C$ of codomain $\Gamma.1$, and every section $c:\Gamma\to\Gamma.C[\star]$ of $P_{C[\star]}$, there is a choice of:
    \begin{itemize}
        \item[$\succ$] (\textit{Elim Rule}) a section: $$\ind^1_c : \Gamma.1\to\Gamma.1.C$$ of the display map $\Gamma.1.C\to\Gamma.1$;
        \item[$\succ$] (\textit{Comp Axiom}) a section: $$\beta^1_c:\Gamma\to\Gamma.\id_{C[\star]}[\ind^1_c[\star];c]$$ of the display map $\Gamma.\id_{C[\star]}[\ind^1_c[\star];c]\to\Gamma$, where the arrow: $$\ind^1_c[\star];c: \Gamma\to\Gamma.C[\star].C[\star]^\liltriangle$$ is built as explained in Subsection \ref{subsection:Encoding parallel terms into a substitution}.
    \end{itemize} Moreover, let us assume that the following \textit{stability conditions}: $$\begin{alignedat}{2}
        1[f]&=1&\quad\quad\ind^1_c[f^\lilbullet]&=\ind^1_{c[f]}\\
        \star[f]&=\star&\beta^1_c[f]&=\beta^1_{c[f]}\\
    \end{alignedat}$$ hold\footnotemark[1] for every arrow $f:\Delta\to\Gamma$.

    \noindent Then we say that $(\mathbf{C},\mathcal{D})$ is \textit{endowed with axiomatic $1$-types}.

    \item Let us assume that, for every object $\Gamma$, there is a choice of:
    \begin{itemize}
        \item[$\succ$] (\textit{Form Rule}) a display map $\Gamma.2\to\Gamma$;
        \item[$\succ$] (\textit{Intro Rule}) sections: $$\bot:\Gamma\to \Gamma.2 \quad\quad\quad\quad \top:\Gamma\to \Gamma.2$$ of the display map $\Gamma.2\to\Gamma$;
    \end{itemize} and that, for every display map $P_C$ of codomain $\Gamma.2$, every section $c:\Gamma\to\Gamma.C[\bot]$ of $P_{C[\bot]}$, and every section $d:\Gamma\to\Gamma.C[\top]$ of $P_{C[\top]}$, there is a choice of:
    \begin{itemize}
        \item[$\succ$] (\textit{Elim Rule}) a section: $$\ind^2_{c,d} : \Gamma.2\to\Gamma.2.C$$ of the display map $\Gamma.2.C\to\Gamma.2$;
        \item[$\succ$] (\textit{Comp Axiom}) sections: $$\beta^{2,\bot}_{c,d}:\Gamma\to\Gamma.\id_{C[\bot]}[\ind^2_{c,d}[\bot];c] \quad\quad\quad\quad \beta^{2,\top}_{c,d}:\Gamma\to\Gamma.\id_{C[\top]}[\ind^2_{c,d}[\top];d]$$ of the display maps $\Gamma.\id_{C[\bot]}[\ind^2_{c,d}[\bot];c]\to\Gamma$ and $\Gamma.\id_{C[\top]}[\ind^2_{c,d}[\top];c]\to\Gamma$ respectively, where the arrows: $$\ind^2_{c,d}[\bot];c: \Gamma\to\Gamma.C[\bot].C[\bot]^\liltriangle \quad\quad\quad\quad \ind^2_{c,d}[\top];c: \Gamma\to\Gamma.C[\top].C[\top]^\liltriangle$$ are built as explained in Subsection \ref{subsection:Encoding parallel terms into a substitution}.
    \end{itemize} Moreover, let us assume that the following \textit{stability conditions}: $$\begin{alignedat}{2}
        2[f]&=2&\quad\quad\ind^2_{c,d}[f^\lilbullet]&=\ind^2_{c[f],d[f]}\\
        \bot[f]&=\bot&\beta^{2,\bot}_{c,d}[f]&=\beta^{2,\bot}_{c[f],d[f]}\\
        \top[f]&=\top&\beta^{2,\top}_{c,d}[f]&=\beta^{2,\top}_{c[f],d[f]}\\
    \end{alignedat}$$ hold\footnotemark[1] for every arrow $f:\Delta\to\Gamma$.

    \noindent Then we say that $(\mathbf{C},\mathcal{D})$ is \textit{endowed with axiomatic $2$-types}.

    \item Let us assume that, for every object $\Gamma$, there is a choice of:
    \begin{itemize}
        \item[$\succ$] (\textit{Form Rule}) a display map $\Gamma.\nat\to\Gamma$;
        \item[$\succ$] (\textit{Intro Rule}) sections: $$\0:\Gamma\to \Gamma.\nat \quad\quad\quad\quad \successivo:\Gamma.\nat\to \Gamma.\nat.\nat^\liltriangle$$ of the display map $\Gamma.\nat\to\Gamma$ and $\Gamma.\nat.\nat^\liltriangle\to\Gamma.\nat$ respectively;
    \end{itemize} and that, for every display map $P_C$ of codomain $\Gamma.\nat$, every section $c:\Gamma\to\Gamma.C[\0]$ of $P_{C[\0]}$, and every section $d:\Gamma.\nat.C\to\Gamma.\nat.C.C[\s P_C]$ of $P_{C[\s P_C]}$---where $\s$ is the composition: 
\[\begin{tikzcd}[column sep=small,row sep=tiny]
	{\Gamma.\nat} &&&& {\Gamma.\nat.\nat^\liltriangle} &&&& {\Gamma.\nat}
	\arrow["{\successivo}"{description}, from=1-1, to=1-5]
	\arrow["{P_\nat^\lilbullet}"{description}, from=1-5, to=1-9]
\end{tikzcd}\]---there is a choice of:
    \begin{itemize}
        \item[$\succ$] (\textit{Elim Rule}) a section: $$\ind^\nat_{c,d} : \Gamma.\nat\to\Gamma.\nat.C$$ of the display map $\Gamma.\nat.C\to\Gamma.\nat$;
        \item[$\succ$] (\textit{Comp Axiom}) sections: $$\beta^{\nat,\0}_{c,d}:\Gamma\to\Gamma.\id_{C[\0]}[\ind^\nat_{c,d}[\0];c] \quad\quad\quad\quad \beta^{\nat,\s}_{c,d}:\Gamma.\nat\to\Gamma.\nat.\id_{C[\s]}[\ind^\nat_{c,d}[\s ];d[\ind^\nat_{c,d}]]$$ of the display maps $\Gamma.\id_{C[\0]}[\ind^\nat_{c,d}[\0];c]\to\Gamma$ and $\Gamma.\nat.\id_{C[\s]}[\ind^\nat_{c,d}[\s];d[\ind^\nat_{c,d}]]\to\Gamma.\nat$ respectively, where the arrows: $$\ind^\nat_{c,d}[\0];c: \Gamma\to\Gamma.C[\0].C[\0]^\liltriangle \quad\quad\quad\quad \ind^\nat_{c,d}[\s ];d[\ind^\nat_{c,d}]: \Gamma.\nat\to\Gamma.\nat.C[\s ].C[\s ]^\liltriangle$$ are built as explained in Subsection \ref{subsection:Encoding parallel terms into a substitution}.
    \end{itemize} Moreover, let us assume that the following \textit{stability conditions}: $$\begin{alignedat}{2}
        \nat[f]&=\nat&\quad\quad\ind^\nat_{c,d}[f^\lilbullet]&=\ind^\nat_{c[f],d[f^{\lilbullet\lilbullet}]}\\
        \0[f]&=\0&\beta^{\nat,\0}_{c,d}[f]&=\beta^{\nat,\0}_{c[f],d[f^{\lilbullet\lilbullet}]}\\
        \successivo[f]&=\successivo&\beta^{\nat,\s}_{c,d}[f^\lilbullet]&=\beta^{\nat,\s}_{c[f],d[f^{\lilbullet\lilbullet}]}\\
    \end{alignedat}$$ hold\footnotemark[1] for every arrow $f:\Delta\to\Gamma$.

    \noindent Then we say that $(\mathbf{C},\mathcal{D})$ is \textit{endowed with axiomatic $\nat$-types}.
\end{itemize}    
\end{definition}

\subsection{Interpreting ATT into a display map category}

Let $(\mathbf{C},\mathcal{D})$ be a display map category endowed with axiomatic $=$-types, axiomatic $\Sigma$-types, axiomatic $\Pi$-types \& axiomatic function extensionality, and axiomatic $0$-, $1$-, $2$-, $\nat$-types---see \Cref{def:Semantics of axiomatic =-types,def:Semantics of axiomatic Sigma-types,def:Semantics of axiomatic Pi-types,def:Semantics of axiomatic Nat-types}, respectively. Then a notion of interepretation of ATT, with $=$-, $\Sigma$-, $\Pi$-, $0$-, $1$-, $2$- $\nat$-types, and function extensionality, in $(\mathbf{C},\mathcal{D})$---in the sense of Section \ref{section:Recap on display map categories}---is defined. The proof of this fact is standard and adapts the argument provided by Streicher \cite{MR1134134} and Hofmann \cite{Hofmann1997}. We also refer the reader to \cite{phdcartmell,Streicher1991,phdhofmann} for further details. The argument consists in defining an a priori partial interpretation function, whose domain consists of the \textit{pre-judgements} of the raw syntax of ATT, and which can be shown by induction on derivations to be well-defined on the actual judgements of ATT. In detail: \begin{itemize}
    \item a context $\gamma:\Gamma$ of ATT is interpreted as an object $\Gamma$ of $\mathbf{C}$;
    \item if $\Gamma$ is the interpretation of the context $\gamma$, then a type judgment $\judge{\gamma:\Gamma}{A(\gamma):\type}$ of ATT is interpreted as a display map $\Gamma.A \to\Gamma$ of $\mathcal{D}$;
    \item if $\Gamma.A \to \Gamma$ is the interpretation of the type judgment $\judge{\gamma:\Gamma}{A(\gamma):\type}$, then a term judgment $\judge{\gamma:\Gamma}{a(\gamma):A(\gamma)}$ of ATT is interpreted as a section $\Gamma\xrightarrow{a}\Gamma.A$ of $\Gamma.A \xrightarrow{P_A}\Gamma$.
\end{itemize} We only mention that the interpretation function satisfies the following expected clauses. Let us assume that---when required---$\Gamma$ is the interpretation of the context $\judgectx{\gamma:\Gamma}$, $P_A$ is the interpretation of the type judgement $\judge{\gamma:\Gamma}{A(\gamma):\type}$, and $P_B$ is the interpretation of the type judgement $\judge{\gamma:\Gamma;\;x:A(\gamma)}{B(\gamma,x):\type}$. Then:

\medskip

\noindent\textit{Clauses for context interpretation.}

\begin{itemize}
    \item the object $\epsilon$ is the interpretation of the context $\judgectx{\blank}$;
    \item the object $\Gamma.A$ is the interpretation of the context $\judgectx{\gamma;\;x:A(\gamma)}$;
\end{itemize}

\medskip

\noindent\textit{Clauses for type judgement interpretation.}

\begin{itemize}
    \item the display map $\Gamma.A.A^\liltriangle.\id_A\to\Gamma.A.A^\liltriangle$ is the interpretation of $\judge{\gamma;\;x,x':A(\gamma)}{x=x':\type}$;
    \item the display map $\Gamma.\sigmad_A^B\to\Gamma$ is the interpretation of $\judge{\gamma}{\Sigma_{x:A(\gamma)}B(\gamma,x):\type}$;
    \item the display map $\Gamma.\pid_A^B\to\Gamma$ is the interpretation of $\judge{\gamma}{\Pi_{x:A(\gamma)}B(\gamma,x):\type}$;
    \item the display map $\Gamma.0\to\Gamma$ is the interpretation of $\judge{\gamma}{0:\type}$;
    \item the display map $\Gamma.1\to\Gamma$ is the interpretation of $\judge{\gamma}{1:\type}$;
    \item the display map $\Gamma.2\to\Gamma$ is the interpretation of $\judge{\gamma}{2:\type}$;
    \item the display map $\Gamma.\nat\to\Gamma$ is the interpretation of $\judge{\gamma}{\nat:\type}$;
\end{itemize}

\medskip

\noindent\textit{Clauses for term judgement interpretation.} \begin{itemize}
    \item the section $\reflex_A$ of $\Gamma.A.\id_A[\delta_A]\to\Gamma.A$ is the interpretation of $\judge{\gamma;\;x:A(\gamma)}{\refl(x):x=x}$;
    \item the section $\pairing_A^B$ of the display map $\Gamma.A.B.\sigmad_A^B[P_AP_B]\to\Gamma.A.B$ is the interpretation of $\judge{\gamma;\;x:A(\gamma);\;y:B(\gamma,x)}{\langle x,y\rangle:\Sigma_{x:A(\gamma)}B(\gamma,x)}$;
    \item the section $\ev_A^B$ of the display map $\Gamma.\pid_A^B.A[P_{\pid_A^B}].B[P_{\pid_A^B}^\lilbullet]\to\Gamma.\pid_A^B.A[P_{\pid_A^B}]$ is the interpretation of $\judge{\gamma;\;z:\Pi_{x:A(\gamma)}B(\gamma,x);\;x:A}{\ev(z,x):B(x)}$;
    \item the section $\star$ of the display map $\Gamma\to\Gamma.1$ is the interpretation of $\judge{\gamma}{\star:1}$;
    \item the section $\bot$ of the display map $\Gamma\to\Gamma.2$ is the interpretation of $\judge{\gamma}{\bot:2}$;
    \item the section $\top$ of the display map $\Gamma\to\Gamma.2$ is the interpretation of $\judge{\gamma}{\top:2}$;
    \item the section $\0$ of the display map $\Gamma\to\Gamma.\nat$ is the interpretation of $\judge{\gamma}{\0:\nat}$;
    \item the section $\successivo$ of the display map $\Gamma.\nat\to\Gamma.\nat.\nat^\liltriangle$ is the interpretation of the term judgement $\judge{\gamma;\,n:\nat}{\s(n):\nat}$.
\end{itemize} Additionally, let us add to the assumptions, for the following two clauses, that $P_C$ is the interpretation of the type judgement $\judge{\gamma;x,x';\;p:x=x'}{C(\gamma,x,x',p):\type}$ and $\Gamma.A\xrightarrow{c}\Gamma.A.C[\refl_A]$ is the interpretation of the term judgement $\judge{\gamma;x}{c(\gamma,x):C(\gamma,x,x,\refl(x))}$. Then:
\begin{itemize}
    \item the section $\j_c$ of the display map $\Gamma.A.A^\liltriangle.\id_A.C\to\Gamma.A.A^\liltriangle.\id_A$ is the interpretation of $\judge{\gamma;x,y;p}{\j(c,\gamma,x,x',p):C(\gamma,x,x',p)}$;
    \item the section $\h_c$ of $\Gamma.A.\id_{C[\refl_A]}[\j_c[\refl_A];c]\to\Gamma.A$ is the interpretation of the term judgement $\judge{\gamma;x}{\h(c,\gamma,x):\j(c,\gamma,x,x,\refl(x))=c(\gamma,x)}$.
\end{itemize} Additionally, let us add to the assumptions, for the following two clauses, that $P_C$ is the interpretation of the type judgement $\judge{\gamma;\;u:\Sigma_{x:A(\gamma)}B(\gamma,x)}{C(\gamma,u):\type}$ and $\Gamma.A.B\xrightarrow{c}\Gamma.A.B.C[\pair_A]$ is the interpretation of the term judgement $\judge{\gamma;\;x;\;y}{c(\gamma,x,y):C(\gamma,\langle x,y\rangle)}$. Then:
\begin{itemize}
    \item the section $\splitt_c$ of $\Gamma.\sigmad_A^B.C\to\Gamma.\sigmad_A^B$ is the interpretation of $\judge{\gamma;u}{\splitt(c,\gamma,u):C(\gamma,u)}$;
    \item the section $\sigma_c$ of the display map $\Gamma.A.B.\id_{C[\pair_A^B]}[\splitt_c[\pair_A^B];c]\to\Gamma.A.B$ is the interpretation of $\judge{\gamma;x;y}{\sigma(c,\gamma,x,y):\splitt(c,\gamma,\langle x,y\rangle)=c(\gamma,x,y)}$.
\end{itemize} Additionally, let us add to the assumptions, for the following two clauses, that $\upsilon:\Gamma.A\to\Gamma.A.B$ is the interpretation of the term judgement $\judge{\gamma;x}{\upsilon(\gamma,x):B(\gamma,x)}$. Then: \begin{itemize}
    \item the section $\abst{\upsilon}$ of $\Gamma\to\Gamma.\pid_A^B$ is the interpretation of $\judge{\gamma}{\lambda x. \upsilon(\gamma,x):\Pi_{x:A(\gamma)}B(\gamma,x)}$;
    \item the section $\beta_q$ of the display map $\Gamma.A\to\Gamma.A.\id_B[\ev_A^B[\abst{\upsilon}^\lilbullet];\upsilon]$ is the interpretation of $\judge{\gamma;x}{\beta(\upsilon,\gamma,x):\ev(\lambda x.\upsilon(\gamma,x),x)=\upsilon(\gamma,x)}$.
\end{itemize} Additionally, let us add to the assumptions, for the following two clauses, that $z,z':\Gamma\to\Gamma.\Pi_{x:A(\gamma)}B(\gamma,x)$ are the interpretations of $\judge{\gamma}{z(\gamma),z'(\gamma):\Pi_{x:A(\gamma)}B(\gamma,x)}$ and that $q:\Gamma\to \Gamma.\pid_A^{\id_B[\ev_A^B[z^\lilbullet];\ev_A^B[z'^\lilbullet]]}$ is the interpretation of $\judge{\gamma}{q(\gamma):\Pi_{x:A(\gamma)}\ev(z(\gamma),x)=\ev(z'(\gamma),x)}$. Then: \begin{itemize}
    \item the section $\funext_q$ of the display map $\Gamma.\id_{\pid_A^B}[z;z']\to\Gamma$ is the interpretation of the term judgement $\judge{\gamma}{\funext(z(\gamma),z'(\gamma),q(\gamma)):z(\gamma)=z'(\gamma)}$;
    \item the section $\beta^\Pi_q$ of the display map $\Gamma.\id_{\pid_A^{\id_B[\ev_A^B[z^\lilbullet];\ev_A^B[z'^\lilbullet]]}}[\;\happly_{z;z'}[\funext_q]\;;\;q\;]\to\Gamma$ is the interpretation of $\judge{\gamma}{\beta^\Pi(z(\gamma),z'(\gamma),q(\gamma)):\lambda x.\ev(\funext(z(\gamma),z'(\gamma),q(\gamma)),x)=q(\gamma)}$.
\end{itemize} Additionally, let us add to the assumptions, for the following clause, that the sections $z,z':\Gamma\to\Gamma.\Pi_{x:A(\gamma)}B(\gamma,x)$ are the interpretations of $\judge{\gamma}{z(\gamma),z'(\gamma):\Pi_{x:A(\gamma)}B(\gamma,x)}$ and that $p:\Gamma\to\Gamma.\id_{\pid_A^B}[z;z']$ is the interpretation of $\judge{\gamma}{p(\gamma):z(\gamma)=z'(\gamma)}$. Then: \begin{itemize}
    \item the section $\eta^\Pi_p$ of the display map $\Gamma\to\Gamma.\id_{\id_{\pid_A^B}[z;z']}[\;p\;;\;\funext_{\happly_{z;z'}[p]}\;]$ is the interpretation of $\judge{\gamma}{\eta^\Pi(z(\gamma),z'(\gamma),p(\gamma)):p(\gamma)=\funext(z(\gamma),z'(\gamma),\lambda x.\ev(p(\gamma),x))}$.
\end{itemize} Additionally, let us add to the assumptions, for the following clause, that $P_C$ is the interpretation of the type judgement $\judge{\gamma;\,x:0}{C(\gamma,x):\type}$. Then: \begin{itemize}
    \item the section $\ind^0$ of $\Gamma.0.C\to\Gamma.0$ is the interpretation of $\judge{\gamma;\,x:0}{\ind^0(\gamma,x):C(x)}$.
\end{itemize} Additionally, let us add to the assumptions, for the following two clauses, that $P_C$ is the interpretation of the type judgement $\judge{\gamma;\,x:1}{C(\gamma,x):\type}$ and $c: \Gamma\to\Gamma.C[\star]$ is the interpretation of the term judgement $\judge{\gamma}{c(\gamma):C(\gamma,\star)}$. Then: \begin{itemize}
    \item the section $\ind^1_c$ of $\Gamma.1.C\to\Gamma.1$ is the interpretation of $\judge{\gamma;\,x:1}{\ind^1(c,\gamma,x):C(x)}$;
    \item the section $\beta^1_c$ of $\Gamma.\id_{C[\star]}[\ind^1_c[\star];c]\to\Gamma$ is the interpretation of the term judgement $\judge{\gamma}{\beta^1(c,\gamma):\ind^1(c,\gamma,\star)=c(\gamma)}$.
\end{itemize} Additionally, let us add to the assumptions, for the following three clauses, that $P_C$ is the interpretation of the type judgement $\judge{\gamma;\,x:2}{C(\gamma,x):\type}$ and $c: \Gamma\to\Gamma.C[\bot]$ is the interpretation of the term judgement $\judge{\gamma}{c(\gamma):C(\gamma,\bot)}$ and $d: \Gamma\to\Gamma.C[\top]$ is the interpretation of the term judgement $\judge{\gamma}{d(\gamma):C(\gamma,\top)}$. Then: \begin{itemize}
    \item the section $\ind^2_{c,d}$ of $\Gamma.2.C\to\Gamma.2$ is the interpretation of $\judge{\gamma;\,x:2}{\ind^2(c,d,\gamma,x):C(\gamma,x)}$;
    \item the section $\beta^{2,\bot}_{c,d}$ of $\Gamma.\id_{C[\bot]}[\ind^2_{c,d}[\bot];c]\to\Gamma$ is the interpretation of the term judgement $\judge{\gamma}{\beta^2_\bot(c,d,\gamma):\ind^2(c,d,\gamma,\bot)=c(\gamma)}$;
    \item the section $\beta^{2,\top}_{c,d}$ of $\Gamma.\id_{C[\top]}[\ind^2_{c,d}[\top];d]\to\Gamma$ is the interpretation of the term judgement $\judge{\gamma}{\beta^2_\top(c,d,\gamma):\ind^2(c,d,\gamma,\top)=d(\gamma)}$.
\end{itemize} Finally, let us add to the assumptions, for the last three clauses, that $P_C$ is the interpretation of the type judgement $\judge{\gamma;\,n:\nat}{C(\gamma,n):\type}$ and $c: \Gamma\to\Gamma.C[\0]$ is the interpretation of the term judgement $\judge{\gamma}{c(\gamma):C(\gamma,\0)}$ and $d: \Gamma.\nat.C\to\Gamma.\nat.C.C[\s P_C]$ is the interpretation of the term judgement $\judge{\gamma;\,n:\nat,y:C(n)}{d(\gamma,n,y):C(\gamma,\s(n))}$. Then: \begin{itemize}
    \item the section $\ind^\nat_{c,d}$ of $\Gamma.\nat.C\to\Gamma.\nat$ is the interpretation of $\judge{\gamma;\,n:\nat}{\ind^\nat(c,d,\gamma,n):C(n)}$;
    \item the section $\beta^{\nat,\0}_{c,d}$ of $\Gamma.\id_{C[\0]}[\ind^\nat_{c,d}[\0];c]\to\Gamma$ is the interpretation of the term judgement $\judge{\gamma}{\beta^\nat_\0(c,d,\gamma):\ind^\nat(c,d,\gamma,\0)=c(\gamma)}$;
    \item the section $\beta^{\nat,\s}_{c,d}$ of $\Gamma.\nat.\id_{C[\s]}[\ind^\nat_{c,d}[\s];d[\ind^\nat_{c,d}]]\to\Gamma.\nat$ is the interpretation of the term judgement $\judge{\gamma;\,n:\nat}{\beta^\nat_\s(c,d,\gamma,n):\ind^\nat(c,d,\gamma,\s(n))=d(\gamma,n,\ind^\nat(c,d,\gamma,n))}$.
\end{itemize}

\medskip

This notion of interpretation of ATT in $(\mathbf{C},\mathcal{D})$ satisfies the following:

\begin{theorem}[Soundness property] \label{soundness property} The interpretation of \textnormal{ATT} in $(\mathbf{C},\mathcal{D})$ is sound, that is:
\begin{itemize}
    \item whenever \textnormal{ATT} infers $\judgectx{\gamma:\Gamma}$, then its interpretation---let us indicate it as $\Gamma$---is defined;
    \item whenever \textnormal{ATT} infers $\judgectx{\gamma:\Gamma} \equiv \judgectx{\gamma':\Gamma'}$, for some contexts $\judgectx{\gamma:\Gamma}$ and $\judgectx{\gamma':\Gamma'}$, then their interpretations $\Gamma$ and $\Gamma'$ coincide;
    \item whenever \textnormal{ATT} infers $\judge{\gamma:\Gamma}{A(\gamma):\type}$, then its interpretation---let us indicate it as $P_A$---is defined;
    \item whenever \textnormal{ATT} infers $\judge{\gamma:\Gamma}{a(\gamma):A(\gamma)}$, then its interpretation---let us indicate it as $a:\Gamma\to\Gamma.A$---is defined;
    \item whenever \textnormal{ATT} infers $\judge{\gamma:\Gamma}{A(\gamma)\equiv A'(\gamma)}$ for some type judgements $\judge{\gamma:\Gamma}{A(\gamma):\type}$ and $\judge{\gamma:\Gamma}{A'(\gamma):\type}$, then their interpretations $P_A$ and $P_{A'}$ coincide;
    \item whenever \textnormal{ATT} infers $\judge{\gamma:\Gamma}{a(\gamma)\equiv a'(\gamma)}$ for some term judgements $\judge{\gamma:\Gamma}{a(\gamma):A(\gamma)}$ and $\judge{\gamma:\Gamma}{a'(\gamma):A(\gamma)}$ and some type judgement $\judge{\gamma:\Gamma}{A(\gamma):\type}$, then their interpretations $a$ and $a'$---as sections of $P_A$---coincide.
\end{itemize}

\proof
The proof of such a result is standard and we refer the reader to \cite{MR1134134,Streicher1991,phdhofmann,Hofmann1997} for details. The argument is by induction on the derivation of type, term, type equality, and term equality judgements.

A preliminary step is proving three properties called \textit{weakening property}, \textit{substitution property}, and \textit{general substitution property}:\begin{itemize}
    \item \textit{Weakening property.} Let $\Gamma$ be the interpretation of $\judgectx{\gamma:\Gamma}$, $P_{\Delta_{k+1}}$ be the interpretation of: $$\judge{\gamma:\Gamma,\delta_1:\Delta_1(\gamma),...,\delta_k(\gamma,\delta_1,...,\delta_{k-1})}{\Delta_{k+1}(\gamma,\delta_1,...,\delta_k):\type}$$ for all $k\in\{0,...,n-1\}$ for some $n\geq 0$, and $P_D$ be the interpretation of $\judge{\gamma:\Gamma}{D(\gamma):\type}$.
    
    Let $P_A$ be the interpretation of $\judge{\gamma:\Gamma,\delta:\Delta(\gamma)}{A(\gamma,\delta):\type}$, and $a$ be the interpretation of $\judge{\gamma:\Gamma,\delta:\Delta(\gamma)}{a(\gamma,\delta):A(\gamma,\delta)}$.
    
    Then the interpretations of: \begin{center}$\judge{\gamma:\Gamma,d:D(\gamma),\delta:\Delta(\gamma)}{A(\gamma,\delta):\type}$\;\;\;and\;\;\;$\judge{\gamma:\Gamma,d:D(\gamma),\delta:\Delta(\gamma)}{a(\gamma,\delta):A(\gamma,\delta)}$\end{center} are defined and are $P_{A[P_D.\Delta]}$ and $a[P_D.\Delta]$ respectively.

    \item \textit{Substitution property.} Let $\Gamma$ be the interpretation of $\judgectx{\gamma:\Gamma}$, $P_D$ be the interpretation of $\judge{\gamma:\Gamma}{D(\gamma):\type}$, $P_{\Delta_{k+1}}$ be the interpretation of: $$\judge{\gamma:\Gamma,d:D(\gamma),\delta_1:\Delta_1(\gamma,d),...,\delta_k(\gamma,d,\delta_1,...,\delta_{k-1})}{\Delta_{k+1}(\gamma,d,\delta_1,...,\delta_k):\type}$$ for all $k\in\{0,...,n-1\}$ for some $n\geq 0$, and $t$ be the interpretation of $\judge{\gamma:\Gamma}{t(\gamma):D(\gamma)}$.
    
    Let $P_A$ be the interpretation of $\judge{\gamma:\Gamma,d:D(\gamma),\delta:\Delta(\gamma,d)}{A(\gamma,d,\delta):\type}$, and $a$ be the interpretation of $\judge{\gamma:\Gamma,d:D(\gamma),\delta:\Delta(\gamma,d)}{a(\gamma,d,\delta):A(\gamma,d,\delta)}$.
    
    Then the interpretations of: \begin{center}$\judge{\gamma:\Gamma,\delta:\Delta(\gamma,t)}{A(\gamma,t,\delta):\type}$\;\;\;and\;\;\;$\judge{\gamma:\Gamma,\delta:\Delta(\gamma,t)}{a(\gamma,t,\delta):A(\gamma,t,\delta)}$\end{center} are defined and are $P_{A[t.\Delta]}$ and $a[t.\Delta]$ respectively.

    \item \textit{General substitution property.} Let $\Gamma$ be the interpretation of $\judgectx{\gamma:\Gamma}$, $P_{\Delta_{k+1}}$ be the interpretation of: $$\judge{\gamma:\Gamma,\delta_1:\Delta_1(\gamma),...,\delta_k(\gamma,\delta_1,...,\delta_{k-1})}{\Delta_{k+1}(\gamma,\delta_1,...,\delta_k):\type}$$ for all $k\in\{0,...,n-1\}$ for some $n\geq 0$, let $\Omega\xrightarrow{f}\Gamma$ be the interpretation of the generalised term judgement $\judge{\omega:\Omega}{f(\omega):\Gamma}$.
    
    Let $P_A$ be the interpretation of $\judge{\gamma:\Gamma,\delta:\Delta(\gamma)}{A(\gamma,\delta):\type}$, and $a$ be the interpretation of $\judge{\gamma:\Gamma,\delta:\Delta(\gamma)}{a(\gamma,\delta):A(\gamma,\delta)}$.
    
    Then the interpretations of: \begin{center}$\judge{\omega:\Omega,\delta:\Delta(f(\omega))}{A(f(\omega),\delta):\type}$\;\;\;and\;\;\;$\judge{\omega:\Omega,\delta:\Delta(f(\omega))}{a(f(\omega),\delta):A(f(\omega),\delta)}$\end{center} are defined and are $P_{A[f.\Delta]}$ and $a[f.\Delta]$ respectively.
\end{itemize} The proof of these properties is by induction on the derivation of type and term judgements. We show the weakening property for the second clause for type judgement interpretation and for the second, the eleventh, and the twelfth clauses for term judgement interpretation (i.e. the clauses involving propositional dependent sums). For the other clauses the argument is analogous. Let the interpretations of: \begin{center}$\judgectx{\gamma:\Gamma}$\\$\judge{\gamma}{A(\gamma):\type}$\\$\judge{\gamma,x:A(\gamma)}{B(\gamma,x):\type}$\\$\judge{\gamma,u}{C(\gamma,u):\type}$\\$\judge{\gamma,x,y}{c(\gamma,x,y):C(\gamma,\langle x,y\rangle)}$\\$\judge{\gamma}{D(\gamma):\type}$\end{center} be $\Gamma$, $P_A$, $P_B$, $P_C$, $c$, and $P_D$ respectively. By inductive hypothesis, the interpretations of $\judge{\gamma,d:D(\gamma)}{A(\gamma):\type}$, $\judge{\gamma,d:D(\gamma),x:A(\gamma)}{B(\gamma,x):\type}$, $\judge{\gamma,d:D(\gamma),u}{C(\gamma,u):\type}$, $\judge{\gamma,d:D(\gamma),x,y}{c(\gamma,x,y):C(\gamma,\langle x,y\rangle)}$ exist and are: \begin{center}$P_{A[P_D]}$,\;\;\;$P_{B[P_D^\lilbullet]}$,\;\;\;$P_{C[P_D^\lilbullet]}$,\;\;\;and\;\;\;$c[P_D^{\lilbullet\lilbullet}]$\end{center} respectively. Then, by definition, the interpretations of: \begin{center} $\judge{\gamma, d:D(\gamma)}{\Sigma_{x:A(\gamma)}B(\gamma,x):\type}$\\$\judge{\gamma,d:D(\gamma),x,y}{\langle x,y\rangle:\Sigma_{x:A(\gamma)}B(\gamma,x)}$\\
$\judge{\gamma:\Gamma,d:D(\gamma),u:\Sigma_{x:A(\gamma)}B(\gamma,x)}{\splitt(c,\gamma,u):C(\gamma,u)}$\\$\judge{\gamma,d:D(\gamma),x,y}{\sigma(c,\gamma,x,y):\splitt(c,\gamma,\langle x,y\rangle)=c(\gamma,x,y)}$ \end{center} exist and are: \begin{center}$P_{\sigmad_{A[P_D]}^{B[P_D^\lilbullet]}}$,\;\;\;$\pairing_{A[P_D]}^{B[P_D^\lilbullet]}$,\;\;\;$\splitt_{c[P_D^{\lilbullet\lilbullet}]}$,\;\;\;and\;\;\;$\sigma_{c[P_D^{\lilbullet\lilbullet}]}$\end{center} respectively. By the \textit{stability conditions} of \Cref{def:Semantics of axiomatic Sigma-types} these morphisms coincide with: \begin{center}$P_{\sigmad_A^B[P_D]}$,\;\;\;$\pairing_A^B[P_D^{\lilbullet\lilbullet}]$,\;\;\;$\splitt_{c}[P_D^\lilbullet]$,\;\;\;and\;\;\;$\sigma_c[P_D^{\lilbullet\lilbullet}]$\end{center} respectively, hence we are done. The proof of the substitution property and the general substitution property are analogous.

These three results are used in the proof by induction of the four points of the statement. 
In the theory ATT, there are fewer rules ending with a type/term equality judgement compared to ITT, hence the proof of the third and the fourth points is contained in the proof of the soundness of the interpretation of ITT, which is contained e.g. in \cite{Hofmann1997}. Analogously, in ATT the rules ending with a type judgement are the same as in ITT, hence we are done also with the first point of the statement. Regarding the second point, there are additional clauses for term judgement interpretation that we need to check (the others are contained in ATT), namely the computation axioms as well as the expansion axiom for function extensionality. Let us check e.g. the computation axiom for $\Sigma$-types, as the others are analogous. Let $\Gamma$ be the interpretation of the context $\judgectx{\gamma:\Gamma}$. By inductive hypothesis, the interpretations of: \begin{center}$\judge{\gamma:\Gamma}{A(\gamma):\type}$\\$\judge{\gamma:\Gamma,x:A(\gamma)}{B(\gamma,x):\type}$\\$\judge{\gamma,u}{C(\gamma,u):\type}$\\$\judge{\gamma,x,y}{c(\gamma,x,y):C(\gamma,\langle x,y\rangle)}$\end{center} are defined. Let us call them $P_A$, $P_B$, $P_C$, and $c$ respectively. Then the section: $$\sigma_c:\Gamma.A.B\to\Gamma.A.B.\id_{C[\pair_A^B]}[\;\splitt_c[\pair_A^B]\;;\;c\;]$$ of $P_{\id_{C[\pair_A^B]}[\;\splitt_c[\pair_A^B]\;;\;c\;]}$ is defined. Hence the interpretation of: $$\judge{\gamma,x,y}{\sigma(c,\gamma,x,y):\splitt(c,\gamma,\langle x,y\rangle)=c(\gamma,x,y)}$$ is defined. \endproof
\end{theorem} In light of this result, we may call such a display map category $(\mathbf{C},\mathcal{D})$ equipped with axiomatic $=$-types, axiomatic $\Sigma$-types, axiomatic $\Pi$-types \& axiomatic function extensionality, and axiomatic $0$-, $1$-, $2$-, $\nat$-types a \textit{model of }ATT. As anticipated, in the next section we show how to provide a categorical presentation of axiomatic type constructors.

\section{Categorical formulation of the semantics of axiomatic type theory}\label{section:Categorical formulation of the semantics of axiomatic dependent type theory}

As mentioned in the introduction, the issue with the formulation presented in the Section \ref{section:Syntactic approach to the semantics of axiomatic dependent type theory} regarding the semantics of axiomatic type formers is that, in this form, it may be regarded as impractical, depending on the context---for example, when seeking concrete models of ATT. As we pointed out, a more practical formulation should rely less on equipping the base structure with choice functions and more on universal properties, which allow these choice functions to be derived automatically. Universal properties are generally easier to identify in mathematical practice. However, the challenge with axiomatic type formers is that one-dimensional universal properties fail to uniquely characterise them. These properties collapse intensional and axiomatic type formers into extensional ones, making it impossible to distinguish between them. Consequently, they effectively characterise only the semantic transcription of the rules of ETT.

This is why the notion we introduce below is 2-dimensional: while 1-cells, as usual, represent substitutions, 2-cells represent (equivalence classes of) \textit{context propositional equalities}, or propositional equalities between substitutions: lists, indicated as $\judge{\gamma:\Gamma}{p(\gamma):f(\gamma)=g(\gamma)}$, of term judgements of the form: 
\begin{itemize}
    \item[] $\judge{\gamma:\Gamma}{p_1(\gamma):f_1(\gamma)=g_1(\gamma)}$
    \item[] $\judge{\gamma:\Gamma}{p_2(\gamma):f_2(\gamma)=p_1(\gamma)^*g_2(\gamma)}$
    \item[] $\judge{\gamma:\Gamma}{p_3(\gamma):f_3(\gamma)=(p_1(\gamma),p_2(\gamma))^*g_3(\gamma)}$
    \item[] \quad ...
    \item[] $\judge{\gamma:\Gamma}{p_n(\gamma):f_n(\gamma)=(p_1(\gamma),...,p_{n-1}(\gamma))^*g_n(\gamma)}$
\end{itemize} where $\judge{\gamma:\Gamma}{f(\gamma):\Delta}$ and $\judge{\gamma:\Gamma}{g(\gamma):\Delta}$ are parallel substitutions of ATT and trasport operations along multiple identity proofs are defined by the elimination rule for axiomatic identity types---see \Cref{sec:completeness} for more details.

Following the approach devised by Garner \cite{MR2525957}, the presence of an additional dimension enables us to impose more fine-grained universal properties on such a structure. These properties allow us to identify axiomatic type formers without conflating them with extensional ones.

\begin{definition}\label{split display map 2-category} A \textit{(split) display map 2-category} $(\boldsymbol{\mathcal{C}},\mathcal{D})$ is a (2,1)-category $\boldsymbol{\mathcal{C}}$ with a chosen 2-terminal object $\epsilon$, together with a class $\mathcal{D}$ of 1-cells of $\boldsymbol{\mathcal{C}}$, called \textit{display maps} and denoted as $\Gamma.A\to \Gamma$ and labeled as $P_A$, such that:
    \begin{itemize}
        \item[$\succ$] for every display map $\Gamma.A\to \Gamma$ and every arrow $f: \Delta \to \Gamma$, there is a choice of a display map $\Delta.A[f]\to\Delta$ such that the square: 
    \[\begin{tikzcd}[column sep=tiny,row sep=small]
	{\Delta.A[f]} & {\Gamma.A} \\
	\Delta & \Gamma
	\arrow[from=1-1, to=1-2]
	\arrow[from=1-1, to=2-1]
	\arrow["\lrcorner"{anchor=center, pos=0.125}, draw=none, from=1-1, to=2-2]
	\arrow[from=1-2, to=2-2]
	\arrow["f"{description}, from=2-1, to=2-2]
    \end{tikzcd}\] is a pullback and a 2-pullback; as the former display map is labeled as $P_A$, the latter will be labeled as $P_{A[f]}$; the arrow $\Delta.A[f]\to \Gamma.A$ will be denoted as $f.A$, or as $f^\lilbullet$ in absence of ambiguity; we will adopt the notation $A^\liltriangle$ in place of $A[P_A]$;
    \item[$\succ$] every display map $\Gamma.A\to \Gamma$ is a \textit{cloven isofibration} i.e. for every 2-cell $p$ as in the diagram below, there is a choice of a 1-cell $\t_g^p$ and of a 2-cell $\tau_g^p$ as in the diagram below, in such a way that the equality:
\[\begin{tikzcd}[sep=tiny]
	\Delta &&&& {\Gamma.A} && \Delta &&&& {\Gamma.A} \\
	&& {} & {} & {} & {=} &&& {} & {} & {} \\
	&&&& \Gamma &&&&&& \Gamma
	\arrow["g", shift left=2, from=1-1, to=1-5]
	\arrow["f"', curve={height=12pt}, from=1-1, to=3-5]
	\arrow[from=1-5, to=3-5]
	\arrow[""{name=0, anchor=center, inner sep=0}, "g", shift left=2, from=1-7, to=1-11]
	\arrow[""{name=1, anchor=center, inner sep=0}, "{\t_g^p}"', shift right, curve={height=6pt}, from=1-7, to=1-11]
	\arrow["f"', curve={height=12pt}, from=1-7, to=3-11]
	\arrow[from=1-11, to=3-11]
	\arrow["p"{pos=0.4}, shorten <=5pt, shorten >=10pt, Rightarrow, from=2-3, to=2-5]
	\arrow[shift right=2, shorten <=14pt, shorten >=10pt, Rightarrow, no head, from=2-9, to=2-11]
	\arrow["{\;\tau_g^p}"', shorten <=2pt, shorten >=2pt, Rightarrow, from=1, to=0]
\end{tikzcd}\] between the 2-cells $p$ and $P_A*\tau_g^p$ holds and the equalities: \begin{center} $\t_{gh}^{p*h}=\t_g^ph$ \quad\quad and \quad\quad $\tau_{gh}^{p*h}=\tau_g^p*h$ \end{center} hold for every 1-cell $h:\Omega\to\Delta$;

    \item[$\succ$] the equalities: $$\begin{alignedat}{2}
        P_{A[1_A]}&=P_A& (1_\Gamma).A&=1_{\Gamma.A}\\
        P_{A[fg]}&=P_{A[f][g]}&\quad\quad (f.A)(g.A[f])&=(fg).A
    \end{alignedat}$$ hold for every choice of composable arrows $f$ and $g$ and every display map $\Gamma.A\to \Gamma$, where $\Gamma$ is the target of $f$;

    \item[$\succ$] for every choice of display maps: \begin{center} $\Gamma.A_1\to\Gamma$,\;\; $\Gamma.A_1.A_2\to\Gamma.A_1$,\;\; ...,\;\; $\Gamma.A_1.\,...\,.A_{n-1}.A_n\to \Gamma.A_1.\,...\,.A_{n-1}$\end{center} and for every 1-cell $\Delta \xrightarrow{f} \Gamma$,
    \footnotetext[2]{Here $A_1.\,...\,.A_n$ and $A_1[f].A_2[f^\lilbullet].A_3[f^{\lilbullet\lilbullet}].\,...\,.A_n[f^{\lilbullet\lilbullet\,...\,\lilbullet}]$ are indicated as $A$ and $A[f]$  respectively, and $P_{A_n}P_{A_{n-1}}...P_{A_1}$ is indicated as $P_A$. Re-indexings on 1- and 2-cells are defined using the 1- and 2-pullback structure, analogously to display map categories.}if we are given a display map $\Gamma.A.C\to\Gamma.A$ \footnotemark[2] and commutative squares: 
\[\begin{tikzcd}[sep=small]
	{\Omega'} && \Omega & {\Omega'} && \Omega \\
	\\
	{\Delta.A[f].C[f.A]} && {\Gamma.A.C} & {\Delta.A[f]} && {\Gamma.A}
	\arrow["{f'}", from=1-1, to=1-3]
	\arrow["{g'}"', from=1-1, to=3-1]
	\arrow["g", from=1-3, to=3-3]
	\arrow["{f'}", from=1-4, to=1-6]
	\arrow["{h'}"', from=1-4, to=3-4]
	\arrow["h", from=1-6, to=3-6]
	\arrow["{f.A.C}", from=3-1, to=3-3]
	\arrow["{f.A}", from=3-4, to=3-6]
\end{tikzcd}\] such that $P_AP_Cg=P_Ah$ and $P_{A[f]}P_Cg'=P_{A[f]}h'$ and a 2-cell: 
\[\begin{tikzcd}[sep=tiny]
	\Omega &&&& {\Gamma.A.C} \\
	&& {} & {} & {} \\
	&&&& {\Gamma.A}
	\arrow["g", from=1-1, to=1-5]
	\arrow["h"', curve={height=12pt}, from=1-1, to=3-5]
	\arrow[from=1-5, to=3-5]
	\arrow["p", shift right, shorten <=10pt, shorten >=10pt, Rightarrow, from=2-3, to=2-5]
\end{tikzcd}\] such that $P_A * p = 1_{P_Ah}$, then the following diagram of 2-cells:
\[\begin{tikzcd}[column sep=scriptsize]
	{\Omega'} && \Omega \\
	\\
	{\Delta.A[f].C[f.A]} && {\Gamma.A.C}
	\arrow["{f'}", from=1-1, to=1-3]
	\arrow[""{name=0, anchor=center, inner sep=0}, "{\t_{g'}^{p[f]}}"', shift right, curve={height=12pt}, from=1-1, to=3-1]
	\arrow[""{name=1, anchor=center, inner sep=0}, "{g'}", shift left, curve={height=-12pt}, from=1-1, to=3-1]
	\arrow[""{name=2, anchor=center, inner sep=0}, "{\t_{g}^{p}}"', shift right, curve={height=12pt}, from=1-3, to=3-3]
	\arrow[""{name=3, anchor=center, inner sep=0}, "g", shift left, curve={height=-12pt}, from=1-3, to=3-3]
	\arrow["{f.A.C}", from=3-1, to=3-3]
	\arrow["{\tau_{g'}^{p[f]}}", shift right, shorten <=6pt, shorten >=6pt, Rightarrow, from=0, to=1]
	\arrow["{\tau_{g}^{p}}", shift right, shorten <=6pt, shorten >=6pt, Rightarrow, from=2, to=3]
\end{tikzcd}\] commutes, i.e. the equalities: $$\t_g^p[f.A]=\t_{g'}^{p[f]}\quad\quad\text{and}\quad\quad\tau_g^p[f.A]=\tau_{g'}^{p[f]}$$ hold.\footnotemark[2]

\end{itemize}
\end{definition}

The notion of a display map 2-category captures the structural aspects of dependent type theory, combined with fragments of the specific logical aspects of ATT. Specifically, the cloven isofibration structure with which the display maps are equipped provides a 2-dimensional interpretation of a consequence of the elimination rule and the computation \textit{axiom} for $=$-types: transport---along one or more identity proofs. In \cite{MR2525957}, display maps in a model of ITT are endowed with a \textit{normal} isofibration structure, which in our setting corresponds to the additional \textit{normality} requirements that $\t_g^{1_{P_Ag}}=g$ and $\tau_g^{1_{P_Ag}}=1_g$. However, dropping these requirements and reducing display maps to mere cloven isofibrations is fundamental to prevent models of ATT from being models of ITT. Let us briefly analyse how the cloven isofibration structure is induced in the context of the syntax of ATT itself.

If we are given a type judgement $\judge{\gamma:\Gamma}{A(\gamma):\type}$ and substitutions: \begin{center}$\judge{\delta:\Delta}{f(\delta):\Gamma}$ \quad and \quad $\judge{\delta:\Delta}{g(\delta):\Gamma.A}$\end{center}---hence $g(\delta)$ is of the form $g_1(\delta):\Gamma\,,\,g_2(\delta):A(g_1(\delta))$---with a context identity proof $\judge{\delta:\Delta}{p(\delta):f(\delta)=g_1(\delta)}$, then $\t_g^p$ is the substitution $\judge{\delta:\Delta}{f(\delta):\Gamma\,,\,p(\delta)^*g_2(\delta):A(f(\delta))}$ and $\tau_g^p$ is the context identity proof: $$\judge{\delta:\Delta}{f(\delta)\,,\,p(\delta)^*g_2(\delta) \quad = \quad g_1(\delta)\,,\,g_2(\delta)}$$ provided by the list: $$\begin{alignedat}{3}
    \judge{\delta:\Delta}& &{p(\delta):\;\;}& &f(\delta)&=g_1(\delta)\\
    \judge{\delta:\Delta}& &{\refl(p(\delta)^*g_2(\delta)):\;\;}& &p(\delta)^*g_2(\delta)&=p(\delta)^*g_2(\delta)
\end{alignedat}$$ of identity proofs. Now, if $f(\delta)\equiv g_1(\delta)$ and $p(\delta)\equiv \refl(g_1(\delta))$, then $\t_g^p$ is $\judge{\delta:\Delta}{g_1(\delta):\Gamma\,,\,\refl(g_1(\delta))^*g_2(\delta):A(g_1(\delta))}$ and $\tau_g^p$ is the list: $$\begin{alignedat}{3}
    \judge{\delta:\Delta}&&{\refl(g_1(\delta)):\;\;}&&g_1(\delta)&=g_1(\delta)\\
    \judge{\delta:\Delta}&&{\refl(\refl(g_1(\delta))^*g_2(\delta)):\;\;}&&\refl(g_1(\delta))^*g_2(\delta)&=\refl(g_1(\delta))^*g_2(\delta)
\end{alignedat}$$ hence in general in this case $\t_g^p$ is \textit{not} $g(\delta)$ and $\tau_g^p$ is \textit{not} its identity 2-cell: in ATT we can infer that $\refl(g_1(\delta))^*g_2(\delta)=g_2(\delta)$---it is a fragment of the computation axiom for $=$-types---but not that $\refl(g_1(\delta))^*g_2(\delta)\equiv g_2(\delta)$. In other words, the display map associated to the type $A$ is a cloven isofibration but not necessarily a normal isofibration.

This fact explains why, in a general display map 2-category, display maps are cloven isofibrations but not normal isofibrations: this is what happens when we consider the syntax of ATT. As we will see---Subsection \ref{subsection:The role of the cloven isofibration structure}---the cloven isofibration structure on display maps enables the interpretation of computation axioms within a display map 2-category. However, the fact that they are generally \textit{not} normal isofibrations is crucial to ensuring that these computation axioms are validated without necessarily validating the corresponding computation rules: this weakening prevents the interpretation of terms $\h_c$ and $\sigma_c$ from collapsing into reflexivities, hence the respective propositional equalities from being judgemental.

\subsection{Axiomatic $\boldsymbol{=}$-types}

We specialise the notion of a display map 2-category to obtain a version that incorporates a 2-dimensional semantic interpretation of axiomatic $=$-types.

\begin{definition}[Semantics of axiomatic $=$-types---categorical formulation]\label{def:Categorical semantics of axiomatic =-types}
    Let $(\boldsymbol{\mathcal{C}},\mathcal{D})$ be a display map 2-category.

    Let us assume that, for every object $\Gamma$ and every display map $P_A$ of codomain $\Gamma$, there is a choice of a display map $\Gamma.A.A^\liltriangle.\id_A\to\Gamma.A.A^\liltriangle$ and of an \textit{arrow object}: 
\[\begin{tikzcd}[sep=scriptsize]
	{\Gamma.A.A^\liltriangle.\id_A} &&& {\Gamma.A} \\
	&&& \Gamma
	\arrow[""{name=0, anchor=center, inner sep=0}, "{P_A^\lilbullet P_{\id_A}}"', shift right=2, from=1-1, to=1-4]
	\arrow[""{name=1, anchor=center, inner sep=0}, "{P_{A^\liltriangle}P_{\id_A}}", shift left=2, from=1-1, to=1-4]
	\arrow["{P_AP_{A^\liltriangle}P_{\id_A}}"', curve={height=12pt}, from=1-1, to=2-4]
	\arrow["{P_A}", from=1-4, to=2-4]
	\arrow["{\;\alpha_A}", shorten <=1pt, shorten >=1pt, Rightarrow, from=1, to=0]
\end{tikzcd}\]
for $P_A$, i.e. a 2-cell $\alpha_A$ of $\boldsymbol{\mathcal{C}}/\Gamma$ into $P_A$ as in the diagram that induces by post-composition a family of isomorphisms of categories: $$(\boldsymbol{\mathcal{C}} / \Gamma)(\;h\;,\;P_{A}P_{A^\liltriangle}P_{\id_A}\;)\;\;\to\;\;(\boldsymbol{\mathcal{C}} / \Gamma)(\;h\;,\;P_{A}\;)\;^\to$$ 2-natural in $h$, where $h$ is any 1-cell $\Delta\to\Gamma$, in such a way that the \textit{stability conditions}: \begin{center}
    $\id_A[f^{\lilbullet\lilbullet}]=\id_{A[f]}$ \quad\quad and \quad\quad $\alpha_A[f]=\alpha_{A[f]}$
\end{center} hold for every 1-cell $f:\Delta\to \Gamma$.

Then we say that $(\boldsymbol{\mathcal{C}},\mathcal{D})$ is \textit{endowed with axiomatic $=$-types}.
\end{definition}

\noindent\textbf{Notation and properties.} Let $a,b$ be sections of a display map $P_A:\Gamma.A\to\Gamma$ and let $p$ be an arrow $a\Rightarrow b$ in the category $(\boldsymbol{\mathcal{C}}/\Gamma)(1_\Gamma,P_A)$ of sections of $P_A$. Thenm being $\alpha_A$ an arrow object, there is a unique arrow: $$\tilde p : \Gamma\to\Gamma.A.A^\liltriangle.\id_A$$ in $(\boldsymbol{\mathcal{C}}/\Gamma)(1_\Gamma,P_{A}P_{A^\liltriangle}P_{\id_A})$ such that $\alpha_A * \tilde p = p$. In particular: $$P_{A^\liltriangle}P_{\id_A}\tilde p=a \quad \text{and} \quad P_{A}^\lilbullet P_{\id_A}\tilde p=b$$ which means that $P_{\id_A}\tilde p =a;b$ i.e. the square: 
\[\begin{tikzcd}[column sep=small,row sep=scriptsize]
	\Gamma && {\Gamma.A.A^\liltriangle.\id_A} \\
	\\
	\Gamma && {\Gamma.A.A^\liltriangle}
	\arrow["{\tilde p}", from=1-1, to=1-3]
	\arrow[equals, from=1-1, to=3-1]
	\arrow[from=1-3, to=3-3]
	\arrow["{a;b}", from=3-1, to=3-3]
\end{tikzcd}\] commutes. We denote as $\{p\}$ the unique section of $P_{\id_A[a;b]}$ such that: 
\[\begin{tikzcd}[column sep=tiny,row sep=scriptsize]
	\Gamma \\
	\\
	&& {\Gamma.\id_A[a;b]} && {\Gamma.A.A^\liltriangle.\id_A} \\
	\\
	&& \Gamma && {\Gamma.A.A^\liltriangle}
	\arrow["{\{p\}}"{description}, from=1-1, to=3-3]
	\arrow["{\tilde p}", curve={height=-12pt}, from=1-1, to=3-5]
	\arrow[curve={height=12pt}, equals, from=1-1, to=5-3]
	\arrow[from=3-3, to=3-5]
	\arrow[from=3-3, to=5-3]
	\arrow[from=3-5, to=5-5]
	\arrow["{a;b}", from=5-3, to=5-5]
\end{tikzcd}\] commutes i.e. the unique section of $P_{\id_A[a;b]}$ such that: \begin{equation}\label{whatcharacterisescurlyparentheses}
    \alpha_A(a;b)^\lilbullet\{p\}=p \quad.
\end{equation} Now, let $f$ be any arrow $\Delta \to \Gamma$ and let us observe that the diagram: 
\[\begin{tikzcd}
	\Delta && \Gamma \\
	{\Delta.\id_{A[f]}[a[f];b[f]]} && {\Gamma.\id_A[a;b]} \\
	{\Delta.A[f].A[f]^\liltriangle.\id_{A[f]}} && {\Gamma.A.A^\liltriangle.\id_A} \\
	{\Delta.A[f]} && {\Gamma.A}
	\arrow["f"{description}, from=1-1, to=1-3]
	\arrow["{\{p\}[f]}"', from=1-1, to=2-1]
	\arrow["{\{p\}}", from=1-3, to=2-3]
	\arrow["{f^\lilbullet}"{description}, from=2-1, to=2-3]
	\arrow["{(a[f];b[f])^\lilbullet}"', from=2-1, to=3-1]
	\arrow["{(a;b)^\lilbullet}", from=2-3, to=3-3]
	\arrow["{f^{\lilbullet\lilbullet\lilbullet}}"{description}, from=3-1, to=3-3]
	\arrow[""{name=0, anchor=center, inner sep=0}, shift right=5, from=3-1, to=4-1]
	\arrow[""{name=1, anchor=center, inner sep=0}, shift left=5, from=3-1, to=4-1]
	\arrow[""{name=2, anchor=center, inner sep=0}, shift right=3, from=3-3, to=4-3]
	\arrow[""{name=3, anchor=center, inner sep=0}, shift left=3, from=3-3, to=4-3]
	\arrow["{f^\lilbullet}"{description}, from=4-1, to=4-3]
	\arrow["{\alpha_A[f]}", shorten <=4pt, shorten >=4pt, Rightarrow, from=0, to=1]
	\arrow["{\alpha_A}", shorten <=2pt, shorten >=2pt, Rightarrow, from=2, to=3]
\end{tikzcd}\] commutes: the upper and lower squares commute by definition, while the middle one commutes because Diagram (\ref{reindexingprimareindexingdopo}) does. By (\ref{whatcharacterisescurlyparentheses}) and since $\alpha_A[f]$ is $\alpha_{A[f]}$, this implies that: \begin{equation}\label{curlyparenthesespreservereindexing}
    \{p\}[f]=\{p[f]\}
\end{equation} for every arrow $p$ in the category of the section of $P_A$ and every substitution $f$ for the context of $P_A$.

\medskip

Let $(\boldsymbol{\mathcal{C}},\mathcal{D})$ be a display map 2-category. Clearly $(\boldsymbol{\mathcal{C}},\mathcal{D})$ is in particular a display map category $(\mathbf{C},\mathcal{D})$, where the category $\mathbf{C}$ is the underlying category of $\boldsymbol{\mathcal{C}}$: the class $\mathcal{D}$ of display maps over $\boldsymbol{\mathcal{C}}$ according to Definition \ref{split display map 2-category} is in particular a class of display maps over $\mathbf{C}$ according to Definition \ref{display map category}.

Now, let us assume that the display map 2-category $(\boldsymbol{\mathcal{C}},\mathcal{D})$ is endowed with axiomatic $=$-types---Definition \ref{def:Categorical semantics of axiomatic =-types}---and let $P_A$ be a display map of codomain $\Gamma$. The next paragraphs show that the display map category $(\mathbf{C},\mathcal{D})$ induced by $(\boldsymbol{\mathcal{C}},\mathcal{D})$ is endowed with appropriate choice functions validating axiomatic $=$-types---as in Definition \ref{def:Semantics of axiomatic =-types}.

\medskip

\noindent\textbf{Form Rule for $\boldsymbol{=}$-types.} We already have a choice of a display map $\Gamma.A.A^\liltriangle.\id_A\to \Gamma.A.A^\liltriangle$.

\medskip

\noindent\textbf{Intro Rule for $\boldsymbol{=}$-types.} Being $\alpha_A$ an arrow object, the identity 2-cell:
\[\begin{tikzcd}[row sep=scriptsize]
	{\Gamma.A} && {\Gamma.A} \\
	&& \Gamma
	\arrow[""{name=0, anchor=center, inner sep=0}, "{1_{\Gamma.A}}"', shift right=3, from=1-1, to=1-3]
	\arrow[""{name=1, anchor=center, inner sep=0}, "{1_{\Gamma.A}}", shift left=2, from=1-1, to=1-3]
	\arrow["{P_A}"', curve={height=12pt}, from=1-1, to=2-3]
	\arrow["{P_A}", from=1-3, to=2-3]
	\arrow["{\;1_{1_{\Gamma.A}}}", shift right=2, shorten <=1pt, shorten >=1pt, Rightarrow, from=1, to=0]
\end{tikzcd}\] of $1_{\Gamma.A}$ in $\boldsymbol{\mathcal{C}}/\Gamma$ factors uniquely through $\alpha_A$ as a 1-cell: 
\[\begin{tikzcd}[sep=scriptsize]
	{\Gamma.A} && {\Gamma.A.A^\liltriangle.\id_A} \\
	&& \Gamma
	\arrow["{\refl_A}", from=1-1, to=1-3]
	\arrow["{P_A}"', curve={height=12pt}, from=1-1, to=2-3]
	\arrow["{P_{\id_A}P_{A^\liltriangle}P_A}", from=1-3, to=2-3]
\end{tikzcd}\] of $\boldsymbol{\mathcal{C}}/\Gamma$ such that $\alpha_A*\refl_A = 1_{1_{\Gamma.A}}$. In particular: $$
    P_{A^\liltriangle}(P_{\id_A}\refl_A)=1_{\Gamma.A} \quad\quad P_A^\lilbullet(P_{\id_A}\refl_A)=1_{\Gamma.A}$$ hence $P_{\id_A}\refl_A=\delta_A$---being: 
\[\begin{tikzcd}[sep=small]
	{\Gamma.A.A^\liltriangle} && {\Gamma.A} \\
	\\
	{\Gamma.A} && \Gamma
	\arrow["{P_A^\lilbullet}", from=1-1, to=1-3]
	\arrow["{P_{A^\liltriangle}}"', from=1-1, to=3-1]
	\arrow["\lrcorner"{anchor=center, pos=0.125}, draw=none, from=1-1, to=3-3]
	\arrow["{P_A}", from=1-3, to=3-3]
	\arrow["{P_A}", from=3-1, to=3-3]
\end{tikzcd}\] a pullback square---and therefore there is unique a section: $$\reflex_A:\Gamma.A\to\Gamma.A.\id_A[\delta_A]$$ of $P_{\id_A[\delta_A]}$ such that $\delta_A^\lilbullet\reflex_A=\refl_A$.

\medskip

\noindent\textbf{Elim Rule for $\boldsymbol{=}$-types.} Let $\Gamma.A.A^\liltriangle.\id_A.C\to\Gamma.A.A^\liltriangle.\id_A$ be a display map and let $c$ be a section $\Gamma.A\to\Gamma.A.C[\refl_A]$ of the display map $P_{C[\refl_A]}$. The pair $(\alpha_A,1_{P_A^\lilbullet P_{\id_A}})$ constitutes an arrow from $\alpha_A$ to $\alpha_A * \refl_A(P_A^\lilbullet P_{\id_A})$ in the category $(\boldsymbol{\mathcal{C}}/\Gamma)(P_AP_{A^\liltriangle}P_{\id_A},P_A)^\to$. Hence, it is induced, by post-composition via $\alpha_A$, by a unique arrow: $$\varphi_A:1_{\Gamma.A.A^\liltriangle.\id_A}\Longrightarrow\refl_A(P_A^\lilbullet P_{\id_A})$$ of the category $(\boldsymbol{\mathcal{C}}/\Gamma)(P_AP_{A^\liltriangle}P_{\id_A},P_AP_{A^\liltriangle}P_{\id_A})$. In other words, the equalities: \begin{center}
    $P_{A^\liltriangle}P_{\id_A}*\varphi_A=\alpha_A$ \quad and \quad $P_A^\lilbullet P_{\id_A}*\varphi_A=1_{P_A^\lilbullet P_{\id_A}}$
\end{center} hold. Moreover: \begin{center}
    $P_{A^\liltriangle}P_{\id_A}*\varphi_A*\refl_A = \alpha_A*\refl_A =1_{1_{\Gamma.A}}=P_{A^\liltriangle}P_{\id_A}*1_{\refl_A}$ \\
    $P_A^\lilbullet P_{\id_A}*\varphi_A*\refl_A =1_{P_A^\lilbullet P_{\id_A}}*\refl_A=1_{1_{\Gamma.A}}=P_A^\lilbullet P_{\id_A}*1_{\refl_A}$
\end{center} therefore $\varphi_A*\refl_A$ and $1_{\refl_A}$ induce, by post-composition via $\alpha_A$, the same arrow $\alpha_A*\refl_A \to \alpha_A*\refl_A$, i.e. $1_{1_{\Gamma.A}}\to 1_{1_{\Gamma.A}}$, in $(\boldsymbol{\mathcal{C}}/\Gamma)(P_A,P_A)^\to$ --- namely the pair $(1_{1_{\Gamma.A}},1_{1_{\Gamma.A}})$. By the universality of $\alpha_A$, we conclude that: $$\varphi_A * \refl_A = 1_{\refl_A}.$$ Now, let $\tilde\j_c:=(\refl_A.C)cP_A^\lilbullet P_{\id_A}$. We observe that: $$P_C\tilde\j_c=\refl_AP_{C[\refl_A]}cP_A^\lilbullet P_{\id_A}=\refl_AP_A^\lilbullet P_{\id_A}$$ hence we can rewrite the codomain of the 2-cell $\varphi_A$ as follows: 
\[\begin{tikzcd}[column sep=large,row sep=small]
	{\Gamma.A.A^\liltriangle.\id_A} && {\Gamma.A.A^\liltriangle.\id_A.C} \\
	&& {\Gamma.A.A^\liltriangle.\id_A}
	\arrow["{\tilde \j_c}", from=1-1, to=1-3]
	\arrow["{\varphi_A\,\,\,}"'{pos=0.9}, shift right=5, shorten <=55pt, Rightarrow, from=1-1, to=1-3]
	\arrow[curve={height=12pt}, equals, from=1-1, to=2-3]
	\arrow["{P_C}", from=1-3, to=2-3]
\end{tikzcd}\] and, using the cloven isofibration structure on $P_C$, we obtain a section $\t_{\tilde\j_c}^{\varphi_A}:\Gamma.A.A^\liltriangle.\id_A\to\Gamma.A.A^\liltriangle.\id_A.C$ of $P_C$, as well as a 2-cell:
\[\begin{tikzcd}[sep=scriptsize]
	{\Gamma.A.A^\liltriangle.\id_A} &&& {\Gamma.A.A^\liltriangle.\id_A.C}
	\arrow[""{name=0, anchor=center, inner sep=0}, "{\tilde\j_c}"', curve={height=18pt}, from=1-1, to=1-4]
	\arrow[""{name=1, anchor=center, inner sep=0}, "{\t_{\tilde\j_c}^{\varphi_A}}", curve={height=-18pt}, from=1-1, to=1-4]
	\arrow["{\;\tau_{\tilde\j_c}^{\varphi_A}}", shorten <=7pt, shorten >=7pt, Rightarrow, from=1, to=0]
\end{tikzcd}\] such that $P_C*\tau_{\tilde\j_c}^{\varphi_A}=\varphi_A$. We define $\j_c:=\t_{\tilde\j_c}^{\varphi_A}$.

\medskip

\noindent\textbf{Comp Axiom for $\boldsymbol{=}$-types.} Referring to the diagram: 
\begin{equation}\label{reindexingid}\begin{tikzcd}[column sep=scriptsize,row sep=small]
	{\Gamma.A} && {\Gamma.A.A^\liltriangle.\id_A} \\
	\\
	{\Gamma.A.C[\refl_A]} && {\Gamma.A.A^\liltriangle.\id_A.C} \\
	\\
	{\Gamma.A} && {\Gamma.A.A^\liltriangle.\id_A}
	\arrow["{\refl_A}"{description}, from=1-1, to=1-3]
	\arrow["{\j_c[\refl_A]}"', shift right=3, from=1-1, to=3-1]
	\arrow["c", shift left=3, from=1-1, to=3-1]
	\arrow[""{name=0, anchor=center, inner sep=0}, "{\tilde\j_c}", shift left=3, from=1-3, to=3-3]
	\arrow[""{name=1, anchor=center, inner sep=0}, "{\j_c}"', shift right=3, from=1-3, to=3-3]
	\arrow["{\refl_A^\lilbullet}"{description}, from=3-1, to=3-3]
	\arrow["{P_{C[\refl_A]}}"', from=3-1, to=5-1]
	\arrow["\lrcorner"{anchor=center, pos=0.125}, draw=none, from=3-1, to=5-3]
	\arrow["{P_C}", from=3-3, to=5-3]
	\arrow["{\refl_A}"{description}, from=5-1, to=5-3]
	\arrow[shorten <=2pt, shorten >=2pt, Rightarrow, from=1, to=0]
\end{tikzcd} \end{equation} and observing\footnotetext[3]{Since $P_C \tilde \j_c \refl_A=P_C\refl_A^\lilbullet cP_A^\lilbullet P_{\id_A}\refl_A=P_C\refl_A^\lilbullet c=\refl_AP_{C[\refl_A]}c=\refl_A1_{\Gamma.A}$ and $P_{C[\refl_A]}c=1_{\Gamma.A}$ and $\refl_A^\lilbullet c=\refl_A^\lilbullet cP_A^\lilbullet P_{\id_A}\refl_A=\j_c \refl_A$.} that $c=\tilde\j_c[\refl_A]$ \footnotemark[3] and that: $$P_C*\tau_{\tilde\j_c}^{\varphi_A}*\refl_A = P_C*\tau_{\tilde\j_c\refl_A}^{\varphi_A*\refl_A}=\varphi_A*\refl_A=1_{\refl_A}={\refl_A}*1_{1_{\Gamma.A}}$$ we conclude, by the 2-universal property of $\Gamma.A.C[\refl_A]$, that there is unique a 2-cell:
\[\begin{tikzcd}[sep=small]
	{\Gamma.A} &&&& {\Gamma.A}
	\arrow[""{name=0, anchor=center, inner sep=0}, "c"', curve={height=12pt}, from=1-1, to=1-5]
	\arrow[""{name=1, anchor=center, inner sep=0}, "{\j_c[\refl_A]}", curve={height=-12pt}, from=1-1, to=1-5]
	\arrow["{h_c}", shift right=2, shorten <=5pt, shorten >=5pt, Rightarrow, from=1, to=0]
\end{tikzcd}\hspace{-0.4em}{.C[\refl_A]}\] such that $P_{C[\refl_A]} *h_c=1_{1_{\Gamma.A}}$ and $\refl_A^\lilbullet*h_c=\tau_{\tilde\j_c}^{\varphi_A}*\refl_A$. Then we define $\h_C$ as the section: $$\{h_c\}:\Gamma.A\to\Gamma.A.\id_{C[\refl_A]}[\j_c[\refl_A];c]$$ satisfying (\ref{whatcharacterisescurlyparentheses}), which we recall to be defined as follows: by the universal property of $\alpha_{C[\refl_A]}$, there is unique an arrow $\tilde h_c:\Gamma.A \to \Gamma.A.C[\refl_A].C[\refl_A]^\liltriangle.\id_{C[\refl_A]}$ over $\Gamma.A$ such that $\alpha_{C[\refl_A]}*\tilde h_c=h_c.$ In particular: $$P_{C[\refl_A]^\liltriangle}P_{\id_{C[\refl_A]}}\tilde h_c=\j_c[\refl_A]\;\;\;\;\;\text{and}\;\;\;\;\;P_{C[\refl_A]}^\lilbullet P_{\id_{C[\refl_A]}}\tilde h_c=c$$ which means that $P_{\id_{C[\refl_A]}}\tilde h_c = \j_c[\refl_A];c$. Hence the diagram:
\[\begin{tikzcd}[column sep=tiny,row sep=small]
	{\Gamma.A} \\
	\\
	{\Gamma.A.\id_{C[\refl_A]}[\j_c[\refl_A];c]} && {\Gamma.A.C[\refl_A].C[\refl_A]^\liltriangle.\id_{C[\refl_A]}} \\
	\\
	{\Gamma.A} && {\Gamma.A.C[\refl_A].C[\refl_A]^\liltriangle}
	\arrow["{\h_c}"', from=1-1, to=3-1]
	\arrow["{\tilde h_c}", curve={height=-18pt}, from=1-1, to=3-3]
	\arrow[from=3-1, to=3-3]
	\arrow[from=3-1, to=5-1]
	\arrow[from=3-3, to=5-3]
	\arrow[""{name=0, anchor=center, inner sep=0}, "{\j_c[\refl_A];c}"', from=5-1, to=5-3]
	\arrow["\lrcorner"{anchor=center, pos=0.125}, draw=none, from=3-1, to=0]
\end{tikzcd}\] commutes for a unique section: $$\{h_c\}:\Gamma.A\to\Gamma.A.\id_{C[\refl_A]}[\j_c[\refl_A];c]$$ of $P_{\id_{C[\refl_A]}[\j_c[\refl_A];c]}$.

\begin{proposition}\label{stability of id in 2d model}
    Let $(\boldsymbol{\mathcal{C}},\mathcal{D})$ be a display map 2-category endowed with axiomatic $=$-types. Referring to the data defined in paragraphs \emph{\textbf{Form Rule, Intro Rule, Elim Rule, and Comp Axiom, for $\boldsymbol{=}$-types}}, the \emph{stability conditions}: $$\begin{alignedat}{2}
        \id_A[f^{\lilbullet\lilbullet}]&=\id_{A[f]}&\quad\quad\j_c[f^{\lilbullet\lilbullet\lilbullet}]&=\j_{c[f^\lilbullet]}\\
        \reflex_A[f^\lilbullet]&=\reflex_{A[f]}&\h_c[f^\lilbullet]&=\h_{c[f^\lilbullet]}\\
    \end{alignedat}$$ of Defintion \ref{def:Semantics of axiomatic =-types} hold for every arrow $\Delta \xrightarrow{f}\Gamma$.

\proof
See Appendix \ref{appendixIII}.
\endproof
\end{proposition}

\subsection{Axiomatic $\boldsymbol{\Sigma}$-types}

We specialise the notion of a display map 2-category with axiomatic $=$-types to obtain a version that incorporates a 2-dimensional semantic interpretation of axiomatic $\Sigma$-types. The condition we require to be satisfied is a generalisation of the condition characterising the semantics of extensional $\Sigma$-types: in place of asking that display maps are closed under composition up to
isomorphism \cite{MR1674451}, we require that they are closed by composition up to homotopy equivalence. As we
will see, the 2-cell components of such a homotopy equivalence induce the interpretation of the corresponding Comp Axiom. We observe a significant difference between the intensional case \cite{MR2525957} and the axiomatic one: such a homotopy equivalence need not consist of a deformation retract.

\begin{definition}[Semantics of axiomatic $\Sigma$-types---categorical formulation]\label{def:Categorical semantics of axiomatic Sigma-types}
    Let $(\boldsymbol{\mathcal{C}},\mathcal{D})$ be a display map 2-category endowed with axiomatic $=$-types.

    Let us assume that, for every object $\Gamma$, every display map $P_A$ of codomain $\Gamma$, and every display map $P_B$ of codomain $\Gamma.A$, there is a choice of a display map $\Gamma.\sigmad_A^B\to\Gamma$ and of an \textit{equivalence}:
\[\begin{tikzcd}[row sep=small]
	{\Gamma.A.B} && {\Gamma.\sigmad_A^B} \\
	\\
	{\Gamma.A} && \Gamma
	\arrow["{(\pair_A^B,\pi_A^B,\eta_A^B,\beta_A^B)}", from=1-1, to=1-3]
	\arrow["{P_B}"', from=1-1, to=3-1]
	\arrow["{P_{\sigmad_A^B}}", from=1-3, to=3-3]
	\arrow["{P_A}"{description}, from=3-1, to=3-3]
\end{tikzcd}\] from $P_AP_B$ to $P_{\sigmad_A^B}$, i.e. a choice of objects: $$\pair_A^B:\Gamma.A.B\to\Gamma.\sigmad_A^B\textnormal{\;\;\;\;and\;\;\;\;}\pi_A^B:\Gamma.\sigmad_A^B\to\Gamma.A.B$$ of $(\boldsymbol{\mathcal{C}}/\Gamma)(P_AP_B,P_{\sigmad_A^B})$ and $(\boldsymbol{\mathcal{C}}/\Gamma)(P_{\sigmad_A^B},P_AP_B)$ respectively and of arrows: $$\eta_A^B:1_{\Gamma.\sigmad_A^B}\Longrightarrow\pair_A^B\pi_A^B\textnormal{\;\;\;\;and\;\;\;\;}\beta_A^B:\pi_A^B\pair_A^B\Longrightarrow1_{\Gamma.A.B}$$ of $(\boldsymbol{\mathcal{C}}/\Gamma)(P_{\sigmad_A^B},P_{\sigmad_A^B})$ and $(\boldsymbol{\mathcal{C}}/\Gamma)(P_AP_B,P_AP_B)$ respectively, in such a way that the \textit{stability conditions}: $$\sigmad_A^B[f]=\sigmad_{A[f]}^{B[f^\lilbullet]}\quad\quad\begin{alignedat}{2}
        \pair_A^B[f]&=\pair_{A[f]}^{B[f^\lilbullet]}&\quad\quad\pi_A^B[f]&=\pi_{A[f]}^{B[f^\lilbullet]}\\
        \eta_A^B[f]&=\eta_{A[f]}^{B[f^\lilbullet]}&\beta_A^B[f]&=\beta_{A[f]}^{B[f^\lilbullet]}\\
    \end{alignedat}$$ hold for every 1-cell $f:\Delta\to\Gamma$.
    
Then we say that $(\boldsymbol{\mathcal{C}},\mathcal{D})$ is \textit{endowed with axiomatic $\Sigma$-types}.
\end{definition}

Let $(\boldsymbol{\mathcal{C}},\mathcal{D})$ be a display map 2-category endowed with axiomatic $=$-types.

Let us assume that $(\boldsymbol{\mathcal{C}},\mathcal{D})$ is endowed with axiomatic $\Sigma$-types---Definition \ref{def:Categorical semantics of axiomatic Sigma-types}---let $P_A$ be a display map of codomain $\Gamma$ and let $P_B$ be a display map of codomain $\Gamma.A$. The next paragraphs show that the display map category $(\mathbf{C},\mathcal{D})$ induced by $(\boldsymbol{\mathcal{C}},\mathcal{D})$ is endowed with appropriate choice functions validatin axiomatic $\Sigma$-types---as in Definition \ref{def:Semantics of axiomatic Sigma-types}.

\medskip

\noindent\textbf{Form Rule for $\boldsymbol{\Sigma}$-types.} We already have a choice of a display map $\Gamma.\sigmad_A^B\to \Gamma$.

\medskip

\noindent\textbf{Intro Rule for $\boldsymbol{\Sigma}$-types.} Since $P_{\sigmad_A^B}\pair_A^B=P_AP_B$ there is unique a section: $$\pairing_A^B:\Gamma.A.B\to\Gamma.A.B.\sigmad_A^B[P_AP_B]$$ of $P_{\sigmad_A^B[P_AP_B]}$ such that $(P_AP_B)^\lilbullet\reflex_A=\pair_A^B$.

\medskip

\noindent\textbf{Elim Rule for $\boldsymbol{\Sigma}$-types.} Let $\Gamma.\sigmad_A^B.C\to\Gamma.\sigmad_A^B$ be a display map and let $c$ be a section $\Gamma.A.B\to\Gamma.A.B.C[\pair_A^B]$ of the display map $P_{C[\pair_A^B]}$. Let $\tilde\splitt_c:=(\pair_A^B)^\lilbullet c\pi_A^B$. We observe that: $$P_C\tilde\splitt_c=\pair_A^B P_{C[\pair_A^B]}c\pi_A^B=\pair_A^B\pi_A^B$$ hence we can rewrite the codomain of the 2-cell $\eta_A^B$ as follows: 
\[\begin{tikzcd}[cramped,column sep=large,row sep=small]
	{\Gamma.\sigmad_A^B} && {\Gamma.\sigmad_A^B.C} \\
	&& {\Gamma.\sigmad_A^B}
	\arrow["{\tilde \splitt_c}", from=1-1, to=1-3]
	\arrow["{\eta_A^B}"'{pos=0.85}, shift right=5, shorten <=55pt, Rightarrow, from=1-1, to=1-3]
	\arrow[curve={height=12pt}, equals, from=1-1, to=2-3]
	\arrow["{P_C}", from=1-3, to=2-3]
\end{tikzcd}\] and, using the cloven isofibration structure on $P_C$, we obtain a section $\t_{\tilde \splitt_c}^{\eta_A^B}:\Gamma.\sigmad_A^B\to\Gamma.\sigmad_A^B.C$ of $P_C$, as well as a 2-cell:
\[\begin{tikzcd}[cramped,column sep=scriptsize]
	{\Gamma.\sigmad_A^B} &&&& {\Gamma.\sigmad_A^B}
	\arrow[""{name=0, anchor=center, inner sep=0}, "{\tilde\splitt_c}"', curve={height=18pt}, from=1-1, to=1-5]
	\arrow[""{name=1, anchor=center, inner sep=0}, "{\t_{\tilde\splitt_c}^{\eta_A^B}}", curve={height=-18pt}, from=1-1, to=1-5]
	\arrow["{\;\tau_{\tilde\splitt_c}^{\eta_A^B}}", shorten <=7pt, shorten >=7pt, Rightarrow, from=1, to=0]
\end{tikzcd}\hspace{-0.1em}{.C}\] such that $P_C\;*\;\tau_{\tilde\splitt_c}^{\eta_A^B}=\eta_A^B$. We define $\splitt_c:=\t_{\tilde\splitt_c}^{\eta_A^B}$.

\medskip

\noindent\textbf{Comp Axiom for $\boldsymbol{\Sigma}$-types.} Let us condider the 2-cell: 
\[\begin{tikzcd}[cramped,column sep=small,row sep=tiny]
	{\splitt_c\pair_A^B} &&&& {\tilde\splitt_c\pair_A^B} &&&& {(\pair_A^B)^\lilbullet c}
	\arrow["{\tau_{\tilde\splitt_c}^{\eta_A^B}\;*\;\pair_A^B}", Rightarrow, from=1-1, to=1-5]
	\arrow["{(\pair_A^B)^\lilbullet c\;*\;\beta_A^B}", Rightarrow, from=1-5, to=1-9]
\end{tikzcd}\] and observe that: $$P_{C}*(\splitt_c\pair_A^B\Longrightarrow(\pair_A^B)^\lilbullet c)= (\eta_A^B*\pair_A^B)(\pair_A^B*\beta_A^B)=1_{\pair_A^B}$$ since, being the given $(\pair_A^B,\pi_A^B,\eta_A^B,\beta_A^B)$ an equivalence, we can assume without loss of generality that it is in fact an adjoint equivalence---up to replacing $\eta_A^B$ with a parallel 2-cell \cite{MR2664622}, \cite{MR0220789}, \cite[Chapter 6]{MR0371990}. Therefore, referring to the diagram: 
\begin{equation}\label{reindexingsigma}\begin{tikzcd}[column sep=scriptsize,row sep=small]
	{\Gamma.A.B} && {\Gamma.A.B} \\
	\\
	{\Gamma.A.B.C[\pair_A^B]} && {\Gamma.\sigmad_A^B.C} \\
	\\
	{\Gamma.A.B} && {\Gamma.\sigmad_A^B}
	\arrow[equals, from=1-1, to=1-3]
	\arrow["{\splitt_c[\pair_A^B]}"', shift right=3, from=1-1, to=3-1]
	\arrow["c", shift left=3, from=1-1, to=3-1]
	\arrow[""{name=0, anchor=center, inner sep=0}, "{(\pair_A^B)^\lilbullet c}", shift left=3, from=1-3, to=3-3]
	\arrow[""{name=1, anchor=center, inner sep=0}, "{\splitt_c\pair_A^B}"', shift right=3, from=1-3, to=3-3]
	\arrow["{(\pair_A^B)^\lilbullet}"{description}, from=3-1, to=3-3]
	\arrow["{P_{C[\pair_A^B]}}"', from=3-1, to=5-1]
	\arrow["\lrcorner"{anchor=center, pos=0.125}, draw=none, from=3-1, to=5-3]
	\arrow["{P_C}", from=3-3, to=5-3]
	\arrow["{\pair_A^B}"{description}, from=5-1, to=5-3]
	\arrow[shorten <=2pt, shorten >=2pt, Rightarrow, from=1, to=0]
\end{tikzcd}\end{equation} by the 2-universal property of $\Gamma.A.B.C[\pair_A^B]$ there is unique a 2-cell:
\[\begin{tikzcd}[sep=small]
	{\Gamma.A.B} &&&& {\Gamma.A.B}
	\arrow[""{name=0, anchor=center, inner sep=0}, "c"', curve={height=12pt}, from=1-1, to=1-5]
	\arrow[""{name=1, anchor=center, inner sep=0}, "{\splitt_c[\pair_A^B]}", curve={height=-12pt}, from=1-1, to=1-5]
	\arrow["{\;s_c}", shorten <=5pt, shorten >=5pt, Rightarrow, from=1, to=0]
\end{tikzcd}\hspace{-0.4em}{.C[\pair_A^B]}\] such that $P_{C[\pair_A^B]}*s_c= 1_{1_{\Gamma.A.B}}$ and $(\splitt_c\pair_A^B\Rightarrow(\pair_A^B)^\lilbullet c)=(\pair_A^B)^\lilbullet * s_c$. As for the paragraph \textbf{Comp Axiom for $\boldsymbol{=}$-types}, we define $\sigma_c$ as the unique section: $$\{s_c\}:\Gamma.A.B\to\Gamma.A.B.\id_{C[\pair_A^B]}[\splitt_c[\pair_A^B];c]$$ of $P_{\id_{C[\pair_A^B]}[\splitt_c[\pair_A^B];c]}$ such that (\ref{whatcharacterisescurlyparentheses}) is satisfied---where $A$ is $C[\pair_A^B]$, $a$ and $b$ are $\splitt_c[\pair_A^B]$ and $c$ respectively, and $p$ is $s_c$.

\begin{proposition}\label{stability of sigma in 2d model}
    Let $(\boldsymbol{\mathcal{C}},\mathcal{D})$ be a display map 2-category endowed with axiomatic $=$-types and axiomatic $\Sigma$-types. Referring to the data defined in paragraphs \emph{\textbf{Form Rule, Intro Rule, Elim Rule, and Comp Axiom, for $\boldsymbol{\Sigma}$-types}}, the \emph{stability conditions}: $$\begin{alignedat}{2}
        \sigmad_A^B[f]&=\sigmad_{A[f]}^{B[f^\lilbullet]}&\quad\quad\splitt_c[f^\lilbullet]&=\splitt_{c[f^{\lilbullet\lilbullet}]}\\
        \pairing_A^B[f^{\lilbullet\lilbullet}]&=\pairing_{A[f]}^{B[f^\lilbullet]}&\sigma_c[f^{\lilbullet\lilbullet}]&=\sigma_{c[f^{\lilbullet\lilbullet}]}\\
    \end{alignedat}$$ of Defintion \ref{def:Semantics of axiomatic Sigma-types} hold for every arrow $\Delta \xrightarrow{f}\Gamma$.

\proof
See Appendix \ref{appendixIII}.
\endproof
\end{proposition}

\subsection{Axiomatic $\boldsymbol{\Pi}$-types}

We specialise the notion of a display map 2-category with axiomatic $=$-types to obtain a version that incorporates a 2-dimensional semantic interpretation of axiomatic $\Pi$-types \& axiomatic function extensionality. Again, the condition we require to be satisfied is a generalisation of the condition characterising the semantics of extensional $\Pi$-types: right adjoints \cite{MR1674451} are demoted to right biadjoints. In particular, we observe the same difference between the intensional case \cite{MR2525957} and the axiomatic one, that we observed for $\Sigma$-types: such a right biadjoint need not be a retract right biadjoint.

\begin{definition}[Semantics of axiomatic $\Pi$-types \& axiomatic function extensionality---categorical formulation]\label{def:Categorical semantics of axiomatic Pi-types}
    Let $(\boldsymbol{\mathcal{C}},\mathcal{D})$ be a display map 2-category endowed with axiomatic $=$-types.

    Let us assume that, for every object $\Gamma$, every display map $P_A$ of codomain $\Gamma$, and every display map $P_B$ of codomain $\Gamma.A$, there is a choice of a display map $\Gamma.\pid_A^B\to\Gamma$ in such a way that: $$\Pi_A^B[f]=\Pi_{A[f]}^{B[f^\lilbullet]}$$ for every arrow $f:\Delta\to\Gamma$, and in such a way that the assignment: $$\begin{aligned}
        \mathcal{D}/\Gamma.A &\to \boldsymbol{\mathcal{C}}/\Gamma\\
        P_B&\mapsto P_{\pid_A^B}
    \end{aligned}$$ constitutes a right $I_A$-relative 2-coadjoint to the 2-functor: $$\begin{aligned}
        \boldsymbol{\mathcal{C}}/\Gamma &\to \boldsymbol{\mathcal{C}}/\Gamma.A\\
        f&\mapsto f^\lilbullet
    \end{aligned}$$ where $I_A$ is the full forgetful 2-functor $\mathcal{D}/\Gamma.A\hookrightarrow\boldsymbol{\mathcal{C}}/\Gamma.A$.
    
    Unfolding this requirement, this means that there is a family of equivalences of categories: $$(\lambda_A^B,\app_A^B,\app_A^B\lambda_A^B\overset{\beta^{A,B}}{\Longrightarrow}1,1\overset{\eta^{A,B}}{\Longrightarrow}\lambda_A^B\app_A^B):(\boldsymbol{\mathcal{C}}/\Gamma.A)(f^\lilbullet,P_B)\simeq (\boldsymbol{\mathcal{C}}/\Gamma)(f,P_{\pid_A^B})$$ 2-natural in $f:\Delta\to\Gamma$ and $P_B$---which means that the diagrams of 1-functors: 
\[\begin{tikzcd}[column sep=tiny]
	{(\boldsymbol{\mathcal{C}}/\Gamma.A)(f^\lilbullet,P_B)} && {(\boldsymbol{\mathcal{C}}/\Gamma)(f,P_{\pid_A^B})} & {(\boldsymbol{\mathcal{C}}/\Gamma.A)(f^\lilbullet,P_B)} && {(\boldsymbol{\mathcal{C}}/\Gamma)(f,P_{\pid_A^B})} \\
	\\
	{(\boldsymbol{\mathcal{C}}/\Gamma.A)((fg)^\lilbullet,P_B)} && {(\boldsymbol{\mathcal{C}}/\Gamma)(fg,P_{\pid_A^B})} & {(\boldsymbol{\mathcal{C}}/\Delta.A[f])(1_{\Delta}^\lilbullet,P_{B[f^\lilbullet]})} && {(\boldsymbol{\mathcal{C}}/\Delta)(1_\Delta,P_{\pid_{A[f]}^{B[f^\lilbullet]}})}
	\arrow["{\lambda_A^B}", shift left, from=1-1, to=1-3]
	\arrow["{(\blank)(g^\lilbullet)}"{description}, from=1-1, to=3-1]
	\arrow["{\app_A^B}", shift left, from=1-3, to=1-1]
	\arrow["{(\blank)(g)}"{description}, from=1-3, to=3-3]
	\arrow["{\lambda_A^B}", shift left, from=1-4, to=1-6]
	\arrow["{1_{\Delta.A[f]};\blank}"{description}, from=1-4, to=3-4]
	\arrow["{\app_A^B}", shift left, from=1-6, to=1-4]
	\arrow["{1_{\Delta};\blank}"{description}, from=1-6, to=3-6]
	\arrow["{\lambda_A^B}", shift left, from=3-1, to=3-3]
	\arrow["{\app_A^B}", shift left, from=3-3, to=3-1]
	\arrow["{\lambda_{A[f]}^{B[f.A]}}", shift left, from=3-4, to=3-6]
	\arrow["{\app_{A[f]}^{B[f.A]}}", shift left, from=3-6, to=3-4]
\end{tikzcd}\] commute and the equalities: $$\begin{alignedat}{2}
        \beta^{A,B}_h g^\lilbullet&=\beta^{A,B}_{hg^\lilbullet}&\quad\quad\quad\quad\quad\quad\quad\quad1_{\Delta.A[f]};\beta^{A,B}_h&=\beta^{A[f],B[f.A]}_{1_\Delta.A[f];h}\\
        \eta^{A,B}_kg&=\eta^{A,B}_{kg}&1_\Delta;\eta^{A,B}_k&=\eta^{A[f],B[f.A]}_{1_\Delta;k}\\
    \end{alignedat}$$ hold, for every $g:\Delta'\to\Delta$, every $h:f^\lilbullet\to P_B$, and every $k:f\to P_{\pid_A^B}$.

    Then we say that $(\boldsymbol{\mathcal{C}},\mathcal{D})$ is \textit{endowed with axiomatic $\Pi$-types \& axiomatic function extensionality}.
\end{definition}

Let $(\boldsymbol{\mathcal{C}},\mathcal{D})$ be a display map 2-category endowed with axiomatic $=$-types.

Let us assume that the display map 2-category $(\boldsymbol{\mathcal{C}},\mathcal{D})$ is endowed with axiomatic $\Pi$-types \& axiomatic function extensionality---Definition \ref{def:Categorical semantics of axiomatic Pi-types}---let $P_A$ be a display map of codomain $\Gamma$ and let $P_B$ be a display map of codomain $\Gamma.A$. The next paragraphs show that the display map category $(\mathbf{C},\mathcal{D})$ induced by $(\boldsymbol{\mathcal{C}},\mathcal{D})$ is endowed with appropriate choice functions validating axiomatic $\Pi$-types \& axiomatic function extensionality---as in Definition \ref{def:Semantics of axiomatic Pi-types}.

\medskip

\noindent\textbf{Form Rule for $\boldsymbol{\Pi}$-types.} We already have a choice of a display map $\Gamma.\Pi_A^B\to\Gamma$.

\medskip

\noindent\textbf{Elim Rule for $\boldsymbol{\Pi}$-types.} The diagram:
\[\begin{tikzcd}
	{(\boldsymbol{\mathcal{C}}/\Gamma.A)(P_{\pid_A^B}^\lilbullet,P_B)} && {(\boldsymbol{\mathcal{C}}/\Gamma)(P_{\pid_A^B},P_{\pid_A^B})} \\
	\\
	{(\boldsymbol{\mathcal{C}}/\Gamma.\pid_A^B.A[P_{\pid_A^B}])(1_{\Gamma.\pid_A^B}^\lilbullet,P_{B[P_{\pid_A^B}^\lilbullet]})} && {(\boldsymbol{\mathcal{C}}/\Gamma.\pid_A^B)(1_{\Gamma.\pid_A^B},P_{{\pid_{A}^{B}}^\liltriangle)}}
	\arrow["{1_{\Gamma.\pid_A^B.A[P_{\pid_A^B}]};\blank}"{description}, from=1-1, to=3-1]
	\arrow["{\app_A^B}"', from=1-3, to=1-1]
	\arrow["{1_{\Gamma.\pid_A^B};\blank}"{description}, from=1-3, to=3-3]
	\arrow["{\app_{A[P_{\pid_A^B}]}^{B[P_{\pid_A^B}^\lilbullet]}}", from=3-3, to=3-1]
\end{tikzcd}\] commutes---observe that ${\pid_A^B}^\liltriangle=\pid_{A}^{B}[P_{\pid_A^B}]=\pid_{A[P_{\pid_A^B}]}^{B[P_{\pid_A^B}^\lilbullet]}$. Hence we may define $\ev_A^B$ as the section: $$\app_{A[P_{\pid_A^B}]}^{B[P_{\pid_A^B}^\lilbullet]}\delta_{\pid_A^B}=\app_{A[P_{\pid_A^B}]}^{B[P_{\pid_A^B}^\lilbullet]}(1_{\Gamma.\pid_A^B};1_{\Gamma.\pid_A^B})=1_{\Gamma.\pid_A^B.A[P_{\pid_A^B}]};\app_A^B1_{\Gamma.\pid_A^B}$$ of $P_{B[P_{\pid_A^B}^\lilbullet]}$.

\medskip

\noindent\textbf{Intro Rule for $\boldsymbol{\Pi}$-types.} If $\upsilon$ is a section of $P_B$ then we define the section $\abst{\upsilon}$ of $P_{\pid_A^B}$ as the image of $\upsilon$ via the map $\lambda_A^B:(\boldsymbol{\mathcal{C}}/\Gamma.A)(1_{\Gamma.A},P_B)\to(\boldsymbol{\mathcal{C}}/\Gamma)(1_\Gamma,P_{\pid_A^B})$.

\medskip

\noindent\textbf{Comp Axiom for $\boldsymbol{\Pi}$-types.} We observe that: $$\ev_A^B[\abst{\upsilon}^\lilbullet]=P_{\pid_A^B}^{\lilbullet\lilbullet}\ev_A^B\abst{\upsilon}^\lilbullet$$ because the diagram:
\[\begin{tikzcd}[column sep=small]
	{\Gamma.A} & {\Gamma.\pid_A^B.A[P_{\pid_A^B}]} \\
	{\Gamma.A.B} & {\Gamma.\pid_A^B.A[P_{\pid_A^B}].B[P_{\pid_A^B}^\lilbullet]} \\
	&& {\Gamma.A.B}
	\arrow["{z^\lilbullet}", from=1-1, to=1-2]
	\arrow["{\ev_A^B[z^\lilbullet]}"', from=1-1, to=2-1]
	\arrow["{\ev_A^B}", from=1-2, to=2-2]
	\arrow["{z^{\lilbullet\lilbullet}}", from=2-1, to=2-2]
	\arrow[curve={height=12pt}, equals, from=2-1, to=3-3]
	\arrow["{P_{\pid_A^B}^{\lilbullet\lilbullet}}"{description}, from=2-2, to=3-3]
\end{tikzcd}\]
commutes for every section $z$ of $P_{\pid_A^B}$. Hence: $$\begin{aligned}
    \ev_A^B[z^\lilbullet]&=P_{\pid_A^B}^{\lilbullet\lilbullet}\ev_A^Bz^\lilbullet\\
    &=P_{\pid_A^B}^{\lilbullet\lilbullet}(1_{\Gamma.\pid_A^B.A[P_{\pid_A^B}]};\app_A^B1_{\Gamma.\pid_A^B})z^\lilbullet\\
    &=(\app_A^B1_{\Gamma.\pid_A^B})z^\lilbullet\\
    &=\app_A^B(1_{\Gamma.\pid_A^B}z)\\
    &=\app_A^Bz
\end{aligned}$$ that is: \begin{equation}\label{evèapp}
    \ev_A^B[z^\lilbullet]=\app_A^Bz
\end{equation} for every section $z$ of $P_{\pid_A^B}$. In particular: $$\ev_A^B[\abst{\upsilon}^\lilbullet]=\app_A^B\abst{\upsilon}=\app_A^B\lambda_A^B\upsilon$$ and therefore $\beta^{A,B}_\upsilon$ is a morphism: $$\ev_A^B[\abst{\upsilon}^\lilbullet]\Rightarrow\upsilon$$ in the category $(\boldsymbol{\mathcal{C}}/\Gamma.A)(1_{\Gamma.A},P_B)$ of the sections of $P_B$. As for the paragraphs \textbf{Comp Axiom for $\boldsymbol{=}$-types} and \textbf{Comp Axiom for $\boldsymbol{\Sigma}$-types}, we define $\beta_\upsilon$ as the unique section: $$\{\beta^{A,B}_\upsilon\}:\Gamma.A\to\Gamma.A.\id_{B}[\ev_A^B[\abst{\upsilon}^\lilbullet];\upsilon]$$ of $P_{\id_{B}[\ev_A^B[\abst{\upsilon}^\lilbullet];\upsilon]}$ such that (\ref{whatcharacterisescurlyparentheses}) is satisfied---where $A$ is $B$, $a$ and $b$ are $\ev_A^B[\abst{\upsilon}^\lilbullet]$ and $\upsilon$ respectively, and $p$ is $\beta^{A,B}_\upsilon$.

\medskip

\noindent\textbf{Stability of $\boldsymbol{\app}$ and $\boldsymbol{\lambda}$ under re-indexing.} By naturalities, if $z$ is a section of $P_{\pid_A^B}$, $\upsilon$ is a section of $P_B$, and $f$ is an arrow $\Delta\to\Gamma$, we observe that: \[\begin{alignedat}{2}
    (\app_A^B z)[f^\lilbullet]&=1_{\Delta.A[f]};(\app_A^Bz)f^\lilbullet& \quad\quad\quad\quad (\lambda_A^B \upsilon)[f]&=1_\Delta;(\lambda_A^B\upsilon)f\\
    &=1_{\Delta.A[f]};\app_A^B(zf)&&=1_\Delta;\lambda_A^B(\upsilon f^\lilbullet)\\
    &=\app_{A[f]}^{B[f^\lilbullet]}(1_{\Delta};zf)&&=\lambda_{A[f]}^{B[f^\lilbullet]}(1_{\Delta.A[f]};\upsilon f^\lilbullet)\\
    &=\app_{A[f]}^{B[f^\lilbullet]}(z[f])&&=\lambda_{A[f]}^{B[f^\lilbullet]}(\upsilon [f^\lilbullet])
\end{alignedat}\] that is: \begin{equation}\label{appandlambdastableunderpullback}
    (\app_A^B z)[f^\lilbullet]=\app_{A[f]}^{B[f^\lilbullet]}(z[f]) \quad\quad\quad\quad\quad\quad (\lambda_A^B \upsilon)[f]=\lambda_{A[f]}^{B[f^\lilbullet]}(\upsilon [f^\lilbullet])\quad.
\end{equation} In particular: $$\begin{aligned}\ev_A^B[f^{\lilbullet\lilbullet}]&=(\app_{A[P_{\pid_A^B}]}^{B[P_{\pid_A^B}^\lilbullet]}\delta_{\pid_A^B})[f^{\lilbullet\lilbullet}]\overset{(\ref{appandlambdastableunderpullback})}{=}\app_{A[P_{\pid_A^B}][f^{\lilbullet}]}^{B[P_{\pid_A^B}^\lilbullet][f^{\lilbullet\lilbullet}]}\delta_{\pid_A^B}[f^{\lilbullet}]
\\&=\app_{A[f][P_{\pid_{A[f]}^{B[f^\lilbullet]}}]}^{B[f^\lilbullet][P_{\pid_{A[f]}^{B[f^\lilbullet]}}^\lilbullet]}\delta_{\pid_A^B[f]}=\app_{A[f][P_{\pid_{A[f]}^{B[f^\lilbullet]}}]}^{B[f^\lilbullet][P_{\pid_{A[f]}^{B[f^\lilbullet]}}^\lilbullet]}\delta_{\pid_{A[f]}^{B[f^\lilbullet]}}\\
&=\ev_{A[f]}^{B[f^\lilbullet]}
\end{aligned}$$ that is: \begin{equation}\label{evstabilesottoreindexing}
    \ev_A^B[f^{\lilbullet\lilbullet}] = \ev_{A[f]}^{B[f^\lilbullet]}\quad.
\end{equation}

\medskip

\noindent\textbf{Characterising \boldsymbol{\happly}.} Referring to the notation of \Cref{howtohapply}, let: 
\[\begin{tikzcd}
	{\Gamma.\pid_A^B.(\pid_A^B)^\liltriangle.\id_{\pid_A^B}} &&& {\Gamma.\pid_A^B.(\pid_A^B)^\liltriangle.\id_{\pid_A^B}.\id_{A[P_{\pid_A^B}\mathsf{s}]}^{B[(P_{\pid_A^B}\mathsf{s})^\lilbullet]}}
	\arrow[""{name=0, anchor=center, inner sep=0}, "{\tilde{\mathsf{s}}}", shift left=3, from=1-1, to=1-4]
	\arrow[""{name=1, anchor=center, inner sep=0}, "{\tilde{\mathsf{t}}}"', shift right=3, from=1-1, to=1-4]
	\arrow["{\tilde\alpha_{\pid_A^B}}", shorten <=2pt, shorten >=2pt, Rightarrow, from=0, to=1]
\end{tikzcd}\hspace{-0.2em}\] be the re-indexing of $\alpha_{\pid_A^B}$ via $P_{\pid_A^B}\mathsf{s}:\Gamma.\pid_A^B.(\pid_A^B)^\liltriangle.\id_{\pid_A^B}\to\Gamma$ and let: $$\happ:\Gamma.\pid_A^B.(\pid_A^B)^\liltriangle.\id_{\pid_A^B}.A[P_{\pid_A^B}\mathsf{s}]\to\Gamma.\pid_A^B.(\pid_A^B)^\liltriangle.\id_{\pid_A^B}.A[P_{\pid_A^B}\mathsf{s}].\id_{B[(P_{\pid_A^B}\mathsf{s})^\lilbullet]}[\app\,\tilde{\mathsf{s}};\app\,\tilde{\mathsf{t}}]$$ be $\{\app\,\tilde\alpha_{\pid_A^B}\}$ i.e. the re-indexing along $\app\,\tilde{\mathsf{s}};\app\,\tilde{\mathsf{t}}$ of the unique factorisation of $\app\,\tilde\alpha_{\pid_A^B}$ through the arrow object $\alpha_{B[(P_{\pid_A^B})^\lilbullet]}$---see (\ref{whatcharacterisescurlyparentheses}). Now, we observe that: $$\app\, \tilde{\mathsf{s}}=\app_{A[P_{\pid_A^B}\mathsf{s}]}^{B[(P_{\pid_A^B}\mathsf{s})^\lilbullet]}\tilde{\mathsf{s}}\overset{(\ref{evèapp})}{=}\ev_{A[P_{\pid_A^B}\mathsf{s}]}^{B[(P_{\pid_A^B}\mathsf{s})^\lilbullet]}[\tilde{\mathsf{s}}^\lilbullet]=\ev_{A}^{B}[(P_{\pid_A^B}\mathsf{s})^{\lilbullet\lilbullet}][\tilde{\mathsf{s}}^\lilbullet]=\ev_A^B[((P_{\pid_A^B}\mathsf{s})^{\lilbullet}\tilde{\mathsf{s}})^\lilbullet]=\ev_A^B[\mathsf{s}^\lilbullet]$$ and analogously $\app\, \tilde{\mathsf{t}}=\ev_A^B[\mathsf{t}^\lilbullet]$. Therefore $\happ$ is an arrow: $$\Gamma.\pid_A^B.(\pid_A^B)^\liltriangle.\id_{\pid_A^B}.A[P_{\pid_A^B}\mathsf{s}]\to\Gamma.\pid_A^B.(\pid_A^B)^\liltriangle.\id_{\pid_A^B}.A[P_{\pid_A^B}\mathsf{s}].\id_{B[(P_{\pid_A^B}\mathsf{s})^\lilbullet]}[\ev_A^B[\mathsf{s}^\lilbullet];\ev_A^B[\mathsf{t}^\lilbullet]]$$ and hence $\lambda\,\happ$ is parallel to $\happly$. We can now state and prove the following:

\begin{theorem}[The happy theorem]\label{happytheorem} Let $(\boldsymbol{\mathcal{C}},\mathcal{D})$ be a display map 2-category endowed with axiomatic $=$-types and axiomatic $\Pi$-types \& axiomatic function extensionality. Then there exists an arrow: $$\happly\Rightarrow\lambda\,\happ$$in the category of the sections of the display map: $$\Gamma.\pid_A^B.(\pid_A^B)^\liltriangle.\id_{\pid_A^B}.\bigg(\pid_{A[P_{\pid_A^B}\mathsf{s}]}^{\id_{B[(P_{\pid_A^B}\mathsf{s})^\lilbullet]}[\ev_A^B[\mathsf{s}^\lilbullet];\ev_A^B[\mathsf{t}^\lilbullet]]}\bigg)\to\Gamma.\pid_A^B.(\pid_A^B)^\liltriangle.\id_{\pid_A^B}$$ which is stable under re-indexing.

\proof
Using the appropriate arrow object, it is enough to build a section of the display map: $$\Gamma.\pid_A^B.(\pid_A^B)^\liltriangle.\id_{\pid_A^B}.\id_{\bigg(\pid_{A[P_{\pid_A^B}\mathsf{s}]}^{\id_{B[(P_{\pid_A^B}\mathsf{s})^\lilbullet]}[\ev_A^B[\mathsf{s}^\lilbullet];\ev_A^B[\mathsf{t}^\lilbullet]]}\bigg)}[\happly;\lambda\,\happ]\to\Gamma.\pid_A^B.(\pid_A^B)^\liltriangle.\id_{\pid_A^B}$$ and, by elimination, such a section can be defined as $\j_c$ if we build a section $c$ of the display map: $$\Gamma.\pid_A^B.\id_{\bigg(\pid_{A[P_{\pid_A^B}\mathsf{s}]}^{\id_{B[(P_{\pid_A^B}\mathsf{s})^\lilbullet]}[\ev_A^B[\mathsf{s}^\lilbullet];\ev_A^B[\mathsf{t}^\lilbullet]]}\bigg)}[\happly;\lambda\,\happ][\refl_{\pid_A^B}]\to\Gamma.\pid_A^B$$ that is: $$\Gamma.\pid_A^B.\id_{\bigg(\pid_{A[P_{\pid_A^B}]}^{\id_{B[P_{\pid_A^B}^\lilbullet]}[\ev_A^B;\ev_A^B]}\bigg)}[\;\happly[\refl_{\pid_A^B}]\;;\;(\lambda\,\happ)[\refl_{\pid_A^B}]\;]\to\Gamma.\pid_A^B\quad.$$ Again, using the appropriate arrow object, it is enough to build an arrow: $$\begin{aligned}
    \happly[\refl_{\pid_A^B}]\Rightarrow(\lambda\,\happ)[\refl_{\pid_A^B}]&\overset{(\ref{appandlambdastableunderpullback})}{=}\lambda(\happ[\refl_{\pid_A^B}^\lilbullet])=\lambda(\{\app\,\tilde\alpha_{\pid_A^B}\}[\refl_{\pid_A^B}^\lilbullet])\\
    &\overset{(\ref{curlyparenthesespreservereindexing})}{=}\lambda\{(\app\,\tilde\alpha_{\pid_A^B})[\refl_{\pid_A^B}^\lilbullet]\}\\
    &\overset{(\ref{appandlambdastableunderpullback})}{=}\lambda\{\app(\tilde\alpha_{\pid_A^B}[\refl_{\pid_A^B}])\}=\lambda\{\app\,1_{\delta_{\pid_A^B}}\}=\lambda\{1_{\app\,\delta_{\pid_A^B}}\}=\lambda\{1_{\ev_A^B}\}
\end{aligned}$$ in the category of the sections of $\Gamma.\pid_A^B.\pid_{A[P_{\pid_A^B}]}^{\id_{B[P_{\pid_A^B}^\lilbullet]}[\ev_A^B;\ev_A^B]}\to\Gamma.\pid_A^B$ and, since in this category there is an arrow: $$\happly[\refl_{\pid_A^B}]=\j_{\abst{\refl_{B[P_{\pid_A^B}^\lilbullet]}[\ev_A^B]}}[\refl_{\pid_A^B}]\Rightarrow\abst{\refl_{B[P_{\pid_A^B}^\lilbullet]}[\ev_A^B]}=\lambda(\refl_{B[P_{\pid_A^B}^\lilbullet]}[\ev_A^B])$$ it is enough to build an arrow: $$\lambda(\refl_{B[P_{\pid_A^B}^\lilbullet]}[\ev_A^B])\Rightarrow \lambda\{1_{\ev_A^B}\}$$ in the same category. Therefore we are left to build an arrow $\refl_{B[P_{\pid_A^B}^\lilbullet]}[\ev_A^B]\Rightarrow \{1_{\ev_A^B}\}$ which amounts to building an arrow: $$\alpha_{B[P_{\pid_A^B}^\lilbullet]}(\ev_A^B;\ev_A^B)^\lilbullet\refl_{B[P_{\pid_A^B}^\lilbullet]}[\ev_A^B]\Rightarrow1_{\ev_A^B}$$ in the category $(\boldsymbol{\mathcal{C}}/\Gamma.A)(1_{\Gamma.\pid_A^B.A[P_{\pid_A^B}]},P_{B[P_{\pid_A^B}^\lilbullet]})^\to$. However: $$\alpha_{B[P_{\pid_A^B}^\lilbullet]}(\ev_A^B;\ev_A^B)^\lilbullet\refl_{B[P_{\pid_A^B}^\lilbullet]}[\ev_A^B]=\alpha_{B[P_{\pid_A^B}^\lilbullet]}\refl_{B[P_{\pid_A^B}^\lilbullet]}\ev_A^B=1_{1_{\Gamma.\pid_A^B.A[P_{\pid_A^B}].B[P_{\pid_A^B}^\lilbullet]}}\ev_A^B=1_{\ev_A^B}$$ hence the identity is such an arrow. In conclusion, an arrow: $$\happly[\refl_{\pid_A^B}]\Rightarrow(\lambda\,\happ)[\refl_{\pid_A^B}]$$ is just the one associated to the computation axiom: $${\h_{\abst{{\refl_{B[P_{\pid_A^B}^\lilbullet]}[\ev_A^B]}}}}$$ for $\happly[\refl_{\pid_A^B}]$---which is stable under re-indexing by \Cref{stability of id in 2d model}. Then: $$\j_{\h_{\abst{{\refl_{B[P_{\pid_A^B}^\lilbullet]}[\ev_A^B]}}}}$$ is a section of: $$\Gamma.\pid_A^B.(\pid_A^B)^\liltriangle.\id_{\pid_A^B}.\id[\happly;\lambda\,\happ]\to\Gamma.\pid_A^B.(\pid_A^B)^\liltriangle.\id_{\pid_A^B}$$hence, as we said, we can define: $$\happly\Rightarrow\lambda\,\happ$$ as: $$\alpha_{\bigg(\pid_{A[P_{\pid_A^B}\mathsf{s}]}^{\id_{B[(P_{\pid_A^B}\mathsf{s})^\lilbullet]}[\ev_A^B[\mathsf{s}^\lilbullet];\ev_A^B[\mathsf{t}^\lilbullet]]}\bigg)}(\happly;\lambda\,\happ)^\lilbullet\j_{\h_{\abst{{\refl_{B[P_{\pid_A^B}^\lilbullet]}[\ev_A^B]}}}}\quad.$$ This is enough to verify that $\happly\Rightarrow\lambda\,\happ$ is stable under every re-indexing for the domain of $P_A$. Happy ending.
\endproof
\end{theorem}

\medskip

\noindent\textbf{Intro Rule for function extensionality.} Let $z,z'$ be sections of the display map $\Gamma.\pid_A^B\to\Gamma$ and let $q$ be a section of the display map: $$\Gamma.\pid_A^{\id_B[\ev_A^B[z^\lilbullet];\ev[z'^\lilbullet]]}\overset{(\ref{evèapp})}{=}\Gamma.\pid_A^{\id_B[\app_A^Bz;\app_A^Bz']}\to\Gamma\;\;.$$ Then $\app_A^{\id_B[\ev_A^B[z^\lilbullet];\ev_A^B[z'^\lilbullet]]} q$ is a section of the display map $\Gamma.A.\id_B[\app_A^Bz;\app_A^Bz']\to\Gamma.A$ and therefore: $$\alpha_B(\app_A^Bz;\app_A^Bz')^\lilbullet \app_A^{\id_B[\ev_A^B[z^\lilbullet];\ev_A^B[z'^\lilbullet]]}q$$ is an arrow $\app_A^Bz\Rightarrow\app_A^Bz'$ in the category $(\boldsymbol{\mathcal{C}}/\Gamma.A)(1_{\Gamma.A},P_B)$ of the sections of $P_B$. Therefore: $$(\eta^{A,B}_{z'})^{-1}\lambda_A^B(\alpha_B(\app_A^Bz;\app_A^Bz')^\lilbullet \app_A^{\id_B[\ev_A^B[z^\lilbullet];\ev_A^B[z'^\lilbullet]]}q)(\eta^{A,B}_z)$$ is an arrow $z\Rightarrow z'$ in the category of the sections of $P_{\pid_A^B}$. We define $\funext_q$ as: $$\{(\eta^{A,B}_{z'})^{-1}\lambda_A^B(\alpha_B(\app_A^Bz;\app_A^Bz')^\lilbullet \app_A^{\id_B[\ev_A^B[z^\lilbullet];\ev_A^B[z'^\lilbullet]]}q)(\eta^{A,B}_z)\}$$ i.e. the re-indexing via $z;z'$ of the unique factorisation of: $$(\eta^{A,B}_{z'})^{-1}\lambda_A^B(\alpha_B(\app_A^Bz;\app_A^Bz')^\lilbullet \app_A^{\id_B[\ev_A^B[z^\lilbullet];\ev_A^B[z'^\lilbullet]]}q)(\eta^{A,B}_z)$$ through $\alpha_{\pid_A^B}$.

\medskip

\noindent\textbf{Comp Axiom for function extensionality.} If we build an arrow $\happly_{z;z'}[\funext_q]\Rightarrow q$ in the category of the sections of $\Gamma.\pid_A^{\id_B[\app_A^Bz;\app_A^Bz']}\to\Gamma$ then $\beta^\Pi_q$ can be defined as the re-indexing via $\happly_{z;z'}[\funext_q];q$ of the unique factorisation of $\happly_{z;z'}[\funext_q]\Rightarrow q$ through $\alpha_{\pid_A^{\id_B[\app_A^Bz;\app_A^Bz']}}$. By \Cref{happytheorem}, we know that: $$\happly_{z;z'}[\funext_q]=\happly[(z;z')^\lilbullet\funext_q]\Rightarrow (\lambda\,\happ)[(z;z')^\lilbullet\funext_q]$$ hence we are left to build $(\lambda\,\happ)[(z;z')^\lilbullet\funext_q]\Rightarrow q$.

Let $p$ be a section of the display map $\Gamma.\id_{\pid_A^B}[z;z']$. Then: $$(\lambda\,\happ)[(z;z')^\lilbullet p]\overset{(\ref{appandlambdastableunderpullback})}{=}\lambda(\happ[((z;z')^\lilbullet p)^\lilbullet])=\lambda(\{\app \,\tilde\alpha_{\pid_A^B}\}[((z;z')^\lilbullet p)^\lilbullet])$$ where we recall that $\{\app \,\tilde\alpha_{\pid_A^B}\}$ denotes the re-indexing along $\app\,\tilde{\mathsf{s}};\app\,\tilde{\mathsf{t}}$ of the unique factorisation of $\app\,\tilde\alpha_{\pid_A^B}$ through the arrow object $\alpha_{B[(P_{\pid_A^B})^\lilbullet]}$. Therefore: $$(\lambda\,\happ)[(z;z')^\lilbullet p]\overset{(\ref{curlyparenthesespreservereindexing})}{=}\lambda\{(\app \,\tilde\alpha_{\pid_A^B})[((z;z')^\lilbullet p)^\lilbullet]\}\overset{(\ref{appandlambdastableunderpullback})}{=}\lambda\{\app(\tilde\alpha_{\pid_A^B}[(z;z')^\lilbullet p])\}$$ and, since: $$\tilde\alpha_{\pid_A^B}[(z;z')^\lilbullet p]=\alpha_{\pid_A^B}(z;z')^\lilbullet p$$ by definition of $\tilde\alpha_{\pid_A^B}$, we deduce that: \begin{equation}\label{funnyhapp}
    (\lambda\,\happ)[(z;z')^\lilbullet p]=\lambda\{\app(\alpha_{\pid_A^B}(z;z')^\lilbullet p)\}
\end{equation} and in particular, when $p$ is $\funext_q$, we obtain that: $$(\lambda\,\happ)[(z;z')^\lilbullet\funext_q]=\lambda\{\app(\;(\eta^{A,B}_{z'})^{-1}\lambda(\;\alpha_B(\app\, z;\app\, z')^\lilbullet \app\,q\;)(\eta^{A,B}_z)\;)\}\quad.$$ Therefore, in order to build $(\lambda\,\happ)[(z;z')^\lilbullet\funext_q]\Rightarrow q$ it is enough to build: $$\{\app(\;(\eta^{A,B}_{z'})^{-1}\lambda(\;\alpha_B(\app\, z;\app\, z')^\lilbullet \app\,q\;)(\eta^{A,B}_z)\;)\}\Rightarrow \app\,q\quad$$ which amounts to building: $$\app(\;(\eta^{A,B}_{z'})^{-1}\lambda(\;\alpha_B(\app\, z;\app\, z')^\lilbullet \app\,q\;)(\eta^{A,B}_z)\;)\Rightarrow\alpha_B(\app\, z;\app\, z')^\lilbullet \app\,q\quad$$ in the category $(\boldsymbol{\mathcal{C}}/\Gamma.A)(1_{\Gamma.A},P_B)^\to$ i.e. a commutative diagram of the form: 
\[\begin{tikzcd}[column sep=scriptsize]
	{\app\,z} && {\app\,z} \\
	\\
	{\app\,z'} && {\app\,z'}
	\arrow[Rightarrow, from=1-1, to=1-3]
	\arrow["{\app(\;(\eta^{A,B}_{z'})^{-1}\lambda(\;\alpha_B(\app\, z;\app\, z')^\lilbullet \app\,q\;)(\eta^{A,B}_z)\;)}"', Rightarrow, from=1-1, to=3-1]
	\arrow["{\alpha_B(\app\, z;\app\, z')^\lilbullet \app\,q}", Rightarrow, from=1-3, to=3-3]
	\arrow[Rightarrow, from=3-1, to=3-3]
\end{tikzcd}\] and the pair $(1_{\app\,z},1_{\app\,z'})$ constitutes such a commutative diagram. We are done.

\medskip

\noindent\textbf{Exp Axiom for function extensionality.} If we build an arrow $p \Rightarrow \funext_{\happly_{z;z'}[p]}$ in the category of the section of $\Gamma.\id_{\pid_A^B}\to \Gamma$ then $\gamma^\Pi_p$ can be defined as the section: $$\{p \Rightarrow \funext_{\happly_{z;z'}[p]}\}$$ of $\Gamma.\id_{\id_{\pid_A^B}}[p;\funext_{\happly_{z;z'}[p]}]\to\Gamma$. By definition $\funext_{\happly_{z;z'}[p]}$ is the section: $$\{(\eta^{A,B}_{z'})^{-1}\lambda(\;\alpha_B(\app\,z;\app\,z')^\lilbullet (\app\,\happly_{z;z'}[p])\;)(\eta^{A,B}_z)\}$$ and, by \Cref{happytheorem}, we know that: $$\lambda\{\app(\alpha_{\pid_A^B}(z;z')^\lilbullet p)\}\overset{(\ref{funnyhapp})}{=}(\lambda\,\happ)[(z;z')^\lilbullet p]\Rightarrow\happly[(z;z')^\lilbullet p]=\happly_{z;z'}[p]$$ hence: $$\{\app(\alpha_{\pid_A^B}(z;z')^\lilbullet p)\}\Rightarrow\app\,\lambda\{\app(\alpha_{\pid_A^B}(z;z')^\lilbullet p)\}\Rightarrow\app\,\happly_{z;z'}[p].$$ Therefore we are left to exhibit an arrow: $$\begin{aligned}
    p \Rightarrow \{(\eta^{A,B}_{z'})^{-1}&\lambda(\;\alpha_B(\app\,z;\app\,z')^\lilbullet \{\app(\alpha_{\pid_A^B}(z;z')^\lilbullet p)\}\;)(\eta^{A,B}_z)\}\\&=\{(\eta^{A,B}_{z'})^{-1}\lambda\app(\alpha_{\pid_A^B}(z;z')^\lilbullet p)(\eta^{A,B}_z)\}\\
    &=\{\alpha_{\pid_A^B}(z;z')^\lilbullet p\}\\
    &=\,p
\end{aligned}$$ and $1_p$ is such an arrow.

\begin{proposition}\label{stability of pi in 2d model}
    Let $(\boldsymbol{\mathcal{C}},\mathcal{D})$ be a display map 2-category endowed with axiomatic $=$-types and axiomatic $\Pi$-types \& axiomatic function extensionality. Referring to the data defined in paragraphs \emph{\textbf{Form Rule, Intro Rule, Elim Rule, and Comp Axiom, for $\boldsymbol{\Pi}$-types}} and \emph{\textbf{Intro Rule, Comp Axiom, and Exp Axiom, for function extensionality}}, the \emph{stability conditions}: $$\begin{alignedat}{2}
        \pid_A^B[f]&=\pid_{A[f]}^{B[f^\lilbullet]}&\quad\quad\abst{\upsilon}[f]&=\abst{\upsilon[f^{\lilbullet}]}\\
        \ev_A^B[f^{\lilbullet\lilbullet}]&=\ev_{A[f]}^{B[f^\lilbullet]}&\beta_\upsilon[f^{\lilbullet}]&=\beta_{\upsilon[f^{\lilbullet}]}\\
    \end{alignedat}$$ and: $$\begin{aligned}
    \funext_q[f]&=\funext_{q[f]}\\
    \beta^\Pi_q[f]&=\beta^\Pi_{q[f]}\\
    \eta^\Pi_q[f]&=\eta^\Pi_{q[f]}
    \end{aligned}$$ of Definition \ref{def:Semantics of axiomatic Pi-types} hold for every arrow $\Delta \xrightarrow{f}\Gamma$.

\proof
See Appendix \ref{appendixIII}.
\endproof
\end{proposition}

\subsection{Axiomatic $\boldsymbol{0}$-, $\boldsymbol{1}$-, $\boldsymbol{2}$-, and $\boldsymbol{\nat}$-types}

We specialise the notion of a display map 2-category with axiomatic $=$-types to obtain a version that incorporates a 2-dimensional semantic interpretation of axiomatic $0$-, $1$-, $2$-, and $\nat$-types. We observe the usual weakening from the intensional case \cite{MR2525957} to the very intensional, axiomatic one: retract bireflections are demoted to mere bireflections.

\begin{definition}[Semantics of axiomatic $0$-, $1$-, $2$-, $\nat$-types---categorical formulation]\label{def:Categorical semantics of axiomatic Nat-types}
    Let $(\boldsymbol{\mathcal{C}},\mathcal{D})$ be a display map 2-category endowed with axiomatic $=$-types. If $\Delta$ is an object of $\mathcal{C}$, we define $(\mathcal{D}/\Delta)^\dagger :=(\mathcal{D}/\Delta)\cup\{1_\Delta\}$.

    \begin{itemize}
        \item Let us assume that, for every object $\Gamma$, there is a choice of: \begin{center}
            a display map $\Gamma.0\to\Gamma$
        \end{center} in such a way that the pair: $$(\;1_{\Gamma.0}\;,\;1_{\circ}\;)$$ is a \textit{bireflection} of $\circ$ along the 2-functor $U:(\mathcal{D}/\Gamma.0)^\dagger\to (\boldsymbol{\mathcal{C}}/\Gamma.0)^0=\mathbf{1}$, being $\circ$ the unique object of the terminal 2-category $\mathbf{1}$. Let us assume that these data are stable under re-indexing.

        Unfolding this requirement, this means that, for every display map $P_C$ of codomain $\Gamma.0$, the functor: $$(U\blank)(1_\circ):(\mathcal{D}/\Gamma.0)^\dagger(1_{\Gamma.0},P_C)\to \mathbf{1}(\circ,\circ)$$ has a chosen pseudo-inverse functor with the corresponding natural isomorphisms between the two compositions and the respective identity functors.

        Then we say that $(\boldsymbol{\mathcal{C}},\mathcal{D})$ is \textit{endowed with axiomatic $0$-types}.
        
        \item Let us assume that, for every object $\Gamma$, there is a choice of: \begin{center}
            a display map $\Gamma.1\to\Gamma$ and a section $\star$ of $\Gamma.1\to\Gamma$
        \end{center} in such a way that the pair: $$(\;1_{\Gamma.1}\;,\;\star\xrightarrow{\star}1_{\Gamma.1}\;)$$ is a \textit{bireflection} of $\star$ along the 2-functor $U:(\mathcal{D}/\Gamma.1)^\dagger\to (\boldsymbol{\mathcal{C}}/\Gamma.1)^1$. Let us assume that these data are stable under re-indexing.

        Unfolding this requirement, this means that, for every display map $P_C$ of codomain $\Gamma.1$, the functor: $$(U\blank)(\star):(\mathcal{D}/\Gamma.1)^\dagger(1_{\Gamma.1},P_C)\to (\boldsymbol{\mathcal{C}}/\Gamma.1)(\star,P_C)$$ has a chosen pseudo-inverse functor with the corresponding natural isomorphisms $\beta^1,\eta^1$ between the two compositions and the respective identity functors.

        Then we say that $(\boldsymbol{\mathcal{C}},\mathcal{D})$ is \textit{endowed with axiomatic $1$-types}.
        

        \item Let us assume that, for every object $\Gamma$, there is a choice of: \begin{center}
            a display map $\Gamma.2\to\Gamma$ and a pair $(\bot,\top)$ of sections of $\Gamma.2\to\Gamma$
        \end{center} in such a way that the pair: $$(\;1_{\Gamma.2}\;,\;(\bot,\top)\xrightarrow{(\bot,\top)}(1_{\Gamma.2},1_{\Gamma.2})\;)$$ is a \textit{bireflection} of $(\bot,\top)$ along the 2-functor $U:(\mathcal{D}/\Gamma.1)^\dagger\to (\boldsymbol{\mathcal{C}}/\Gamma.1)^2$. Let us assume that these data are stable under re-indexing.

        Unfolding this requirement, this means that, for every display map $P_C$ of codomain $\Gamma.2$, the functor: $$(U\blank)(\bot,\top):(\mathcal{D}/\Gamma.2)^\dagger(1_{\Gamma.2},P_C)\to (\boldsymbol{\mathcal{C}}/\Gamma.2)^2((\bot,\top),(P_C,P_C))$$ has a chosen pseudo-inverse functor with the corresponding natural isomorphisms $\beta^2,\eta^2$ between the two compositions and the respective identity functors.

        Then we say that $(\boldsymbol{\mathcal{C}},\mathcal{D})$ is \textit{endowed with axiomatic $2$-types}.

        \item Let us assume that, for every object $\Gamma$, there is a choice of: \begin{center}
            a display map $\Gamma.\nat\to\Gamma$ a section $\0$ of $\Gamma.\nat\to\Gamma$ and an arrow $\s:\Gamma.\nat\to\Gamma.\nat$ over $\Gamma$
        \end{center} in such a way that the functor: $$(U\blank)(\0,\s P_C):(\mathcal{D}/\Gamma.\nat)^\dagger(1_{\Gamma.\nat},P_C)\to (\boldsymbol{\mathcal{C}}/\Gamma.\nat)^2((\0,\s P_C),(P_C,P_C))$$ is an equivalence, being $U$ the 2-functor $(\mathcal{D}/\Gamma.\nat)^\dagger\to (\boldsymbol{\mathcal{C}}/\Gamma.\nat)^2$. Let us assume that these data are stable under re-indexing.

        Unfolding this requirement, this means that, for every display map $P_C$ of codomain $\Gamma.\nat$, the functor: $$(U\blank)(\0,\s P_C):(\mathcal{D}/\Gamma.\nat)^\dagger(1_{\Gamma.\nat},P_C)\to (\boldsymbol{\mathcal{C}}/\Gamma.\nat)^2((\0,\s P_C),(P_C,P_C))$$ has a chosen pseudo-inverse functor with the corresponding natural isomorphisms $\beta^\nat,\eta^\nat$ between the two compositions and the respective identity functors.

        Then we say that $(\boldsymbol{\mathcal{C}},\mathcal{D})$ is \textit{endowed with axiomatic $\nat$-types}.
    \end{itemize}
\end{definition}

Let $(\boldsymbol{\mathcal{C}},\mathcal{D})$ be a display map 2-category endowed with axiomatic $=$-types.

Let us assume that the display map 2-category $(\boldsymbol{\mathcal{C}},\mathcal{D})$ is endowed with axiomatic $0$-types (resp. axiomatic $1$-types, axiomatic $2$-types, axiomatic $\nat$-types)---Definition \ref{def:Categorical semantics of axiomatic Nat-types}---and let $\Gamma$ be a context. As usual, we can define appropriate choice functions \textbf{Form Rule, Elim Rule, for $\boldsymbol{0}$-types} and \textbf{Form Rule, Intro Rule, Elim Rule, and Comp Axiom, for $\boldsymbol{1}$-types, $\boldsymbol{2}$-types, $\boldsymbol{\nat}$-types}---as in Definition \ref{def:Semantics of axiomatic Nat-types}---using the endowment of the structure $(\boldsymbol{\mathcal{C}},\mathcal{D})$ with axiomatic $0$-types, axiomatic $1$-types, axiomatic $2$-types, axiomatic $\nat$-types. We provide a proof sketch of this fact.

\medskip

\noindent\textbf{Form Rule for $\boldsymbol{0}$-types, $\boldsymbol{1}$-types, $\boldsymbol{2}$-types, $\boldsymbol{\nat}$-types} and \textbf{Intro Rule for $\boldsymbol{1}$-types, $\boldsymbol{2}$-types} are already defined.

\medskip

\noindent\textbf{Intro Rule for $\boldsymbol{\nat}$-types.} We define $\successivo$ as the unique section of $P_{\nat^\liltriangle}$ whose post-composition via $P_\nat^\lilbullet$ is $\s$.

\medskip

\noindent\textbf{Elim Rule for $\boldsymbol{0}$-types, $\boldsymbol{1}$-types, $\boldsymbol{2}$-types, $\boldsymbol{\nat}$-types.} If $P_C$ is an arrow of codomain $\Gamma.0$, then we define $\ind^0$ as the image of $1_\circ$ via the chosen pseudo-inverse of $(U\blank)(1_\circ)$. If $P_C$ is an arrow of codomain $\Gamma.1$ and $c$ is a section of $P_{C[\star]}$, then we define $\ind^1_c$ as the image of $\star^\lilbullet c$ via the chosen pseudo-inverse of $(U\blank)(\star)$. If $P_C$ is an arrow of codomain $\Gamma.2$ and $c$ is a section of $P_{C[\bot]}$ and $d$ is a section of $P_{C[\top]}$, then we define $\ind^2_{c,d}$ as the image of $(\bot^\lilbullet c,\top^\lilbullet d)$ via the chosen pseudo-inverse of $(U\blank)(\bot,\top)$. If $P_C$ is an arrow of codomain $\Gamma.\nat$ and $c$ is a section of $P_{C[\0]}$ and $d$ is a section of $P_{C[\s P_C]}$, then we define $\ind^\nat_{c,d}$ as the image of $(\0^\lilbullet c,(\s P_C)^\lilbullet d)$ via the chosen pseudo-inverse of $(U\blank)(\0,\s P_C)$.

\medskip

\noindent\textbf{Comp Axiom for $\boldsymbol{1}$-types, $\boldsymbol{2}$-types, $\boldsymbol{\nat}$-types.} We have a chosen 2-cell $(U\ind^1_c)(\star) \Rightarrow \star^\lilbullet c$, i.e. $\ind^1_c \star \Rightarrow \star^\lilbullet c$, over $\Gamma.1$. Therefore, referring to the diagram: 
\begin{equation}\label{reindexingone}\begin{tikzcd}[column sep=scriptsize,row sep=small]
	\Gamma && \Gamma \\
	\\
	{\Gamma.C[\star]} && {\Gamma.1.C} \\
	\\
	\Gamma && {\Gamma.1}
	\arrow[equals, from=1-1, to=1-3]
	\arrow["{\ind^1_c[\star]}"', shift right=3, from=1-1, to=3-1]
	\arrow["c", shift left=3, from=1-1, to=3-1]
	\arrow[""{name=0, anchor=center, inner sep=0}, "{\star^\lilbullet c}", shift left=3, from=1-3, to=3-3]
	\arrow[""{name=1, anchor=center, inner sep=0}, "{\ind^1_c\star}"', shift right=3, from=1-3, to=3-3]
	\arrow["{\star^\lilbullet}"{description}, from=3-1, to=3-3]
	\arrow["{P_{C[\star]}}"', from=3-1, to=5-1]
	\arrow["\lrcorner"{anchor=center, pos=0.125}, draw=none, from=3-1, to=5-3]
	\arrow["{P_C}", from=3-3, to=5-3]
	\arrow["{\star}"{description}, from=5-1, to=5-3]
	\arrow[shorten <=2pt, shorten >=2pt, Rightarrow, from=1, to=0]
\end{tikzcd}\end{equation} by the 2-universal property of $\Gamma.C[\star]$ there is unique a 2-cell:
\[\begin{tikzcd}[sep=small]
	{\Gamma} &&&& {\Gamma}
	\arrow[""{name=0, anchor=center, inner sep=0}, "c"', curve={height=12pt}, from=1-1, to=1-5]
	\arrow[""{name=1, anchor=center, inner sep=0}, "{\ind^1_c[\star]}", curve={height=-12pt}, from=1-1, to=1-5]
	\arrow["{\;b^1_c}", shorten <=5pt, shorten >=5pt, Rightarrow, from=1, to=0]
\end{tikzcd}\hspace{-0.4em}{.C[\star]}\] such that $P_{C[\star]}b^1_c= 1_{1_{\Gamma}}$ and $(\ind^1_c \star \Rightarrow \star^\lilbullet c)=(\star)^\lilbullet b^1_c$. As for the paragraphs \textbf{Comp Axiom for $\boldsymbol{=}$-types} and \textbf{Comp Axiom for $\boldsymbol{\Sigma}$-types}, we define $\beta^1_c$ as the unique section: $$\{b^1_c\}:\Gamma.1\to\Gamma.1.\id_{C[\star]}[\ind^1_c[\star];c]$$ of $P_{\id_{C[\star]}[\ind^1_c[\star];c]}$ such that (\ref{whatcharacterisescurlyparentheses}) is satisfied---where $A$ is $C[\star]$, $a$ and $b$ are $\ind^1_c[\star]$ and $c$ respectively, and $p$ is $b^1_c$. The definition of the sections $\beta^{2,\bot}_{c,d}$, $\beta^{2,\top}_{c,d}$, $\beta^{\nat,\0}_{c,d}$, $\beta^{\nat,\s}_{c,d}$ is completely analogous.

\begin{proposition}\label{stability of nat in 2d model}
    Let $(\boldsymbol{\mathcal{C}},\mathcal{D})$ be a display map 2-category endowed with axiomatic $0$-types (resp. axiomatic $1$-types, axiomatic $2$-types, axiomatic $\nat$-types). Referring to the data defined in paragraphs \emph{\textbf{Form Rule, Elim Rule, for $\boldsymbol{0}$-types}} and \emph{\textbf{Form Rule, Intro Rule, Elim Rule, and Comp Axiom, for $\boldsymbol{1}$-types, $\boldsymbol{2}$-types, $\boldsymbol{\nat}$-types}}, the \emph{stability conditions} of Definition \ref{def:Semantics of axiomatic Nat-types} hold for every arrow $\Delta \xrightarrow{f}\Gamma$.

\proof
See Appendix \ref{appendixIII}.
\endproof
\end{proposition}

\medskip

By \Cref{stability of id in 2d model}, \Cref{stability of sigma in 2d model}, \Cref{stability of pi in 2d model}, and \Cref{stability of nat in 2d model}, we infer the following:

\begin{theorem}\label{from a 2-dimensional model to a 1-dimensional one}
    Every \emph{display map 2-category} $(\boldsymbol{\mathcal{C}},\mathcal{D})$ endowed with axiomatic $=$-types, axiomatic $\Sigma$-types, axiomatic $\Pi$-types \& axiomatic function extensionality, and axiomatic $0$-, $1$-, $2$-, $\nat$-types---Definitions \ref{def:Categorical semantics of axiomatic =-types}, \ref{def:Categorical semantics of axiomatic Sigma-types}, \ref{def:Categorical semantics of axiomatic Pi-types}, and \ref{def:Categorical semantics of axiomatic Nat-types}---induces a \textnormal{display map category} $(\mathbf{C},\mathcal{D})$ endowed with axiomatic $=$-types, axiomatic $\Sigma$-types, axiomatic $\Pi$-types \& axiomatic function extensionality, and axiomatic $0$-, $1$-, $2$-, $\nat$-types---Definitions \ref{def:Semantics of axiomatic =-types}, \ref{def:Semantics of axiomatic Sigma-types}, \ref{def:Semantics of axiomatic Pi-types}, and \ref{def:Semantics of axiomatic Nat-types}---i.e. a model of \emph{ATT}.
\end{theorem}

\noindent By combining \Cref{from a 2-dimensional model to a 1-dimensional one} with \Cref{soundness property}, we deduce that an interpretation of ATT in any display map 2-category---endowed with axiomatic $=$-types, axiomatic $\Sigma$-types, axiomatic $\Pi$-types \& axiomatic function extensionality, and axiomatic $0$-, $1$-, $2$-, $\nat$-types---is well-defined and sound.

\subsection{The role of the cloven isofibration structure}\label{subsection:The role of the cloven isofibration structure}

Let $(\boldsymbol{\mathcal{C}},\mathcal{D})$ be a display map 2-category equipped with axiomatic $=$-types and axiomatic $\Sigma$-types (e.g.). As proven in the paragraph \textbf{Elim Rule for $=$-types}, the structure of axiomatic $=$-types associated with the display maps in the class $\mathcal{D}$ enables the construction of an \textit{elimination pseudo-term} $\tilde\j_c$ \textit{for $\boldsymbol{=}$-types}---for every section $c$ of the display map $P_{C[\refl_A]}$. This is not a genuine section of $P_C$, as the composition $P_C\tilde\j_c$ equals the identity of $\Gamma.A.A^\liltriangle.\id_A$ merely up to the 2-cell $\varphi_A$, making it the interpretation of a pseudo-term of type $C$. Similarly, the structure of axiomatic $\Sigma$-types on display maps allows for the construction of an \textit{elimination pseudo-term} $\tilde\splitt_c$ \textit{for $\Sigma$-types}---for every section $c$ of the display map $P_{C[\pair_A^B]}$---as shown in the paragraph \textbf{Elim Rule for $\boldsymbol{\Sigma}$-types}. Once again, this is the interpretation of a pseudo-term in the sense that the composition $P_C\tilde\splitt_c$ is the identity of $\Gamma.\sigmad_A^B$ only up to the 2-cell $\eta_A^B$.

The structure of cloven isofibrations on display maps plays a crucial role at this point, as it enables us to “strictify” the pseudo-terms $\tilde\j_c$ and $\tilde\splitt_c$---by transporting them back along the context identity proofs $\varphi_A$ and $\eta_A^B$---into the genuine elimination terms $\j_c$ and $\splitt_c$, respectively. In other words, these become actual sections of the corresponding display maps involved, at the cost of introducing additional 2-cells $\j_c \Rightarrow \tilde\j_c$ and $\splitt_c \Rightarrow \tilde\splitt_c$.

Paragraphs \textbf{Comp Axiom for $\boldsymbol{=}$-types} and \textbf{$\boldsymbol{\Sigma}$-types} show how, essentially, these two 2-cells represent the respective computation axioms: they lead, through re-indexing and further exploiting the arrow object structure---this time the one of the display maps $P_{C[\refl_A]}$ and $P_{C[\pair_A^B]}$ respectively---to the identification of sections $\h_c$ and $\sigma_c$ of the display maps associated to the types: $$\id_{C[\refl_A]}[\j_c[\refl_A];c]\text{\quad and \quad}\id_{C[\pair_A^B]}[\splitt_c[\pair_A^B];c]$$ respectively. These sections, therefore, provide interpretations for the terms of the computation axioms for $=$-types and $\Sigma$-types, respectively.

As mentioned above, since the isofibration structure on the display maps is cloven but not necessarily normal, this ensures that the computation axioms are satisfied, but not necessarily that the computation rules are. In other words, $\j_c[\refl_A]$ and $\splitt_c[\pair_A^B]$ do not, in general, coincide with the respective term $c$. To clarify this point, let us start by briefly analysing what happens if display maps in $(\boldsymbol{\mathcal{C}},\mathcal{D})$ are normal, as isofibrations. Referring to the diagram (\ref{reindexingid}), if $P_C$ is normal, then: \[\begin{alignedat}{3}
    \j_c\refl_A&=t_{\tilde\j_c\refl_A}^{\varphi_A*\refl_A}&&=t_{\tilde\j_c\refl_A}^{1_{\refl_A}}&&=\tilde\j_c\refl_A\\
    (\j_c \Rightarrow \tilde\j_c)*\refl_A&=\tau_{\tilde\j_c\refl_A}^{\varphi_A*\refl_A}&&=\tau_{\tilde\j_c\refl_A}^{1_{\refl_A}}&&= 1_{\tilde\j_c\refl_A}
\end{alignedat}\] implying that $\j_c[\refl_A]$ is in fact $c$ and that $h_c$ is the identity 1-cell of $c$, respectively. The latter condition, in turn, implies that $\h_c$ is $\reflex_{C[\refl_A]}[c]$. Similarly, referring to the diagram (\ref{reindexingsigma}), if $P_C$ is normal and, additionally, the 2-cell $\beta_A^B$ is the identity 2-cell---hence $(\pair_A^B)^\lilbullet_c$ coincides with $\tilde\splitt_c\pair_A^B$---then: \[\begin{alignedat}{3}
    \splitt_c\pair_A^B&=t_{\tilde\splitt_c\pair_A^B}^{\eta_A^B*\pair_A^B}&&=t_{\tilde\splitt_c\pair_A^B}^{1_{\pair_A^B}}&&=\tilde\splitt_c\pair_A^B\\
    (\splitt_c \Rightarrow \tilde\splitt_c)*\pair_A^B&=\tau_{\tilde\splitt_c\pair_A^B}^{\eta_A^B*\pair_A^B}&&=\tau_{\tilde\splitt_c\pair_A^B}^{1_{\pair_A^B}}&&= 1_{\tilde\splitt_c\pair_A^B}
\end{alignedat}\] implying, as before, that $\splitt_c[\pair_A^B]$ is $c$ and that $\sigma_c$ is $\reflex_{C[\pair_A^B]}[c]$. In this case, we conclude that the display map 2-category $(\boldsymbol{\mathcal{C}},\mathcal{D})$ is a model of the intensional version of these type formers---and similar considerations apply to the other type formers. In general, we have the following:
\begin{theorem}
    Let $(\boldsymbol{\mathcal{C}},\mathcal{D})$ be a display map 2-category endowed with axiomatic $0$-, $1$-, $2$-, $\nat$-, $=$-, $\Sigma$-, $\Pi$-types, and function extensionality---Definitions \ref{def:Categorical semantics of axiomatic =-types}, \ref{def:Categorical semantics of axiomatic Sigma-types}, \ref{def:Categorical semantics of axiomatic Pi-types}, and \ref{def:Categorical semantics of axiomatic Nat-types}. Let us assume that:\begin{itemize}
        \item the (families of) 2-cells $\beta^B_A,\;\beta^{A,B},\;\beta^1,\;\beta^2,\;\beta^\nat$ are identity 2-cells, and
        \item the display maps in $\mathcal{D}$ are normal, as isofibrations.
    \end{itemize} Then $(\boldsymbol{\mathcal{C}},\mathcal{D})$ (induces a display map category $(\mathbf{C},\mathcal{D})$ which) is a model of ITT.
\end{theorem} However, in the general case where display maps are \textit{not} required to be normal as isofibrations, but only cloven, this situation is not necessary: this is what we show in the next section.

\section{Revisiting the groupoid model}\label{sec:Revisiting the groupoid model}

The groupoid model of ITT, \cite{HofmannStreicher-Thegroupoidmodelref,MR1686862,MR4442637}, was devised by Hofmann and Streicher as a model of ITT that does not satisfy the uniqueness of identity proofs rule, providing a semantic proof that the latter is not admissible in ITT. It is considered a prelude towards the development of homotopy type theory. Additionally, it can be considered the quintessential example of a model under the notion of semantics introduced by Garner \cite{MR2525957}, as every other model according to this notion of semantics can, in a certain sense, be viewed as a generalisation of it. Notably, the display maps of the groupoid model formulated as a display map 2-category are in fact normal isofibrations.

The groupoid model, celebrated for its fundamental and illuminating nature, has inspired other authors to propose modifications aimed at developing models for generalised versions of dependent type theories. Notably, North \cite{MR4043275} and, from a different but related perspective, Altenkirch and Neumann \cite{altenkirchneumann} considered a weakening of this structure based on demoting groupoids to mere categories, thereby obtaining a foundational model for directed type theory. Our approach follows a similar trajectory. We propose a weakening of the groupoid model with the aim of demoting display maps to mere cloven isofibrations. However, we ensure to choose them in such a way that the re-indexing choices remain split, thereby preserving the integrity of an authentic model.

Let $\grpd$ be the (2,1)-category of groupoids, functors, and natural transformations---i.e. natural isomorphisms---with a chosen groupoid $1$ made of one object and one morphism, as a specified 2-terminal object of $\grpd$. Let $\mathcal{D}$ be the class of 1-cells whose elements are the functors of the form $P_A:\Gamma.A\to\Gamma$ obtained by applying the Grothendieck construction to the \textit{pseudofunctors} $(A,\phi^A,\psi^A):\Gamma \to \grpd$ where $\Gamma$ is a groupoid, and $\phi^A$ and $\psi^A$ are the \textit{coherent} families of natural isomorphisms $\phi^A_{{p},{q}}: A_qA_p\Rightarrow A_{qp}$ and $\psi_\gamma^A:A_{1_{\gamma}}\Rightarrow1_{A_\gamma}$---where $\gamma$ is any object of $\Gamma$ and $p$ and $q$ any composable arrows of $\Gamma$---making $A$ into a pseudofunctor. Hence $\Gamma.A$ is the groupoid whose objects are the pairs $(\gamma,x)$ where $\gamma$ is an object of $\Gamma$ and $x$ is an object of $A_\gamma$, and the arrows $(\gamma,x)$ to $(\gamma',x')$ are the pairs $({p_1},{p_2})$ where ${p_1}$ is an arrow $\gamma\to\gamma'$ and ${p_2}$ is an arrow $A_{p_1} x\to x'$. Composition of arrows, identity arrows, and inverses to arrows, are defined using the families $\phi^A$ and $\psi^A$: associativity, unitalities, and invertibilities follow from their coherence laws.\footnotemark[4] \footnotetext[4]{See Appendix \ref{sec:further details on the weakened groupoid model} for additional details.} The functor $P_A$ is the projection on the first component.

\medskip

\noindent\textbf{Cloven isofibration structure on display maps.} If we are given a 2-cell $\pi$ of the form: 
\[\begin{tikzcd}[sep=tiny]
	\Delta &&&& {\Gamma.A} \\
	&& {} & {} & {} \\
	&&&& \Gamma
	\arrow["{g\,=\,(g_1,g_2)}", from=1-1, to=1-5]
	\arrow["f"', curve={height=12pt}, from=1-1, to=3-5]
	\arrow["{P_A}", from=1-5, to=3-5]
	\arrow["\pi"{pos=0.4}, shift right, shorten <=5pt, shorten >=10pt, Rightarrow, from=2-3, to=2-5]
\end{tikzcd}\] then the mappings: \begin{center} $\delta \mapsto (f\delta\;,\; A_{\pi_\delta^{-1}}g_2\delta)$\;\;\;\;\;\;and
$(\delta\xrightarrow{p}\delta') \mapsto (f p\;,\; A_{f p}A_{\pi_\delta^{-1}}g_2\delta \xrightarrow{\phi^A}A_{(f p)\pi_\delta^{-1}}g_2\delta=A_{\pi_{\delta'}^{-1}(g_1p)}g_2\delta\xrightarrow{(\phi^A)^{-1}}A_{\pi_{\delta'}^{-1}}A_{(g_1p)}g_2\delta\xrightarrow{A_{\pi_{\delta'}^{-1}}g_2p}A_{\pi_{\delta'}^{-1}}g_2\delta' )$
\end{center} define a functor $\t_{g}^\pi:\Delta \to \Gamma.A$ whose post-composition via $P_A$ is $f$. Moreover, the arrow: $$(\pi_\delta\;,\; A_{\pi_\delta}A_{\pi_\delta^{-1}}g_2\delta\xrightarrow{\phi^A}A_{1_{g_1\delta}}g_2\delta\xrightarrow{\psi^A}g_2\delta)$$ from $(f\delta,A_{\pi_\delta^{-1}}g_2\delta)$ to $(g_1\delta,g_2\delta)$ is the $\delta$-component of a natural isomorphism $\tau_g^\pi:\t_{g}^\pi\Rightarrow g$ whose post-composition via $P_A$ is $\pi$.

These choices of $\t_g^\pi$ and $\tau_g^\pi$ define a cloven isofibration structure on $P_A$, which in general is \textit{not} normal: if $\pi$ is the identity and hence $f$ coincides with $g_1$, then the mapping $\delta \mapsto (f\delta\;,\; A_{\pi_\delta^{-1}}g_2\delta)$ coincides with $\delta \mapsto (g_1\delta\;,\; A_{1_{g_1\delta}}g_2\delta)$ which in general is different from $(g_1\delta,g_2\delta)$, since the pseudofunctor $A$ can be non-strict, hence in general $\t_g^\pi$ is different form $g$. Therefore, the takeaway is that expanding the class of semantic types to encompass all pseudofunctors into the 2-category 
$\grpd$---not just the strict functors used as semantic types in the original groupoid model---allows the class $\mathcal{D}$ of display maps to include cloven isofibrations that are not necessarily normal.

\medskip

\noindent\textbf{Re-indexing structure and arrow object structure on display maps.} Without delving too deeply into the details\footnotemark[4], we state that the 1-cells of $\mathcal{D}$ can be endowed with a 1- and 2-dimensional re-indexing structure, as well as with a structure of arrow objects, such that the compatibilities outlined in the third and fourth points of Definition \ref{split display map 2-category} are satisfied. We only mention that the re-indexing of a display map $P_A:\Gamma.A\to\Gamma$ along a 1-cell $f:\Delta\to\Gamma$ of $\grpd$ is obtained by applying the Grothendieck construction to the pseudofunctor $\begin{tikzcd}[cramped,column sep=small,row sep=tiny]
	\Delta && \Gamma && \grpd,
	\arrow["f"{description}, from=1-1, to=1-3]
	\arrow["A"{description}, from=1-3, to=1-5]
\end{tikzcd}$ where $A$ is the pseudofunctor inducing $P_A$. We also mention that $P_{\id_A}$ is induced by the (pseudo)functor $\id_A:\Gamma.A.A^\liltriangle\to\grpd$ mapping: $$\begin{aligned}
    (\gamma,x,y)&\mapsto A_\gamma(x,y)\\(p_1,p_2,p_3)&\mapsto(p\mapsto (x'\xrightarrow{p_2^{-1}}A_{p_1}x\xrightarrow{A_{p_1}p}A_{p_1}y\xrightarrow{p_3}y'))
\end{aligned}$$---where the homset $A_\gamma(x,y)$ is endowed with the trivial groupoid structure---and that the $(\gamma,x,y,p)$-component of the 2-cell $\alpha_A:{P_{A^\liltriangle}P_{\id_A}} \Rightarrow {P_A^\lilbullet P_{\id_A}}$ is the arrow: $$(\;1_\gamma\;,\;A_{1_\gamma}x\xrightarrow{\psi^A}x\xrightarrow{p}y\;):(\gamma,x)\to(\gamma,y).$$ Following the construction at paragraphs \textbf{Elim Rule} and \textbf{Comp Axiom for $=$-types}, one can now reconstruct the choice functions $\refl$, $\varphi$, and $\j$ and observe that $\j_c$ acts \textit{on objects} as the mapping: $$(\gamma,x,y,p)\mapsto(\gamma,x,y,p,C_{\varphi_A^{-1}}c(\gamma,y))$$ whenever $(\gamma,x)\mapsto(\gamma,x,c(\gamma,x))$ is the action on object of a section $c$ of $P_{C[\refl_A]}$, for some display map $P_C$ over $\Gamma.A.A^\liltriangle.\id_A$. This allows to conclude that $\j_c[\refl_A]$ acts on objects as the mapping: $$(\gamma,x)\mapsto(\gamma,x,C_{1_{(\gamma,x,x,1_x)}}c(\gamma,x))$$ and therefore, unless $C$ was chosen to be a strict functor $\Gamma.A.A^\liltriangle.\id_A\to\grpd$, in general $\j_c[\refl_A]$ and $c$ \textit{will not coincide}. We conclude that:

\begin{theorem}
    The pair $(\grpd,\mathcal{D})$ is a display map 2-category endowed with axiomatic $=$-type, which, as a display map category---see Theorem \ref{from a 2-dimensional model to a 1-dimensional one}--- is a model of axiomatic $=$-types that does not validate Comp Rule for $=$-types.

    In particular Comp Rule for $=$-types is not admissible in \emph{ATT}.
\end{theorem}

\section{Syntactic display map 2-category and completeness}\label{sec:completeness}

\Cref{from a 2-dimensional model to a 1-dimensional one} and \Cref{soundness property} show us that the class of display map 2-categories---endowed with axiomatic $0$-, $1$-, $2$-, $\nat$-, $=$-, $\Sigma$-, $\Pi$-types, and function extensionality---constitute a sound semantics of ATT. However, in light of the following proposition, this semantics cannot be complete, since the rule \ref{discreteness} is not admissible in ATT.

\begin{proposition}\label{discreteness valid} Every display map 2-category $(\boldsymbol{\mathcal{C}},\mathcal{D})$ endowed with axiomatic $=$-types is a model of the following: \begin{equation}\label{discreteness}
    \inferrule{\judge{\gamma:\Gamma}{A(\gamma):\type}}{\judge{\gamma:\Gamma;\;x,y:A(\gamma);\;p,q:x=y;\;\alpha:p=q}{p\equiv q}}
\tag{Disc}
\end{equation} \emph{discreteness} rule.

\proof
If $\Gamma$, $P_A$, $x$, $y$, $p$, $q$, $\alpha$ are the interpretations of: \begin{center} $\judgectx{\gamma:\Gamma}$\\$\judge{\gamma:\Gamma}{A(\gamma):\type}$\;\;\;\;\;\;\;\;\;\;\;\;\;\;\;\;$\judge{\gamma:\Gamma}{x(\gamma),y(\gamma):A(\gamma)}$\\$\judge{\gamma:\Gamma}{p(\gamma),q(\gamma):x(\gamma)=y(\gamma)}$\;\;\;\;\;\;\;\;\;\;\;\;\;\;\;\;$\judge{\gamma:\Gamma}{\alpha(\gamma):p(\gamma)=q(\gamma)}$
\end{center} respectively, we are left to prove that $(\boldsymbol{\mathcal{C}},\mathcal{D})$ models the term equality judgement: $$\judge{\gamma:\Gamma}{p(\gamma)\equiv q(\gamma)}$$ i.e. that $p$ and $q$ coincide as arrows of $\mathcal{C}$. We observe that, by post-composing: \begin{center}$(x;y)^\lilbullet p$ and $(x;y)^\lilbullet q$\end{center} by: $$P_{A^\liltriangle}P_{\id_A} \overset{\alpha_A}{\Longrightarrow}P_A^\lilbullet P_{\id_A}$$ we obtain two objects, $\alpha_A(x;y)^\lilbullet p$ and $\alpha_A(x;y)^\lilbullet q$, of $(\mathcal{C}/\Gamma)(\Gamma,\Gamma.A)^\to$ of the form: \begin{center}
    $P_{A^\liltriangle}P_{\id_A}(x;y)^\lilbullet p\Longrightarrow P_A^\lilbullet P_{\id_A}(x;y)^\lilbullet p$\\$P_{A^\liltriangle}P_{\id_A}(x;y)^\lilbullet q\Longrightarrow P_A^\lilbullet P_{\id_A}(x;y)^\lilbullet q$
\end{center} respectively, hence both of the form: $$x \Longrightarrow y$$ because $P_{A^\liltriangle}P_{\id_A}(x;y)^\lilbullet p=x=P_{A^\liltriangle}P_{\id_A}(x;y)^\lilbullet q$ and $P_A^\lilbullet P_{\id_A}(x;y)^\lilbullet p=y=P_A^\lilbullet P_{\id_A}(x;y)^\lilbullet q$.

Moreover the arrow: $$(x;y)^\lilbullet\alpha_{\id_A[x;y]}(p;q)^\lilbullet\alpha$$ of $(\mathcal{C}/\Gamma)(\Gamma,\Gamma.A.A^\liltriangle.\id_A)$ from $(x;y)^\lilbullet p$ to $(x;y)^\lilbullet q$ induces, by post-composition via the pair $(P_{A^\liltriangle}P_{\id_A},P_A^\lilbullet P_{\id_A})$, an arrow: $$(\;\;P_{A^\liltriangle}P_{\id_A}(x;y)^\lilbullet\alpha_{\id_A[x;y]}(p;q)^\lilbullet\alpha\;\;,\;\;P_{A}^\lilbullet P_{\id_A}(x;y)^\lilbullet\alpha_{\id_A[x;y]}(p;q)^\lilbullet\alpha\;\;)$$ from $\alpha_A(x;y)^\lilbullet p$ to $\alpha_A(x;y)^\lilbullet q$ in $(\mathcal{C}/\Gamma)(\Gamma,\Gamma.A)^\to$.

However, we observe that: \[\begin{aligned}
    P_{\id_A[x;y]}\alpha_{\id_A[x;y]}(p;q)^\lilbullet\alpha&=1_{P_{\id_A[x;y]}P_{\id_A[x;y]^\liltriangle}P_{\id_{\id_A[x;y]}}}(p;q)^\lilbullet\alpha\\
    &=1_{1_\Gamma}*P_{\id_A[x;y]}P_{\id_A[x;y]^\liltriangle}P_{\id_{\id_A[x;y]}}(p;q)^\lilbullet\alpha\\
    &=1_{1_\Gamma}*P_{\id_A[x;y]}P_{\id_A[x;y]^\liltriangle}(p;q)P_{\id_{\id_A[x;y]}[p;q]}\alpha\\
    &=1_{1_\Gamma}*1_\Gamma\\
    &=1_{1_\Gamma}
\end{aligned}\] and that: $$P_{A^\liltriangle}P_{\id_A}(x;y)^\lilbullet=xP_{\id_A[x;y]}\;\;\;\;\;\;\;\;\textnormal{and}\;\;\;\;\;\;\;\;P_A^\lilbullet P_{\id_A}(x;y)^\lilbullet=yP_{\id_A[x;y]}$$ therefore the pair: $$(\;\;P_{A^\liltriangle}P_{\id_A}(x;y)^\lilbullet\alpha_{\id_A[x;y]}(p;q)^\lilbullet\alpha\;\;,\;\;P_A^\lilbullet P_{\id_A}(x;y)^\lilbullet\alpha_{\id_A[x;y]}(p;q)^\lilbullet\alpha\;\;)$$ which we understood to be an arrow from $\alpha_A(x;y)^\lilbullet p$ to $\alpha_A(x;y)^\lilbullet q$ in $(\mathcal{C}/\Gamma)(\Gamma,\Gamma.A)^\to$, coincides with $(x*1_{1_\Gamma},y*1_{1_\Gamma})$ i.e. $(1_x,1_y)$. In other words, we proved that the diagram: 
\[\begin{tikzcd}[column sep=scriptsize]
	x && x \\
	\\
	y && y
	\arrow["{1_x}"{description}, Rightarrow, from=1-1, to=1-3]
	\arrow["{\alpha_A(x;y)^\lilbullet p}"', Rightarrow, from=1-1, to=3-1]
	\arrow["{\alpha_A(x;y)^\lilbullet q}", Rightarrow, from=1-3, to=3-3]
	\arrow["{1_y}"{description}, Rightarrow, from=3-1, to=3-3]
\end{tikzcd}\] commutes, hence $\alpha_A(x;y)^\lilbullet p=\alpha_A(x;y)^\lilbullet q$ which implies that $(x;y)^\lilbullet p=(x;y)^\lilbullet q$, by the universal property of $\alpha_A$. Being: 
\[\begin{tikzcd}[column sep=small]
	{\Gamma.\id_A[x,y]} && {\Gamma.A.A^\liltriangle.\id_A} \\
	\\
	\Gamma && {\Gamma.A.A^\liltriangle}
	\arrow["{(x;y)^\lilbullet}", from=1-1, to=1-3]
	\arrow[from=1-1, to=3-1]
	\arrow["\lrcorner"{anchor=center, pos=0.125}, draw=none, from=1-1, to=3-3]
	\arrow[from=1-3, to=3-3]
	\arrow["{x;y}", from=3-1, to=3-3]
\end{tikzcd}\] a pullback, we conclude that $p=q$ and we are done. \endproof
\end{proposition}

In particular, whenever $(\boldsymbol{\mathcal{C}},\mathcal{D})$ is a display map 2-category endowed with axiomatic $=$-types, axiomatic $\Sigma$-types, axiomatic $\Pi$-types \& axiomatic function extensionality, and axiomatic $0$-, $1$-, $2$-, $\nat$-types---see \Cref{def:Categorical semantics of axiomatic =-types,def:Categorical semantics of axiomatic Sigma-types,def:Categorical semantics of axiomatic Pi-types,def:Categorical semantics of axiomatic Nat-types}, respectively, we deduce the following refinement of \Cref{soundness property}. \begin{cor}[Soundness property] The interpretation of \textnormal{ATT $+$ \ref{discreteness}} in $(\boldsymbol{\mathcal{C}},\mathcal{D})$ is sound, that is:
\begin{itemize}
    \item whenever \textnormal{ATT} infers $\judgectx{\gamma:\Gamma}$, then its interpretation---let us indicate it as $\Gamma$---is defined;
    \item whenever \textnormal{ATT} infers $\judgectx{\gamma:\Gamma} \equiv \judgectx{\gamma':\Gamma'}$, for some contexts $\judgectx{\gamma:\Gamma}$ and $\judgectx{\gamma':\Gamma'}$, then their interpretations $\Gamma$ and $\Gamma'$ coincide;
    \item whenever \textnormal{ATT} infers $\judge{\gamma:\Gamma}{A(\gamma):\type}$, then its interpretation---let us indicate it as $P_A$---is defined;
    \item whenever \textnormal{ATT} infers $\judge{\gamma:\Gamma}{a(\gamma):A(\gamma)}$, then its interpretation---let us indicate it as $a:\Gamma\to\Gamma.A$---is defined;
    \item whenever \textnormal{ATT} infers $\judge{\gamma:\Gamma}{A(\gamma)\equiv A'(\gamma)}$ for some type judgements $\judge{\gamma:\Gamma}{A(\gamma):\type}$ and $\judge{\gamma:\Gamma}{A'(\gamma):\type}$, then their interpretations $P_A$ and $P_{A'}$ coincide;
    \item whenever \textnormal{ATT} infers $\judge{\gamma:\Gamma}{a(\gamma)\equiv a'(\gamma)}$ for some term judgements $\judge{\gamma:\Gamma}{a(\gamma):A(\gamma)}$ and $\judge{\gamma:\Gamma}{a'(\gamma):A(\gamma)}$ and some type judgement $\judge{\gamma:\Gamma}{A(\gamma):\type}$, then their interpretations $a$ and $a'$---as sections of $P_A$---coincide.
\end{itemize}
\end{cor}

From \cite[Propositions 3.2.5]{spadettothesis} we know that we can define a (2,1)-category $\boldsymbol{\mathcal{C}}_{\text{ATT $+$ \ref{discreteness}}}$ with a chosen 2-terminal object from ATT $+$ \ref{discreteness}, whose 0- and 1-cells are the contexts and the substitutions, respectively, of ATT $+$ \ref{discreteness}---identified up to renaming their free variables and up to componentwise equality judgement---and whose 2-cells: $$\judge{\gamma:\Gamma}{f(\gamma):\Delta} \quad\Rightarrow\quad \judge{\gamma:\Gamma}{g(\gamma):\Delta}\quad\quad\quad:\quad\quad\quad \judgectx{\gamma:\Gamma}\quad\to\quad \judgectx{\delta:\Delta}$$ are context propositional equalities i.e. lists, indicated as $\judge{\gamma:\Gamma}{p(\gamma):f(\gamma)=g(\gamma)}$, of term judgements of the form: 
\begin{itemize}
    \item[] $\judge{\gamma:\Gamma}{p_1(\gamma):f_1(\gamma)=g_1(\gamma)}$
    \item[] $\judge{\gamma:\Gamma}{p_2(\gamma):f_2(\gamma)=p_1(\gamma)^*g_2(\gamma)}$
    \item[] $\judge{\gamma:\Gamma}{p_3(\gamma):f_3(\gamma)=(p_1(\gamma),p_2(\gamma))^*g_3(\gamma)}$
    \item[] \quad ...
    \item[] $\judge{\gamma:\Gamma}{p_n(\gamma):f_n(\gamma)=(p_1(\gamma),...,p_{n-1}(\gamma))^*g_n(\gamma)}$
\end{itemize} where transport operations along multiple identity proofs are defined as explained in \cite[Subsection 3.2.I]{spadettothesis}. Additionally \cite[Propositions 3.3.1, 3.3.2, and 3.3.3]{spadettothesis} show how to provide $\boldsymbol{\mathcal{C}}_{\text{ATT $+$ \ref{discreteness}}}$ with the structure $\mathcal{D}$ of a display map 2-category. Finally \cite[Propositions 3.3.4, 3.3.5, 3.3.6, and 3.3.7]{spadettothesis} show that $(\boldsymbol{\mathcal{C}}_{\text{ATT $+$ \ref{discreteness}}},\mathcal{D})$ is endowed with axiomatic $=$- and $\Sigma$-types. Analogously, one can prove that it is also endowed with axiomatic $0$-, $1$-, $2$-, $\nat$-, $\Pi$-types, and function extensionality. Since the interpretation in $(\boldsymbol{\mathcal{C}}_{\text{ATT $+$ \ref{discreteness}}},\mathcal{D})$ of any judgement of ATT $+$ \ref{discreteness} is itself a (list of) derivable judgement of ATT $+$ \ref{discreteness}---as so is every 1-cell of $(\boldsymbol{\mathcal{C}}_{\text{ATT $+$ \ref{discreteness}}},\mathcal{D})$---we deduce the following: \begin{theorem}[Completeness property]
    The interpretation of the theory \emph{ATT $+$ \ref{discreteness}} in display map 2-categories---endowed with axiomatic $0$-, $1$-, $2$-, $\nat$-, $=$-, $\Sigma$-, $\Pi$-types, and function extensionality---is complete, that is: \begin{itemize}
        \item whenever the interpretation of a context, or of a type judgement, or of a term judgement in \emph{ATT $+$ \ref{discreteness}} is defined in every such display map 2-category, then \emph{ATT $+$ \ref{discreteness}} derives it, and
        \item any two derivable contexts, or types, or terms in \emph{ATT $+$ \ref{discreteness}}, whose interpretations coincide in every such display map 2-category, are provably definitionally equal in \emph{ATT $+$ \ref{discreteness}}.
    \end{itemize} 
\end{theorem}

\section{Conclusion and future work}\label{sec:conclusion and future work}

In this paper we provided a 2-categorical procedure to construct models of ATT. We applied this procedure to identify a weakening of the groupoid model that models ATT without believing Comp Rule for $=$-types. The semantics provided by this class of models, however, is not complete with respect to ATT. Similarly to Garner's notion of semantics for ITT, every model validates the discreteness rule---namely, every model believes that every type is a 1-type. This arises due to the arrow object structure in display maps---see the proof of \cref{discreteness valid}. This phenomenon occurs in both the original and our weakened version of the groupoid model, as well as in the versions of the category model proposed in \cite{MR4043275} and \cite{altenkirchneumann}---using directed $=$-types, or \textit{hom-types}, in place of $=$-types.

Hence it is natural to ask whether it is possible to interpret ATT in tricategories---according to the notion presented in \cite{MR3076451}---with a class of display maps in such a way that the interpretation is complete without adding the discreteness rule to ATT. In a display map 2-category, two parallel 2-cells are either distinct or coincide. However, in a \textit{display map 3-category}, we may talk about homotopies between two 2-cells, or propositional equalities between identity proofs, without necessarily collapsing provably identical identity proofs into the same 2-cell. Moreover, every identity proof between identity proofs is itself in particular an identity proof: this leads to the idea that the third dimension should not be explicitly used to define the interpretation, which, like for display map 2-categories, will only rely on 2-cells---to be later converted into sections using the arrow object. Yet, as mentioned, we will not have the issue of having to collapse every two parallel 2-cells that are not distinct into the same 2-cell. This is why we believe that the interpretation of ATT in appropriate display map 3-categories might be complete.

This work belongs to the field of research in the categorical semantics of dependent type theories, along with their variants and generalisations---such as ATT itself---and aims to advance one of the main objectives of the area: to establish a concept of a general semantics that is both broad enough to encompass highly general notions of dependent type theories and functional enough to be effectively applied in their study. We argue that this work can be contextualised within the bicategorical approach to the semantics of dependent type theory proposed by Ahrens, North, and van der Weide \cite{MR4537008,MR4673069}. This approach aims to provide a general framework for describing the semantics of the structural part of dependent type theory, as well as its generalisation towards directed type theory. Specifically, we believe that display map 2-categories constitute a particular instance of the notion of a \textit{display map bicategory} introduced in that work, under the condition that reductions are symmetric. In fact, in a display map 2-category, the 2-cells are independent of the notion of a $=$-type former and can, a priori, be regarded as \textit{reduction judgements}---although not directed: it is the semantic $=$-type former that enables us to interpret and convert them as identity proofs.

We are also interested in investigating whether Garner's idea of encapsulating the semantics of intensional $=$-types in arrow objects---applied in this article to encode the semantics of axiomatic $=$-types---could be adapted, through a suitable generalisation of our notion of display map 2-category, to encode the semantics of directed identity types---whether intensional or axiomatic---as proposed by North \cite{MR4043275} or by Altenkirch and Neumann \cite{altenkirchneumann}. This would aim to address a problem posed by Ahrens, North, and van der Weide \cite{MR4673069} themselves: extending the notion of comprehension bicategory to accommodate the interpretation of the hom-type former à la North.

More generally, we believe there is still much to explore regarding the semantics of various extensions and variations of dependent type theory. This is especially true in relation to the weakenings of the identity type constructor---both with respect to its elimination rule and its computation rule---that are currently prominent, as well as their interrelationships. We believe the notion of semantics introduced in this article provides a contribution to this important line of research in type theory.

\bigskip

\bigskip

\noindent\textbf{Acknowledgements.} This research was supported by a School of Mathematics full-time EPSRC Doctoral Training Partnership Studentship 2019/2020, by the Italian MUR PRIN 2022 “STENDHAL”, and by eOTP RECIPROG.

The author is grateful to his supervisors Nicola Gambino and Federico Olimpieri, as well as to Benedikt Ahrens, Dani\"el Otten, Marino Miculan, and Niels van der Weide for useful discussions on the subject.

\newpage

\appendix

\section{Type-checking of the stability conditions}\label{appendixI}

This section begins with the proof of the fact, mentioned in Subsection \ref{subsection:Encoding parallel terms into a substitution}, that in a display map category $(\mathbf{C},\mathcal{D})$ the square: \begin{equation}\label{reindexingprimareindexingdopo}\begin{tikzcd}[column sep=small,row sep=scriptsize]
	\Delta & \Gamma \\
	{\Delta.A[f].A[f]^{\liltriangle}} & {\Gamma.A.A^\liltriangle}
	\arrow["f", from=1-1, to=1-2]
	\arrow["{a[f];b[f]}"', from=1-1, to=2-1]
	\arrow["{a;b}", from=1-2, to=2-2]
	\arrow["{f^{\lilbullet\lilbullet}}", from=2-1, to=2-2]
\end{tikzcd}\end{equation} commutes if $f$ is a given arrow $\Delta\to\Gamma$ and $a$ and $b$ are sections of a given display map $P_A$.

By building the term $b[P_Aa]$, i.e. $b$ itself---display map categories are assumed to be split---according to its definition:
\[\begin{tikzcd}[column sep=small, row sep=scriptsize]
	\Gamma && {\Gamma.A} && \Gamma \\
	\\
	{\Gamma.A=\Gamma.A[P_Aa]} && {\Gamma.A.A^\liltriangle} && {\Gamma.A} \\
	{} \\
	\Gamma && {\Gamma.A} && \Gamma
	\arrow["a", from=5-1, to=5-3]
	\arrow["{P_A}", from=5-3, to=5-5]
	\arrow["a", from=1-1, to=1-3]
	\arrow["{P_A}", from=1-3, to=1-5]
	\arrow["{P_A}"{description}, from=3-5, to=5-5]
	\arrow["{P_A\,=\,P_{A[P_Aa]}}"{description}, from=3-1, to=5-1]
	\arrow["b"{description}, from=1-5, to=3-5]
	\arrow["{b[P_A]}"{description}, from=1-3, to=3-3]
	\arrow["{b\,=\,b[P_Aa]}"{description}, from=1-1, to=3-1]
	\arrow["{a^\lilbullet}", from=3-1, to=3-3]
	\arrow["{P_A^\lilbullet}", from=3-3, to=3-5]
	\arrow["{P_{A^\liltriangle}}"{description}, from=3-3, to=5-3]
	\arrow["\lrcorner"{anchor=center, pos=0.125}, draw=none, from=3-3, to=5-5]
	\arrow["\lrcorner"{anchor=center, pos=0.125}, draw=none, from=3-1, to=5-3]
\end{tikzcd}\] we infer that the upper left-hand square commutes. And since: \begin{center}
    $P_{A^\liltriangle}b[P_A]a=a$ \quad and \quad $P_A^\lilbullet a^\lilbullet b=(P_Aa)^\lilbullet b=b$
\end{center} we conclude the following two presentations: $$\begin{aligned}
    a;b&=(\begin{tikzcd}[column sep=large]
	\Gamma & {\Gamma.A} & {\Gamma.A.A^\liltriangle}
	\arrow["b"{description}, from=1-1, to=1-2]
	\arrow["{a^\lilbullet}"{description}, from=1-2, to=1-3]\end{tikzcd})\\
    &=(\begin{tikzcd}[column sep=large]
	\Gamma & {\Gamma.A} & {\Gamma.A.A^\liltriangle}
	\arrow["a"{description}, from=1-1, to=1-2]
	\arrow["{b[P_A]}"{description}, from=1-2, to=1-3]
\end{tikzcd})
\end{aligned}$$ of the arrow $a;b$.

In order to verify that the diagram (\ref{reindexingprimareindexingdopo}) commutes, we use the first presentation: $$\Delta\xrightarrow{b[f]}\Delta.A[f]\xrightarrow{a[f]^\lilbullet}\Delta.A[f].A[f]^\liltriangle$$ for $a[f];b[f]$, and the second one: $$\Gamma\xrightarrow{a}\Gamma.A\xrightarrow{b[P_A]}\Gamma.A.A^\liltriangle$$ for $a;b$. We are left to verify that: \[\begin{tikzcd}[column sep=tiny]
	\Delta &&&& \Gamma && {\Gamma.A} \\
	\\
	{\Delta.A[f]} &&& {\Delta.A[f].A[f]^\liltriangle} & {\Delta.A[f].A^\liltriangle[f^\lilbullet]} && {\Gamma.A.A^\liltriangle}
	\arrow["{b[P_A]}"{description}, from=1-7, to=3-7]
	\arrow["{f^{\lilbullet\lilbullet}}", from=3-5, to=3-7]
	\arrow[Rightarrow, no head, from=3-4, to=3-5]
	\arrow["{b[P_A(af)]}"{description}, from=1-1, to=3-1]
	\arrow["{a[f]^\lilbullet}", from=3-1, to=3-4]
	\arrow["f", from=1-1, to=1-5]
	\arrow["a", from=1-5, to=1-7]
\end{tikzcd}\] commutes, since $b[P_A(af)]=b[f]$. Now, since the diagram: \[\begin{tikzcd}
	\Delta &&& \Gamma &&& {\Gamma.A} \\
	\\
	{\Delta.A[f]} &&& {\Gamma.A} &&& {\Gamma.A.A^\liltriangle}
	\arrow["{b[P_A]}"{description}, from=1-7, to=3-7]
	\arrow["{b[P_A(af)]}"{description}, from=1-1, to=3-1]
	\arrow["f", from=1-1, to=1-4]
	\arrow["a", from=1-4, to=1-7]
	\arrow["{f^\lilbullet}", from=3-1, to=3-4]
	\arrow["{a^\lilbullet}", from=3-4, to=3-7]
\end{tikzcd}\] commutes---and this is true because the equality $a^\lilbullet f^\lilbullet=(af)^\lilbullet$ holds---we are done if we verify that: \[\begin{tikzcd}[sep=scriptsize]
	{\Delta.A[f]} &&& {\Gamma.A} \\
	\\
	{\Delta.A[f].A[f]^\liltriangle} &&& {\Gamma.A.A^\liltriangle}
	\arrow["{a^\lilbullet}"{description}, from=1-4, to=3-4]
	\arrow["{f^\lilbullet}", from=1-1, to=1-4]
	\arrow["{f^{\lilbullet\lilbullet}}", from=3-1, to=3-4]
	\arrow["{a[f]^\lilbullet}"{description}, from=1-1, to=3-1]
\end{tikzcd}\] commutes. This is actually true as: $$f^{\lilbullet\lilbullet}a[f]^\lilbullet=(f^\lilbullet a[f])^\lilbullet=(af)^\lilbullet=a^\lilbullet f^\lilbullet$$ hence we are done.

\medskip

In what follows, we provide full verification of the type-checking of the stability conditions of \Cref{def:Semantics of axiomatic =-types,def:Semantics of axiomatic Sigma-types}. Verification is completely analogous for \Cref{def:Semantics of axiomatic Pi-types} and \Cref{def:Semantics of axiomatic Nat-types}.

\medskip

\noindent\textbf{Type-checking of the stability conditions of Definition \ref{def:Semantics of axiomatic =-types}}

Let us assume that $(\mathbf{C},\mathcal{D})$ is endowed with the choice functions $\id$, $\reflex$, $\j$, $\h$ as in \textit{Form Rule}, \textit{Intro Rule}, \textit{Elim Rule}, \textit{Comp Axiom} of Definition \ref{def:Semantics of axiomatic =-types}. Let $f$ be an arrow $\Delta\to\Gamma$ and let $P_A$ be a display map of codomain $\Gamma$.

\medskip

\noindent\textbf{Stability of $\boldsymbol{\id}$.} The codomain of $P_{\id_{A[f]}}$ is by definition $\Delta.A[f].A[f]^\liltriangle$. We observe that the codomain of $P_{\id_A[f^{\lilbullet\lilbullet}]}$ is $\Delta.A[f].A^\liltriangle[f^\lilbullet]$, which coincides with $\Delta.A[f].A[f]^\liltriangle$ since: $$A^\liltriangle[f^\lilbullet]=A[P_A][f.A]=A[f][P_{A[f]}]=A[f]^\liltriangle$$ being $(\mathbf{C},\mathcal{D})$ is split. As $P_{\id_{A[f]}}$ and $P_{\id_A[f^{\lilbullet\lilbullet}]}$ have the same codomain, the equality: $$\id_A[f^{\lilbullet\lilbullet}]=\id_{A[f]}$$ is meaningful. From now on, let us assume that it is satisfied.

\medskip

\noindent\textbf{Stability of $\boldsymbol{\reflex}$.} By definition $\reflex_{A[f]}$ is a section of the display map $\Delta.A[f].\id_{A[f]}[\delta_{A[f]}]\to \Delta.A[f]$. We observe that $\reflex_A[f^\lilbullet]$ is a section of the display map $\Delta.A[f].\id_{A}[\delta_A][f^\lilbullet]\to \Delta.A[f]$, and: $$\id_{A}[\delta_A][f^\lilbullet]=\id_{A}[f^{\lilbullet\lilbullet}][\delta_{A[f]}]=\id_{A[f]}[\delta_{A[f]}].$$ Therefore $\reflex_A[f^\lilbullet]$ is itself a section of $\Delta.A[f].\id_{A[f]}[\delta_{A[f]}]\to\Delta.A[f]$, hence the equality: $$\reflex_A[f^\lilbullet]=\reflex_{A[f]}$$ is meaningful. From now on, let us assume that it is satisfied. Observe that this implies that the square:
\[\begin{tikzcd}[column sep=small,row sep=scriptsize]
	{\Delta.A[f]} && {\Gamma.A} \\
	\\
	{\Delta.A[f].A[f]^\liltriangle.\id_{A[f]}} && {\Delta.A.A^\liltriangle.\id_{A}}
	\arrow["{f^\lilbullet}", from=1-1, to=1-3]
	\arrow["{\refl_{A[f]}}"{description}, from=1-1, to=3-1]
	\arrow["{\refl_A}"{description}, from=1-3, to=3-3]
	\arrow["{f^{\lilbullet\lilbullet\lilbullet}}", from=3-1, to=3-3]
\end{tikzcd}\]
commutes.

\medskip

\noindent\textbf{Stability of $\boldsymbol{\j}$.} Let $P_C$ be a display map of codomain $\Gamma.A.A^\liltriangle.\id_A$ and let $c$ be a section of $\Gamma.A.C[\refl_A]\to\Gamma.A$. As $\j_c$ is by definition a section of $P_C$, then $\j_c[f^{\lilbullet\lilbullet\lilbullet}]$ is a section of $P_{C[f^{\lilbullet\lilbullet\lilbullet}]}$. Now, we observe that $P_{C[f^{\lilbullet\lilbullet\lilbullet}]}$ is a display map of codomain: $$\Delta.A[f].A^\liltriangle[f^\lilbullet].\id_A[f^{\lilbullet\lilbullet}]=\Delta.A[f].A[f]^\liltriangle.\id_{A[f]}$$ and $c[f^\lilbullet]$ is a section of $\Delta.A[f].C[\refl_A][f^\lilbullet]\to\Delta.A[f]$ i.e. of: $$\Delta.A[f].C[f^{\lilbullet\lilbullet\lilbullet}][\refl_{A[f]}]\to\Delta.A[f].$$ Therefore $\j_{c[f^\lilbullet]}$ is well defined and it is a section of $P_{C[f^{\lilbullet\lilbullet\lilbullet}]}$. As $\j_c[f^{\lilbullet\lilbullet\lilbullet}]$ and $\j_{c[f^\lilbullet]}$ are both sections of $P_{C[f^{\lilbullet\lilbullet\lilbullet}]}$, the equality: $$\j_c[f^{\lilbullet\lilbullet\lilbullet}]=\j_{c[f^\lilbullet]}$$ is meaningful. From now on, let us assume that it is satisfied.

\medskip

\noindent\textbf{Stability of $\boldsymbol{\h}$.} Finally, we observe that by definition $\h_{c[f^\lilbullet]}$ is a section of the display map: $${\Delta.A[f].\id_{C[f^{\lilbullet\lilbullet\lilbullet}][\refl_{A[f]}]}[\j_{c[f^\lilbullet]}[\refl_{A[f]}];c[f^\lilbullet]]}\to\Delta.A[f].$$ Analogously, by definition $\h_c$ is a section of the display map: $$\Gamma.A.\id_{C[\refl_A]}[\j_c[\refl_A];c]\to\Gamma.A$$ hence $\h_c[f^\lilbullet]$ is a section of the display map: $${\Delta.A[f].\id_{C[\refl_A]}[\j_c[\refl_A];c][f^\lilbullet]}\to\Delta.A[f].$$ Now, we observe that: 
\[\begin{tikzcd}[sep=scriptsize]
	{\Delta.A[f].\id_{C[f^{\lilbullet\lilbullet\lilbullet}][\refl_{A[f]}]}[\j_{c[f^\lilbullet]}[\refl_{A[f]}];c[f^\lilbullet]]} \\
	{\Delta.A[f].\id_{C[\refl_A][f^\lilbullet]}[\j_{c}[f^{\lilbullet\lilbullet\lilbullet}][\refl_{A[f]}];c[f^\lilbullet]]} \\
	{\Delta.A[f].\id_{C[\refl_A][f^\lilbullet]}[\j_{c}[\refl_A][f^\lilbullet];c[f^\lilbullet]]} \\
	{\Delta.A[f].\id_{C[\refl_A]}[f^\lilbullet.C[\refl_A].C[\refl_A]^\liltriangle][\j_{c}[\refl_A][f^\lilbullet];c[f^\lilbullet]]} \\
	{\Delta.A[f].\id_{C[\refl_A]}[\j_c[\refl_A];c][f^\lilbullet]}
	\arrow[equals, from=1-1, to=2-1]
	\arrow[equals, from=2-1, to=3-1]
	\arrow[equals, from=3-1, to=4-1]
	\arrow["{\text{\,Diagram (\ref{reindexingprimareindexingdopo})}}", equals, from=4-1, to=5-1]
\end{tikzcd}\] where the last equality follows from an instance of the commutative diagram (\ref{reindexingprimareindexingdopo}). We conclude that $\h_c[f^\lilbullet]$ and $\h_{c[f^\lilbullet]}$ are sections of the same display map, hence the equality: $$\h_c[f^\lilbullet]=\h_{c[f^\lilbullet]}$$ is meaningful.

\medskip

\noindent\textbf{Type-checking of the stability conditions of Definition \ref{def:Semantics of axiomatic Sigma-types}}

Let us assume that $(\mathbf{C},\mathcal{D})$ is endowed with axiomatic $=$-types and with the choice functions $\sigmad$, $\pairing$, $\splitt$, $\sigma$ as in \textit{Form Rule}, \textit{Intro Rule}, \textit{Elim Rule}, \textit{Comp Axiom} of Definition \ref{def:Semantics of axiomatic Sigma-types}. Let $f$ be an arrow $\Delta\to\Gamma$, let $P_A$ be a display map of codomain $\Gamma$, and let $P_B$ be a display map of codomain $\Gamma.A$.

\medskip

\noindent\textbf{Stability of $\boldsymbol{\sigmad}$.} The codomain of both $P_{\sigmad_{A[f]}^{B[f^\lilbullet]}}$ and $P_{\sigmad_A^B[f]}$ is by definition $\Delta$. Hence the equality: $$\sigmad_A^B[f]=\sigmad_{A[f]}^{B[f^\lilbullet]}$$ is meaningful. From now on, let us assume that it is satisfied.

\medskip

\noindent\textbf{Stability of $\boldsymbol{\pairing}$.} By definition $\pairing_{A[f]}^{B[f^\lilbullet]}$ is a section of the display map: $$\Delta.A[f].B[f^\lilbullet].\sigmad_{A[f]}^{B[f^\lilbullet]}[P_{A[f]}P_{B[f^\lilbullet]}]\to \Delta.A[f].B[f^\lilbullet].$$ We observe that $\pairing_A^B[f^{\lilbullet\lilbullet}]$ is a section of the display map: $$\Delta.A[f].B[f^\lilbullet].\sigmad_A^B[P_AP_B][f^{\lilbullet\lilbullet}]\to \Delta.A[f].B[f^\lilbullet]$$ and that additionally $\sigmad_A^B[P_AP_B][f^{\lilbullet\lilbullet}]=\sigmad_A^B[f][P_{A[f]}P_{B[f^\lilbullet]}]=\sigmad_{A[f]}^{B[f^\lilbullet]}[P_{A[f]}P_{B[f^\lilbullet]}].$ Therefore $\pairing_A^B[f^{\lilbullet\lilbullet}]$ is itself a section of: $$\Delta.A[f].B[f^\lilbullet].\sigmad_{A[f]}^{B[f^\lilbullet]}[P_{A[f]}P_{B[f^\lilbullet]}]\to \Delta.A[f].B[f^\lilbullet]$$ hence the equality: $$\pairing_A^B[f^{\lilbullet\lilbullet}]=\pairing_{A[f]}^{B[f^\lilbullet]}$$ is meaningful. From now on, let us assume that it is satisfied. Observe that this implies that the square:
\[\begin{tikzcd}[column sep=small,row sep=scriptsize]
	{\Delta.A[f].B[f.A]} && {\Gamma.A.B} \\
	\\
	{\Delta.\sigmad_{A[f]}^{B[f^\lilbullet]}} && {\Gamma.\sigmad_A^B}
	\arrow["{f^{\lilbullet\lilbullet}}", from=1-1, to=1-3]
	\arrow["{\pair_{A[f]}^{B[f^\lilbullet]}}"{description}, from=1-1, to=3-1]
	\arrow["{\pair_A^B}"{description}, from=1-3, to=3-3]
	\arrow["{f^{\lilbullet}}", from=3-1, to=3-3]
\end{tikzcd}\]
commutes.

\medskip

\noindent\textbf{Stability of $\boldsymbol{\splitt}$.} Let $P_C$ be a display map of codomain $\Gamma.\sigmad_A^B$ and let $c$ be a section of $\Gamma.A.B.C[\pair_A^B]\to\Gamma.A.B$. As $\splitt_c$ is by definition a section of $P_C$, then $\splitt_c[f^{\lilbullet}]$ is a section of $P_{C[f^{\lilbullet}]}$. Now, we observe that $P_{C[f^{\lilbullet}]}$ is a display map of codomain: $$\Delta.\sigmad_A^B[f]=\Delta.\sigmad_{A[f]}^{B[f^\lilbullet]}$$ and that the arrow $c[f^{\lilbullet\lilbullet}]$ is a section of the display map $\Delta.A[f].B[f^\lilbullet].C[\pair_A^B][f^{\lilbullet\lilbullet}]\to\Delta.A[f].B[f^\lilbullet]$ i.e. of: $$\Delta.A[f].B[f^\lilbullet].C[f^{\lilbullet}][\pair_{A[f]}^{B[f^\lilbullet]}]\to\Delta.A[f].B[f^\lilbullet].$$ Therefore $\splitt_{c[f^{\lilbullet\lilbullet}]}$ is well defined and it is a section of $P_{C[f^{\lilbullet}]}$. As $\splitt_c[f^{\lilbullet}]$ and $\splitt_{c[f^{\lilbullet\lilbullet}]}$ are both sections of $P_{C[f^{\lilbullet}]}$, the equality: $$\splitt_c[f^{\lilbullet}]=\splitt_{c[f^{\lilbullet\lilbullet}]}$$ is meaningful. From now on, let us assume that it is satisfied.

\medskip

\noindent\textbf{Stability of $\boldsymbol{\sigma}$.} Finally, we observe that by definition $\sigma_{c[f^{\lilbullet\lilbullet}]}$ is a section of the display map: $${\Delta.A[f].B[f^\lilbullet].\id_{C[f^{\lilbullet}][\pair_{A[f]}^{B[f^\lilbullet]}]}[\splitt_{c[f^{\lilbullet\lilbullet}]}[\pair_{A[f]}^{B[f^\lilbullet]}];c[f^{\lilbullet\lilbullet}]]}\to\Delta.A[f].B[f^\lilbullet].$$ Analogously, by definition $\sigma_c$ is a section of the display map: $$\Gamma.A.B.\id_{C[\pair_A^B]}[\splitt_c[\pair_A^B];c]\to\Gamma.A.B$$ hence $\sigma_c[f^{\lilbullet\lilbullet}]$ is a section of the display map: $${\Delta.A[f].B[f^\lilbullet].\id_{C[\pair_A^B]}[\splitt_c[\pair_A^B];c][f^{\lilbullet\lilbullet}]}\to\Delta.A[f].B[f^\lilbullet].$$ Now, we observe that: 
\[\begin{tikzcd}[cramped,sep=scriptsize]
	{\Delta.A[f].B[f^\lilbullet].\id_{C[f^\lilbullet][\pair_{A[f]}^{B[f^\lilbullet]}]}[\splitt_{c[f^{\lilbullet\lilbullet}]}[\pair_{A[f]}^{B[f^\lilbullet]}];c[f^{\lilbullet\lilbullet}]]} \\
	{\Delta.A[f].B[f^\lilbullet].\id_{C[\pair_{A}^{B}][f^{\lilbullet\lilbullet}]}[\splitt_{c[f^{\lilbullet\lilbullet}]}[\pair_{A[f]}^{B[f^\lilbullet]}];c[f^{\lilbullet\lilbullet}]]} \\
	{\Delta.A[f].B[f^\lilbullet].\id_{C[\pair_{A}^{B}][f^{\lilbullet\lilbullet}]}[\splitt_{c}[f^\lilbullet][\pair_{A[f]}^{B[f^\lilbullet]}];c[f^{\lilbullet\lilbullet}]]} \\
	{\Delta.A[f].B[f^\lilbullet].\id_{C[\pair_{A}^{B}][f^{\lilbullet\lilbullet}]}[\splitt_{c}[\pair_A^B][f^{\lilbullet\lilbullet}];c[f^{\lilbullet\lilbullet}]]} \\
	{\Delta.A[f].B[f^\lilbullet].\id_{C[\pair_{A}^{B}]}[f^{\lilbullet\lilbullet}.C[\pair_{A}^{B}].C[\pair_{A}^{B}]^\liltriangle][\splitt_{c}[\pair_A^B][f^{\lilbullet\lilbullet}];c[f^{\lilbullet\lilbullet}]]} \\
	{\Delta.A[f].B[f^\lilbullet].\id_{C[\pair_{A}^{B}]}[\splitt_{c}[\pair_A^B];c][f^{\lilbullet\lilbullet}]}
	\arrow[equals, from=1-1, to=2-1]
	\arrow[equals, from=2-1, to=3-1]
	\arrow[equals, from=3-1, to=4-1]
	\arrow[equals, from=4-1, to=5-1]
	\arrow["{\text{\,Diagram (\ref{reindexingprimareindexingdopo})}}", equals, from=5-1, to=6-1]
\end{tikzcd}\] where the last equality follows from an instance of the commutative diagram (\ref{reindexingprimareindexingdopo}). We conclude that $\sigma_c[f^{\lilbullet\lilbullet}]$ and $\sigma_{c[f^{\lilbullet\lilbullet}]}$ are sections of the same display map, hence the equality: $$\sigma_c[f^{\lilbullet\lilbullet}]=\sigma_{c[f^{\lilbullet\lilbullet}]}$$ is meaningful.

\section{Stability conditions in display map 2-categories}\label{appendixIII}

This section contains the proof of Proposition \ref{stability of id in 2d model}. The proof of Proposition \ref{stability of sigma in 2d model} is completely analogous, and uses the \textit{stability conditions} in Definition \ref{def:Categorical semantics of axiomatic Sigma-types} in place of the ones in Definition \ref{def:Categorical semantics of axiomatic =-types}. We also provide the proof of \Cref{stability of pi in 2d model} for $\pid$, $\ev$, $\abst{}$, and $\beta$. The proof of \Cref{stability of pi in 2d model} for $\funext$, $\beta^\Pi$, and $\eta^\Pi$ follows the same argument together with \Cref{happytheorem}. Additionally, the proof of \Cref{stability of nat in 2d model} is analogous to the one of \Cref{stability of id in 2d model} and \Cref{stability of sigma in 2d model}.

\medskip

Let $(\boldsymbol{\mathcal{C}},\mathcal{D})$ be a display map 2-category endowed with axiomatic $=$-types. Let $f$ be an arrow $\Delta\to\Gamma$ and let $P_A$ be a display map of codomain $\Gamma$.

\medskip

\noindent\textbf{Stability of $\boldsymbol{\id}$.} The first stability condition $\id_A[f^{\lilbullet\lilbullet}]=\id_{A[f]}$ holds by definition.

\medskip

\noindent\textbf{Stability of $\boldsymbol{\reflex}$.} Now, the following diagram: 
\[\begin{tikzcd}[column sep=small,row sep=scriptsize]
	{\Delta.A[f].A[f]^\liltriangle.\id_{A[f]}} && {\Gamma.A.A^\liltriangle.\id_A} \\
	\\
	{\Delta.A[f]} && {\Gamma.A}
	\arrow["{f^{\lilbullet\lilbullet\lilbullet}}"{pos=0.3}, from=1-1, to=1-3]
	\arrow[""{name=0, anchor=center, inner sep=0}, shift right=5, from=1-1, to=3-1]
	\arrow[""{name=1, anchor=center, inner sep=0}, shift left=5, from=1-1, to=3-1]
	\arrow[""{name=2, anchor=center, inner sep=0}, shift right=3, from=1-3, to=3-3]
	\arrow[""{name=3, anchor=center, inner sep=0}, shift left=3, from=1-3, to=3-3]
	\arrow["{f^{\lilbullet}}"', from=3-1, to=3-3]
	\arrow["{\alpha_A[f]}", shorten <=4pt, shorten >=4pt, Rightarrow, from=0, to=1]
	\arrow["{\alpha_A}", shorten <=2pt, shorten >=2pt, Rightarrow, from=2, to=3]
\end{tikzcd}\] commutes and $\alpha_A[f]=\alpha_{A[f]}$, by definition. Therefore: \begin{center}
$\alpha_A*(f^{\lilbullet\lilbullet\lilbullet}\refl_{A[f]})=(\alpha_A*f^{\lilbullet\lilbullet\lilbullet})*\refl_{A[f]}=(f^\lilbullet*\alpha_A[f])*\refl_{A[f]}=(f^\lilbullet*\alpha_{A[f]})*\refl_{A[f]}=f^\lilbullet*(\alpha_{A[f]}*\refl_{A[f]})=1_{f^\lilbullet}$
\end{center} and: \begin{center}
    $\alpha_A*(\refl_Af^\lilbullet)=(\alpha_A*\refl_A)*f^\lilbullet=1_{f^\lilbullet}$
\end{center} hence the diagram: 
\[\begin{tikzcd}[column sep=small,row sep=scriptsize]
	{\Delta.A[f]} && {\Gamma.A} \\
	\\
	{\Delta.A[f].A[f]^\liltriangle.\id_{A[f]}} && {\Gamma.A.A^\liltriangle.\id_A}
	\arrow["{f^\lilbullet}", from=1-1, to=1-3]
	\arrow["{\refl_{A[f]}}"', from=1-1, to=3-1]
	\arrow["{\refl_A}", from=1-3, to=3-3]
	\arrow["{f^{\lilbullet\lilbullet\lilbullet}}"{pos=0.3}, from=3-1, to=3-3]
\end{tikzcd}\] commutes being $\alpha_A$ an arrow object. As usual, the commutativity of this diagram is equivalent to the stability condition $\reflex_A[f^\lilbullet]=\reflex_{A[f]}$.

\medskip

\noindent\textbf{Stability of $\boldsymbol{\j}$.} Let $P_C$ be a display map of codomain $\Gamma.A.A^\liltriangle.\id_A$ and let $c$ be a section of $P_{C[\pair_A^B]}$. Let us observe that the diagram: 
\[\begin{tikzcd}[cramped,sep=scriptsize]
	{\Delta.A[f].A[f]^\liltriangle.\id_{A[f]}} && {\Gamma.A.A^\liltriangle.\id_{A}} \\
	{\Delta.A[f].A[f]^\liltriangle} && {\Gamma.A.A^\liltriangle} \\
	{\Delta.A[f]} && {\Gamma.A} \\
	{\Delta.A[f].C[f^{\lilbullet\lilbullet\lilbullet}][\refl_{A[f]}]} && {\Gamma.A.C[\refl_A^\lilbullet]} \\
	{\Delta.A[f].A[f]^\liltriangle.\id_{A[f]}.C[f^{\lilbullet\lilbullet\lilbullet}]} && {\Gamma.A.A^\liltriangle.\id_{A}.C}
	\arrow["{f^{\lilbullet\lilbullet\lilbullet}}", from=1-1, to=1-3]
	\arrow["{P_{\id_{A[f]}}}"', from=1-1, to=2-1]
	\arrow["{P_{\id_A}}", from=1-3, to=2-3]
	\arrow["{f^{\lilbullet\lilbullet}}"{description}, from=2-1, to=2-3]
	\arrow["{P_{A[f]}^{\lilbullet}}"', from=2-1, to=3-1]
	\arrow["{P_A^{\lilbullet}}", from=2-3, to=3-3]
	\arrow["{f^{\lilbullet}}"{description}, from=3-1, to=3-3]
	\arrow["{c[f^{\lilbullet}]}"', from=3-1, to=4-1]
	\arrow["c", from=3-3, to=4-3]
	\arrow["{f^{\lilbullet}.C[\refl_A^\lilbullet]}"{description}, from=4-1, to=4-3]
	\arrow["{\refl_{A[f]}.C[f^{\lilbullet\lilbullet\lilbullet}]}"', from=4-1, to=5-1]
	\arrow["{\refl_A.C}", from=4-3, to=5-3]
	\arrow["{f^{\lilbullet\lilbullet\lilbullet}.C}"', from=5-1, to=5-3]
\end{tikzcd}\] commutes and that the two columns yield: $$P_{\id_{A[f]}}\text{ 
\quad and 
\quad }P_{\id_{A}}$$ when post-composed by $P_{\id_{A[f]}}P_{C[f^{\lilbullet\lilbullet\lilbullet}]}$ and by $P_{\id_{A}}P_C$ respectively.

Moreover, let us observe that: $$\begin{aligned}
    f^{\lilbullet}*(P_{A[f]^\liltriangle}P_{\id_{A[f]}}*\varphi_A[f])&=P_{A^\liltriangle}P_{\id_A}*(f^{\lilbullet\lilbullet\lilbullet} * \varphi_A[f])\\
    &=(P_{A^\liltriangle}P_{\id_A}*\varphi_A)*f^{\lilbullet\lilbullet\lilbullet}\\
    &=\alpha_A*f^{\lilbullet\lilbullet\lilbullet}\\
    \textnormal{ \; }&\textnormal{ \; }\\
    P_{A[f]}*(P_{A[f]^\liltriangle}P_{\id_{A[f]}}*\varphi_A[f])&=P_{A[f]}P_{A[f]^\liltriangle}P_{\id_{A[f]}}*\varphi_A[f]\\
    &=1_{P_{A[f]}P_{A[f]^\liltriangle}P_{\id_{A[f]}}}
\end{aligned}$$ which implies that: $$P_{A[f]^\liltriangle}P_{\id_{A[f]}}*\varphi_A[f] = \alpha_{A}[f]=\alpha_{A[f]}.$$ Analogously: $$\begin{aligned}
    f^{\lilbullet}*(P_{A[f]}^\lilbullet P_{\id_{A[f]}}*\varphi_A[f])&=P_A^\lilbullet P_{\id_A}*(f^{\lilbullet\lilbullet\lilbullet} * \varphi_A[f])\\
    &=(P_A^\lilbullet P_{\id_A}*\varphi_A)*f^{\lilbullet\lilbullet\lilbullet}\\
    &=1_{P_A^\lilbullet P_{\id_A}}*f^{\lilbullet\lilbullet\lilbullet}\\
    \textnormal{ \; }&\textnormal{ \; }\\
    P_{A[f]}*(P_{A[f]}^\lilbullet P_{\id_{A[f]}}*\varphi_A[f])&=P_{A[f]}P_{A[f]}^\lilbullet P_{\id_{A[f]}}*\varphi_A[f]\\
    &=1_{P_{A[f]}P_{A[f]}^\lilbullet P_{\id_{A[f]}}}
\end{aligned}$$ which implies that: $$P_{A[f]}^\lilbullet P_{\id_{A[f]}}*\varphi_A[f] = 1_{P_{A[f]}^\lilbullet P_{\id_{A[f]}}}.$$ As: $$\begin{aligned}
    P_{A[f]^\liltriangle}P_{\id_{A[f]}}*\varphi_A[f] &= \alpha_{A[f]}\\
    P_{A[f]}^\lilbullet P_{\id_{A[f]}}*\varphi_A[f] &= 1_{P_{A[f]}^\lilbullet P_{\id_{A[f]}}}
\end{aligned}$$ it must be that $\varphi_A[f]=\varphi_{A[f]}$ by the universal property of $\alpha_{A[f]}$.

Therefore: $$\begin{aligned}
    \j_c[f^{\lilbullet\lilbullet\lilbullet}]&=\t_{(\refl_A.C) c P_A^\lilbullet P_{\id_A}}^{\varphi_A}[f^{\lilbullet\lilbullet\lilbullet}]\\
    &=\t_{(\refl_{A[f]}.C[f^{\lilbullet\lilbullet\lilbullet}])c[f^{\lilbullet}]P_{A[f]}^\lilbullet P_{\id_{A[f]}}}^{\varphi_A[f]}\\
    &=\t_{(\refl_{A[f]}.C[f^{\lilbullet\lilbullet\lilbullet}])c[f^{\lilbullet}]P_{A[f]}^\lilbullet P_{\id_{A[f]}}}^{\varphi_{A[f]}}\\
    &=\j_{c[f^\lilbullet]}.
\end{aligned}$$

\medskip

\noindent\textbf{Stability of $\boldsymbol{\h}$.} Analogously: $$\tau_{(\refl_A.C) c P_A^\lilbullet P_{\id_A}}^{\varphi_A}[f^{\lilbullet\lilbullet\lilbullet}]=\tau_{(\refl_{A[f]}.C[f^{\lilbullet\lilbullet\lilbullet}])c[f^{\lilbullet}]P_{A[f]}^\lilbullet P_{\id_{A[f]}}}^{\varphi_{A[f]}}$$ hence: $$\begin{aligned}
    h_c[f^{\lilbullet}]&=\tau_{(\refl_A.C) c P_A^\lilbullet P_{\id_A}}^{\varphi_A}[\refl_A][f^\lilbullet]\\
    &=\tau_{(\refl_A.C) c P_A^\lilbullet P_{\id_A}}^{\varphi_A}[f^{\lilbullet\lilbullet\lilbullet}][\refl_{A[f]}]\\
    &=\tau_{(\refl_{A[f]}.C[f^{\lilbullet\lilbullet\lilbullet}])c[f^{\lilbullet}]P_{A[f]}^\lilbullet P_{\id_{A[f]}}}^{\varphi_{A[f]}}[\refl_{A[f]}]\\
    &=h_{c[f^{\lilbullet}]}.
\end{aligned}$$
In particular, since $\alpha_{C[\refl_A]}[f^\lilbullet] * \tilde h_c[f^\lilbullet]$ is a section of $P_{C[f^{\lilbullet\lilbullet\lilbullet}][\refl_{A[f]}]}$ such that: $$f^\lilbullet.C[\refl_A]*(\alpha_{C[\refl_A]}[f^\lilbullet] * \tilde h_c[f^\lilbullet])=(\alpha_{C[\refl_A]} * \tilde h_c)*f^\lilbullet=h_c*f^\lilbullet$$ it must be that $\alpha_{C[\refl_A]}[f^\lilbullet] * \tilde h_c[f^\lilbullet] = h_c[f^\lilbullet]$, which implies that: $$\tilde h_c[f^\lilbullet]=\tilde h_{c[f^\lilbullet]}$$ as $h_c[f^\lilbullet]=h_{c[f^\lilbullet]}$ and $\alpha_{C[\refl_A]}[f^\lilbullet]=\alpha_{C[\refl_A][f^\lilbullet]}=\alpha_{C[f^{\lilbullet\lilbullet\lilbullet}][\refl_{A[f]}]}$.

Finally, referring to the diagram:
\[\begin{tikzcd}[cramped,column sep=tiny,row sep=scriptsize]
	& {} \\
	\\
	{} & {\id_{C[\refl_A]}[\j_c[\refl_A];c]} & {C[\refl_A].C[\refl_A]^\liltriangle.\id_{C[\refl_A]}} \\
	\\
	{\id_{C'}[\j_{c[f^\lilbullet]}[\refl_{A[f]}];c[f^\lilbullet]]} & {C'.C'^\liltriangle.\id_{C'}} \\
	& {} & {C[\refl_A].C[\refl_A]^\liltriangle} \\
	\\
	{} & {C'.C'^\liltriangle}
	\arrow["{\h_c}"{description}, dashed, from=1-2, to=3-2]
	\arrow["{\tilde h_c}"{description, pos=0.2}, from=1-2, to=3-3]
	\arrow[from=3-1, to=1-2]
	\arrow["{\h_{c[f^\lilbullet]}}"{description}, from=3-1, to=5-1]
	\arrow["{\tilde h_{c[f^\lilbullet]}}"{description, pos=0.2}, from=3-1, to=5-2]
	\arrow[dashed, from=3-2, to=3-3]
	\arrow[dashed, from=3-2, to=6-2]
	\arrow[from=3-3, to=6-3]
	\arrow[dashed, from=5-1, to=3-2]
	\arrow[from=5-1, to=5-2]
	\arrow[from=5-1, to=8-1]
	\arrow[from=5-2, to=3-3]
	\arrow[shift left=3, from=5-2, to=8-2]
	\arrow["{\j_c[\refl_A];c}"'{pos=0.7}, dashed, from=6-2, to=6-3]
	\arrow[from=8-1, to=6-2]
	\arrow["{\j_{c[f^\lilbullet]}[\refl_{A[f]}];c[f^\lilbullet]}"', from=8-1, to=8-2]
	\arrow[from=8-2, to=6-3]
\end{tikzcd}\]---where we omitted the context $\Gamma.A$ in the objects of the front side and the context $\Delta.A[f]$ in the ones of the back side, and where $C':=C[f^{\lilbullet\lilbullet\lilbullet}][\refl_{A[f]}]$---a priori every atomic diagram commutes apart from: 
\begin{equation}\label{piazza}
    \begin{tikzcd}[cramped, column sep=small, row sep=scriptsize]
	{\Delta.A[f]} & {\Gamma.A} \\
	\\
	{\Delta.A[f].\id_{C'}[\j_{c[f^\lilbullet]}[\refl_{A[f]}];c[f^\lilbullet]]} & {\Gamma.A.\id_{C[\refl_A]}[\j_c[\refl_A];c]}
	\arrow[from=1-1, to=1-2]
	\arrow["{\h_{c[f^\lilbullet]}}"{description}, from=1-1, to=3-1]
	\arrow["{\h_c}"{description}, dashed, from=1-2, to=3-2]
	\arrow[dashed, from=3-1, to=3-2]
\end{tikzcd}
\end{equation} and we are going to prove that it actually does. In fact, we observe that: $$\begin{aligned}
    P_{\id_{C[\refl_A]}}(f^\lilbullet.\id_{C[\refl_A]}[\j_C[\refl_A];c]\h_{c[f^\lilbullet]})&=f^\lilbullet P_{\id_{C'}[\j_{c[f^\lilbullet]}[\refl_{A[f]}];c[f^\lilbullet]]}h_{c[f^\lilbullet]}\\
    &=f^\lilbullet\\
    &=P_{\id_{C[\refl_A]}}(\h_cf^\lilbullet)
\end{aligned}$$ and: $$\begin{aligned}
    (\j_c[\refl_A];c).\id_{C[\refl_A]}(f^\lilbullet.\id_{C[\refl_A]}[\j_C[\refl_A];c]\h_{c[f^\lilbullet]})&=f^\lilbullet.C[\refl_A].C[\refl_A]^\liltriangle.\id_{C[\refl_A]}\j_{c[f^\lilbullet]}[\refl_{A[f]};c[f^\lilbullet]].\id_{C'}\h_{c[f^\lilbullet]}\\
    &=f^\lilbullet.C[\refl_A].C[\refl_A]^\liltriangle.\id_{C[\refl_A]}\tilde h_{c[f.A]}\\
    &=\tilde h_cf^\lilbullet\\
    &= (\j_c[\refl_A];c).\id_{C[\refl_A]}(\h_cf^\lilbullet).
\end{aligned}$$ hence, by the universal property of $\Gamma.A.\id_{C[\refl_A]}[\j_c[\refl_A];c]$, the diagram (\ref{piazza}) commutes. Therefore: $$\h_c[f^\lilbullet] = \h_{c[f^\lilbullet]}$$ and we are done.

\medskip

Let $(\boldsymbol{\mathcal{C}},\mathcal{D})$ be a display map 2-category endowed with axiomatic $=$-types and axiomatic $\Pi$-types \& axiomatic function extensionality. Let $f$ be an arrow $\Delta\to\Gamma$, let $P_A$ be a display map of codomain $\Gamma$ and let $P_B$ be a display map of codomain $\Gamma.A$.

\medskip

\noindent\textbf{Stability of $\boldsymbol{\pid}$.} The first stability condition $\pid_A^B[f]=\pid_{A[f]}^{B[f^\lilbullet]}$ holds by definition.

\medskip

\noindent\textbf{Stability of $\boldsymbol{\ev}$.} We verified the second stability condition $\ev_A^B[f^{\lilbullet\lilbullet}]=\ev_{A[f]}^{B[f^\lilbullet]}$ in (\ref{evstabilesottoreindexing}) in the paragraph \textbf{Stability of $\boldsymbol{\app}$ and $\boldsymbol{\lambda}$ under re-indexing.}

\medskip

\noindent\textbf{Stability of $\boldsymbol{\abst{}}$.} Let $\upsilon$ be a section of $P_B$. Then: \begin{center}$\abst{\upsilon}=\lambda_A^B\upsilon$ \quad and \quad $\abst{\upsilon[f^\lilbullet]}=\lambda_{A[f]}^{B[f^\lilbullet]}(\upsilon[f^\lilbullet])$\end{center} by definition. Hence: $$\abst{\upsilon}[f]=(\lambda_A^B\upsilon)[f]\overset{(\ref{appandlambdastableunderpullback})}{=}\lambda_{A[f]}^{B[f^\lilbullet]}(\upsilon[f^\lilbullet])=\abst{\upsilon[f^\lilbullet]}$$ and we are done.

\medskip

\noindent\textbf{Stability of $\boldsymbol{\beta}$.} By definition: \begin{center}
    $\beta_\upsilon=\{\beta^{A,B}_\upsilon\}$ \quad and \quad $\beta_{\upsilon[f^\lilbullet]}=\{\beta^{A[f],B[f^\lilbullet]}_{\upsilon[f^\lilbullet]}\}$
\end{center} hence: $$\begin{aligned}
    \beta_\upsilon[f^\lilbullet]&=\{\beta^{A,B}_\upsilon\}[f^\lilbullet]\\
    &=\{\beta^{A,B}_\upsilon[f^\lilbullet]\}\\
    &=\{1;\beta^{A,B}_\upsilon f^\lilbullet\}\\
    &=\{1;\beta^{A,B}_{\upsilon f^\lilbullet}\}\\
    &=\{\beta^{A[f],B[f^\lilbullet]}_{1;\upsilon f^\lilbullet}\}\\
    &=\{\beta^{A[f],B[f^\lilbullet]}_{\upsilon [f^\lilbullet]}\}\\
    &=\beta_{\upsilon[f^\lilbullet]}
\end{aligned}$$ where the fourth and the fifth equalities hold by \Cref{def:Categorical semantics of axiomatic Pi-types}. We are done.

\section{Further details on the content of Section \ref{sec:Revisiting the groupoid model}}\label{sec:further details on the weakened groupoid model} We provide further details on the structure that the display map 2-category $(\grpd,\mathcal{D})$, defined in Section \ref{sec:Revisiting the groupoid model}, is endowed with.

We stipulated that a display map $P_A:\Gamma.A\to\Gamma$ of the class $\mathcal{D}$ is obtained by applying the Grothendieck construction to a pseudofunctor $(A,\phi^A,\psi^A):\Gamma\to\grpd$, where $\Gamma$ is a groupoid, and $\phi^A$ and $\psi^A$ are the \textit{coherent} families of natural isomorphisms $\phi^A_{{p},{q}}: A_qA_p\rightarrow A_{qp}$ and $\psi_\gamma^A:A_{1_{\gamma}}\rightarrow1_{A_\gamma}$---where $\gamma$ is any object of $\Gamma$ and $p$ and $q$ any composable arrows of $\Gamma$---making $A$ into a pseudofunctor. We recall that the coherence laws are the equalities: \begin{center}
    $\phi^A_{p,rq}(\phi^A_{q,r}*A_p)=\phi^A_{qp,r}(A_r*\phi^A_{p,q})$\\
    \text{ }\\
    $A_p*\psi_\gamma^A=\phi^A_{1_\gamma,p}$\quad\quad$\psi_{\gamma'}^A*A_p=\phi^A_{p,1_{\gamma'}}$
\end{center} for every choice of arrows $p:\gamma\to\gamma'$, $q:\gamma'\to\gamma''$, and $r:\gamma''\to\gamma'''$ in $\Gamma$.

This means that $\Gamma.A$ is the groupoid whose objects are the pairs $(\gamma,x)$ where $\gamma$ is an object of $\Gamma$ and $x$ is an object of $A_\gamma$, and the arrows $(\gamma,x)$ to $(\gamma',x')$ are the pairs $({p_1},{p_2})$ where ${p_1}$ is an arrow $\gamma\to\gamma'$ and ${p_2}$ is an arrow $A_{p_1} x\to x'$.

In detail, the composition of two arrows: $$(\gamma,x)\xrightarrow{({p_1},{p_2})}(\gamma',x')\xrightarrow{({q_1},{q_2})}(\gamma'',x'')$$ is: $$({q_1}{p_1}, A_{{q_1}{p_1}}x\xrightarrow{(\phi^A_{{p_1},{q_1}})^{-1}x}A_{{q_1}}A_{p_1} x \xrightarrow{A_{{q_1}}{p_2}}A_{{q_1}}x'\xrightarrow{{q_2}}x'')$$ the identity arrow of the object $(\gamma,x)$ is $(1_\gamma,\psi_\gamma^Ax)$, and the inverse of $({p_1},{p_2})$ exists and is the pair: \begin{center}
    $({p_1}^{-1},$\\
    \text{ }\\
    $A_{{p_1}^{-1}}x' \xrightarrow{A_{{p_1}^{-1}}{p_2}^{-1}}A_{{p_1}^{-1}}A_{p_1} x \xrightarrow{\phi^A_{{p_1},{p_1}^{-1}}x} A_{{p_1}^{-1}{p_1}}x =$\\
    \text{ }\\
    $= A_{1_\gamma}x \xrightarrow{\psi_\gamma^A x}x)$
\end{center} and associativity, unitalities, and invertibilities follow from the coherence laws of the families $\phi^A$ and $\psi^A$. Hence $\Gamma.A$ is a groupoid. The functor $\Gamma.A\to\Gamma$ maps $(\gamma,x)$ to $\gamma$ and $({p_1},{p_2})$ to ${p_1}$.

Regarding the re-indexing structure on display maps, we observed that, for every 1-cell $f:\Delta\to\Gamma$ in $\grpd$ and every pseudofunctor $(A,\phi^A,\psi^A):\Gamma\to\grpd$, the triple: $$(Af,\phi^A_{f-,f-},\psi^A_{f-})$$ defines a pseudofunctor $\Delta\to\grpd$. We stipulated that the associated Grothendieck construction: $$P_{Af}:\Delta.Af\to\Delta$$ denoted as $P_{A[f]}:\Delta.A[f]\to\Delta$, is the re-indexing of $P_A$ along $f$. The mappings:\begin{center} $(\;\delta : \Delta\;,\;x : (Af)_\delta\;)\mapsto (\;f\delta\;,\;x\;)$ and $(\;\delta\xrightarrow{p_1}\delta'\;,\;(Af)_{p_1}x\xrightarrow{p_2}x')\mapsto (\;f p_1\;,\;p_2\;)$\end{center} define a functor $\Delta.A[f] \xrightarrow{f.A} \Gamma.A$ such that the square: 
\[\begin{tikzcd}[column sep=scriptsize]
	{\Delta.A[f]} & {\Gamma.A} \\
	\Delta & \Gamma
	\arrow["{f.A}", from=1-1, to=1-2]
	\arrow[from=1-1, to=2-1]
	\arrow[from=1-2, to=2-2]
	\arrow["f", from=2-1, to=2-2]
\end{tikzcd}\] is verified to be a pullback and a 2-pullback. In fact, if we are given functors $g_1:\Omega\to\Delta$ and $(g_2,g_3):\Omega\to\Gamma.A$ in $\grpd$ such that $f g_1=g_2$ then the unique factorisation functor $\Omega \to \Delta.A[f]$ will be given by the pairing $(g_1,g_3)$. Similarly one verifies the 2-dimensional factorisation property.

Moreover, as the triples: \begin{center} $(\;A(fg)\;,\;\phi^A_{fg_-,fg_-}\;,\;\psi^A_{fg_-}\;)$ and $(\;(Af)g\;,\;\phi'^{Af}_{g_-,g_-}\;,\;\psi'^{Af}_{g_-}\;)$
\end{center}---where $\phi'^{Af} =\phi^A_{f_-,f_-}$ and $\psi'^{Af}=\psi^A_{f_-}$---coincide whenever $f$ and $g$ are composable 1-cells, the display maps $P_{A[fg]}$ and $P_{A[f][g]}$ coincide. Similarly, one verifies that $P_{A[1_\Gamma]}$ coincides with $P_A$, and that the cloven isofibration structure on display maps that we provided in Section \ref{sec:Revisiting the groupoid model} is compatible with this re-indexing choice in the sense of Definition \ref{split display map 2-category}.

Finally, we observed that the mappings: $$\begin{aligned}
    (\gamma,x,y)&\mapsto A_\gamma(x,y)\\(p_1,p_2,p_3)&\mapsto(p\mapsto (x'\xrightarrow{p_2^{-1}}A_{p_1}x\xrightarrow{A_{p_1}p}A_{p_1}y\xrightarrow{p_3}y'))
\end{aligned}$$---where the homset $A_\gamma(x,y)$ is endowed with the trivial groupoid structure---determine a pseudofunctor---in fact a strict functor---$\Gamma.A.A[P_A]\to\grpd$ that we indicate as $\id_A$. Hence we obtained a display map $\Gamma.A.A[P_A].\id_A\to\Gamma.A.A[P_A]$. Additionally, we stated that the 2-cell: 
\[\begin{tikzcd}[sep=scriptsize]
	{\Gamma.A.A^\liltriangle.\id_A} &&& {\Gamma.A} \\
	&&& \Gamma
	\arrow[""{name=0, anchor=center, inner sep=0}, "{P_A^\lilbullet P_{\id_A}}"', shift right=2, from=1-1, to=1-4]
	\arrow[""{name=1, anchor=center, inner sep=0}, "{P_{A^\liltriangle}P_{\id_A}}", shift left=2, from=1-1, to=1-4]
	\arrow["{P_AP_{A^\liltriangle}P_{\id_A}}"', curve={height=12pt}, from=1-1, to=2-4]
	\arrow["{P_A}", from=1-4, to=2-4]
	\arrow["{\;\alpha_A}", shorten <=1pt, shorten >=1pt, Rightarrow, from=1, to=0]
\end{tikzcd}\] whose $(\gamma,x,y,p)$-component is the arrow: $$(\;1_\gamma\;,\;A_{1_\gamma}x\xrightarrow{\psi^A}x\xrightarrow{p}y\;):(\gamma,x)\to(\gamma,y)$$ constitutes an arrow object, which can be verified to be compatible with the re-indexing choice in the sense of Definition \ref{def:Categorical semantics of axiomatic =-types}.

We also mentioned that, following the construction at paragraphs \textbf{Elim Rule} and \textbf{Comp Axiom for $\boldsymbol{=}$-types}, one can reconstruct the choice functions $\refl$, $\phi$, and $\j$ in order to observe that $\j_c[\refl_A]$ does not need to coincide with $c$, if $c$ is a section of $P_{C[\refl_A]}$ for a given display map $P_C$ over $\Gamma.A.A^\liltriangle.\id_A$. In detail, the setions $\j_c[\refl_A]$ and $c$ will be different when $P_C$ is not normal, i.e. when $C$ is not a strict functor to $\grpd$.

For sake of completeness, we mention that the functor: $$\refl_A:\Gamma.A\to\Gamma.A.A^\liltriangle.\id_A$$ acts \textit{on objects} as the mapping: $$(\gamma,x)\mapsto(\gamma,x,x,1_x)$$ and that the $(\gamma,x,y,p)$-component of $\varphi_A$ is the arrow: $$(1_\gamma,p\psi^A_\gamma x,\psi^A_\gamma x,1_{1_y})$$ from $(\gamma,x,y,p)$ to $(\gamma,y,y,1_y)$. This can be used to verify that, whenever the mapping $(\gamma,x)\mapsto (\gamma,x,c(\gamma,x))$ is the action on objects of a section $c$ of $P_{C[\refl_A]}$ for some display map $P_C$ over $\Gamma.A.A[P_A].\id_A$, then the functor: $$\j_c:\Gamma.A.A[P_A].\id_A\to\Gamma.A.A[P_A].\id_A.C$$ acts on objects as the mapping: $$(\gamma,x,y,p)\mapsto(\gamma,x,y,p,C_{\varphi_A^{-1}}c(\gamma,y))$$ i.e. as $$(\gamma,x,y,p)\mapsto(\gamma,x,y,p,C_{(1_\gamma,p^{-1}\psi^A_\gamma y,\psi^A_\gamma y,1_{p})}c(\gamma,y))$$ hence the functor $\j_c[\refl_A]:\Gamma.A\to\Gamma.A.C[\refl_A]$ acts on objects as the mapping: $$(\gamma,x)\mapsto(\gamma,x,C_{1_{(\gamma,x,x,1_x)}}c(\gamma,x)).$$ Therefore, as we said, unless $C$ was chosen to be a strict functor $\Gamma.A.A[P_A].\id_A\to\grpd$, in general $\j_c[\refl_A]$ and $c$ will not coincide.

\newpage

\addcontentsline{toc}{section}{References}
\bibliographystyle{plain}
\bibliography{references}

\end{document}